\documentclass[a4paper,11pt]{article}
\usepackage{amssymb,amsthm,amsmath}
\usepackage{color,graphicx}
\usepackage{mathrsfs}  
\usepackage{pifont}
\usepackage{float}
\usepackage{url}
\usepackage{verbatim}
\usepackage[titletoc,title]{appendix}

\usepackage[T1]{fontenc}
\usepackage{lmodern}

\makeatletter

\@addtoreset{equation}{section}
\makeatother

\newcommand{\R}{{\mathbb{R}}}
\newcommand{\C}{{\mathbb{C}}}
\newcommand{\N}{{\mathbb{N}}}
\newcommand{\Z}{{\mathbb{Z}}}

\DeclareMathOperator{\sech}{sech}
\DeclareMathOperator{\Ai}{Ai}
\DeclareMathOperator{\Real}{Re}

\newcommand{\ZT}[1]{\mathcal{Z}\{#1\}} 
\newcommand{\IZT}[1]{\mathcal{Z}^{-1}\{#1\}} 

\definecolor{mygreen}{rgb}{0.11,0.7,0.22}
\definecolor{mred}{rgb}{1,0.,0.}
\definecolor{mcyan}{rgb}{1.,0.5,0.}
\definecolor{mgreen}{rgb}{0.,0.5,0.}
\definecolor{mblue}{rgb}{0.,0.,1.}

\newtheorem{theorem}{Theorem}[section]

\newtheorem{remark}{Remark}[section]

\title{Discrete Artificial Boundary Conditions\\ for the Korteweg-de Vries Equation}
\author{\bf{C. Besse} \\
        Institut de Math{\'e}matiques de Toulouse UMR5219,\\
        Universit\'e de Toulouse; CNRS,\\
        UPS IMT, F-31062 Toulouse Cedex 9, France.\\
	\texttt{Christophe.Besse@math.univ-toulouse.fr} \\\\\
{\bf M. Ehrhardt} \\
	Bergische Universit\"at Wuppertal,\\ 
  	Fachbereich Mathematik und Naturwissenschaften,\\
  	Angewandte Mathematik - Numerische Analysis,\\ 
  	Gau{\ss}strasse 20, 42119 Wuppertal, Germany.\\
        \texttt{ehrhardt@math.uni-wuppertal.de} \\\\
{\bf I. Lacroix-Violet} \\
	Laboratoire Paul Painlev\'e, CNRS UMR 8524, Universit\'e Lille 1,\\ 
  	59655 Villeneuve d'Ascq Cedex, France.\\
	\texttt{Ingrid.Violet@math.univ-lille1.fr}}
\date{\today}
\begin{document}

\maketitle

\begin{abstract}
  We consider the derivation of continuous and fully discrete artificial boundary conditions for linearized Korteweg-de Vries equation. We show that we can obtain them for any constant velocities and any dispersion. The discrete artificial boundary conditions are provided for two different numerical schemes. In both continuous and discrete case, the boundary conditions are non local with respect to time variable. We propose fast evaluations of discrete convolutions. We present various numerical tests which show the effectivness of the artificial boundary conditions.
\end{abstract}


\section{Introduction}
\emph{Korteweg-de Vries (KdV) equations} are typical dispersive nonlinear 
partial differential equations (PDEs).  Zabusky and Kruskal
\cite{ZaKr65} observed that KdV equation owns wave-like solutions
which can retain their initial  forms after collision with another
wave. This led them to name these solitary wave solutions ``solitons''. 

These special solutions were observed and investigated for the first 
time in 1834 by Scott Russell \cite{Eil98,Rus37}. 
Later in 1895, Korteweg and de Vries \cite{KdV95} showed that the soliton could be expressed 
as a solution of a rather simple one-dimensional nonlinear PDE describing small amplitude waves
in a narrow and shallow channel of water \cite{AC91}:
\begin{equation}\label{originalKdVL}
\frac{\partial\eta}{\partial\tau} =\frac{3}{2}\frac{\partial}{\partial\xi}
\left(\frac{1}{2}\eta^2+\frac{2}{3}\alpha\eta+\frac{1}{3}\sigma\frac{\partial^2\eta}{\partial\xi^2} \right),
\qquad \sigma=\frac{1}{3}h^3-\frac{Th}{\rho g},\qquad \tau \in \R^+, \quad\xi \in \R,
\end{equation}
where $\alpha$ is some constant, $g$ denotes the gravitational constant, $\rho$ is the density,
$T$ the surface tension and $\eta=\eta(\xi,\tau)$ denotes 
the surface displacement of the wave above the undisturbed water level $h$.
The equation \eqref{originalKdVL} can be written in non-dimensional, simplified 
 form by the transformation \cite{AC91}:
\begin{equation}\label{trafoKdV}
t=\frac{1}{2}\sqrt\frac{{g}{h\sigma}}\tau,\qquad
x=-\frac{\xi}{\sqrt{\sigma}}, \qquad
u=\frac{1}{2}\eta+\frac{1}{3}\alpha
\end{equation}
to obtain the usual KdV equation
\begin{equation}\label{simpleKdVL}
u_t + 6 u u_x + u_{xxx} = 0,\qquad t \in \R^+, \quad x \in \R,
\end{equation}
(subscripts $x$ and $t$ denoting partial differentiations), 
with the soliton solution is given by
\begin{equation}\label{simpleKdVLsoliton}
u(x,t)  = \gamma \sech^2\bigl(\beta(x-ct)\bigr),\qquad t \in \R^+, \quad x \in \R.
\end{equation} 
The KdV equation \eqref{simpleKdVL} has a broad range of applications \cite{GGKM67}:
description of the asymptotic behaviour 
of small- but finite-amplitude shallow-water waves \cite{KdV95}, 
hydromagnetic waves in a cold plasma, 
ion-acoustic waves \cite{WaTa66}, 
interfacial electrohydrodynamics \cite{GHPV07},
internal wave in the coastal ocean \cite{OsSt89},
water wave power stations \cite{Cruz08},
acoustic waves in an anharmonic crystal \cite{Zab63},
or pressure pulse propagation in blood vessels \cite{KuCh06}.

In this paper, we focus on the {\em linearized KdV equation} (also known as generalized Airy equation)
 in one space dimension
\begin{equation}\label{eqKdVL}
u_t + U_1 u_x + U_2 u_{xxx} = h(t,x),\qquad t \in \R^+, \quad x \in \R,
\end{equation} 
where $h$ stands for a source term and $U_1$ and $U_2$ are real constants such that {
$U_{1} \in \R$ and $U_2>0$}. 
Recall that for $U_1=0$ and $U_2=1$ we recover the case considered 
by Zheng, Wen \& Han \cite{ZhWeHa2008}.
Although the PDE~\eqref{eqKdVL} looks very simple, it has a lot of applications,
e.g.\  Whitham \cite{Whit74} used it for the modelling of the 
propagation of long waves in the shallow water equations, see also \cite{Wri05}. 

We emphasize the fact that the restriction of the solution to 
equation \eqref{eqKdVL} to a finite interval  is not periodic. 
Thus concerning the numerical simulation, we cannot use the 
FFT method and we consider instead the equation set on an interval and supplemented with 
specially designed boundary conditions.
Since the linear PDE \eqref{eqKdVL} is defined on an unbounded domain, 
one has to confine the unbounded domain in a numerical finite
computational domain for simulation.
A common used method in such situation consists in reducing the computational domain by introducing 
{\em artificial boundary conditions}. 
Such artificial boundary conditions are constructed with the goal to approximate 
the exact solution on the whole domain restricted to the computational one. 
They are called {\em absorbing boundary conditions (ABCs)} if they lead to a well-posed initial 
boundary value problem where some energy is absorbed at the boundary. 
If the approximate solution coincides on the computational domain with the 
exact solution on the whole domain, they are called {\em transparent boundary conditions (TBCs)}. 
See \cite{review} for a review on the techniques used to construct such transparent or 
artificial boundary conditions for the Schr\"odinger equation. 

The linearity property of equation \eqref{eqKdVL} allows to use many analytical tools 
such as the Laplace transform. 
Using this tool, Zheng, Wen \& Han \cite{ZhWeHa2008} derived the exact TBCs 
for equation \eqref{eqKdVL} at fixed boundary points and then 
obtained an initial boundary value problem "equivalent" to the problem in the whole space domain. 
Moreover, using a dual Petrov-Galerkin scheme \cite{She03} the authors proposed a numerical approximation 
of this initial boundary value problem. 
Thus the derivation in \cite{ZhWeHa2008} of the adapted boundary conditions 
is carried out at the continuous level and then discretized afterwards.
Recently, Zhang, Li and Wu \cite{ZLW14} revisited the approach of Zheng, Wen \& Han  \cite{ZhWeHa2008}
and proposed a fast approximation of the exact TBCs based on Pad\'e approximation of the 
Laplace-transformed TBCs.

In this paper we will follow a different strategy: 
we first discretize the equation \eqref{eqKdVL} with respect to time and space 
and then derive the suitable artificial boundary conditions
for the fully discrete problem using the
$\mathcal{Z}$-transformation. The goal of this paper is therefore to derive
analogous conditions of the transparent boundary conditions obtained by the authors 
in \cite{ZhWeHa2008} but in the fully discrete case. 
These discrete artificial boundary conditions are superior since
they are by construction perfectly adapted to the used interior scheme
and thus retain the stability properties of the underlying
discretization method and theoretically do not produce any reflections 
when compared to the discrete whole space solution. 
However, there will be some small errors induced by the numerical root finding routine and the 
numerical inverse $\mathcal{Z}$-transformation and also later due to the fast sum-of-exponentials approximation.
Let us finally remark that there exists also 
an alternative approach in this ``discrete spirit'', 
namely to use discrete multiple scales, following the work of Schoombie \cite{Sch92}.

The paper is organized as follow. 
In Section~\ref{continuous} we use the ideas of Zheng, Wen \& Han \cite{ZhWeHa2008} 
to obtain the TBCs for the linearized KdV equation \eqref{eqKdVL} 
and we briefly recall the results given in \cite{ZhWeHa2008} for the special case
$U_1=0$ and $U_2=1$. 
In Section~\ref{Fullydisc} we present an appropriate space and time 
discretization and explain the procedure to derive 
the artificial boundary conditions for 
the purely discrete problem mimicking the ideas presented in 
Section~\ref{continuous}. 
Since exact ABCs are too time-consuming, especially for higher 
dimensional problems, we propose in Section~\ref{s:expo} to use a 
sum-of-exponentials approach \cite{AES03}, to speed
up the (approximate) computation of the discrete convolutions at the boundaries.
Finally, in Section~\ref{num} we present some
numerical benchmark examples from the literature to illustrate our findings.

\section{Transparent boundary conditions for the continuous case}\label{continuous}
The motivation for this section is twofold. 
First, we briefly recall from the literature the construction of TBCs for the 
1D linearized KdV equation \eqref{eqKdVL} 
for the special case $U_1=0$ and $U_2=1$ and the well-posedness 
of the resulting initial boundary value problem \cite{ZhWeHa2008}. 

Secondly, we extend the derivation of TBCs
to the generalized case {$U_1 \in \R$} and $U_2>0$; these results will serve us
as a guideline for the completely discrete case in Section~\ref{Fullydisc}. 

To do so, we consider the Cauchy problem
\begin{align}
        u_t+U_1 u_x+ U_2 u_{xxx}&=h(t,x), \qquad t \in \R^+,\quad x \in \R, \label{eqCauchy1} \\
        u(0,x)&=u_0(x), \qquad x \in \R, \label{eqCauchy2} \\
        u &\to 0, \qquad x \to \pm\infty, \label{eqCauchy3}
\end{align}
where (for simplicity) the initial function $u_0$ and the source term $h$ are assumed to be compactly 
supported in a finite computational interval $[a,b]$, with $a<b$ and 
where {$U_1\in \R$} and 
$U_2>0$ are given constants. 
For the construction of TBCs in the case of non-compactly
supported initial data we refer the interested reader to \cite{Eh2008}.

\begin{remark} \label{rmk:U2}
  {Changing the sign of $U_2$ means to reverse the time direction, we therefore only consider the positive case. Moreover, if one reads the equation as $u_t/U_2+(U_1/U_2)u_x+u_{xxx}=h(t,x)/U_2$, we can perform the change of unknown $v(U_2t,x)=u(t,x)$ and the equation for $v$ becomes
$$
v_{\tilde{t}}+c v_x+v_{xxx}=g(\tilde{t},x), \quad \tilde{t}=U_2 t, \ c=U_1/U_2,\ g(\tilde{t},x)=h(\tilde{t}/U_2,x)/U_2.
$$
The dispersive constant $U_2$ can therefore be eliminated from the generalized Airy equation and be considered to be equal to $1$.}
\end{remark}

The construction of (continuous) artificial boundary conditions 
associated to problem 
\eqref{eqCauchy1}--\eqref{eqCauchy3} is established by considering the problem 
on the complementary of $[a,b]$, {\it i.e.}
\begin{align}
	u_t+U_1 u_x+ U_2 u_{xxx}&=0, \qquad t \in \R^+,\quad x<a \quad\text{or}\quad x>b, \label{TBC1} \\
	u(0,x)&=0, \qquad x <a \quad\text{or}\quad x>b, \label{TBC2} \\
	u &\to 0, \quad x \to \pm \infty. \label{TBC3}
\end{align}
Denoting by $\widehat{u}=\widehat{u}(s,x)$ the Laplace transform in time 
of the function $u=u(t,x)$,  
we obtain from \eqref{TBC1} the \emph{transformed exterior problems}
\begin{align}
s\widehat{u} + U_1\widehat{u}_x + U_2\widehat{u}_{xxx} &= 0, \qquad  x<a 
                           \quad\text{or}\quad x>b,\label{TBC4} \\
 \widehat{u} &\to 0, \quad x \to \pm \infty, \label{TBC4bis}
\end{align}
where $s\in\C$, with $\Real(s)>0$, stands for the argument of the transformation, 
i.e.\ the dual time variable.
The general solutions of the ODE \eqref{TBC4} are given explicitly 
by 
\begin{equation}\label{TBC5}
\widehat{u}(s,x)=c_1(s) \,e^{\lambda_1(s)x}+c_2(s) \,e^{\lambda_2(s)x}+c_3(s) \,e^{\lambda_3(s)x},
 \qquad x<a \quad\text{or}\quad x>b,
\end{equation}
where $\lambda_1(s)$, $\lambda_2(s)$, $\lambda_3(s)$ denote the roots of 
the (depressed) cubic equation
\begin{equation}\label{cubic_cont}
s+U_1 \lambda + U_2 \lambda^3 =0.
\end{equation}
{The three solutions are given by
\begin{equation}
  \label{eq:roots}
  \lambda_k(s)=\omega^{k-1}\zeta(s)-\frac{1}{3}\frac{U_1}{U_2}\frac{1}{\omega^{k-1}\zeta(s)},\quad k=1,2,3,
\end{equation}
where $\omega=\exp(2 i\pi/3)$ and  
$$
\zeta(s)=-\frac{1}{2^{1/3}}\left (
  \frac{s}{U_2}+\sqrt{\left(\frac{s}{U_2}\right)^2+\frac{4}{27}\left
      (\frac{U_1}{U_2}\right )^3} \right)^{1/3}.
$$
\begin{theorem}
\label{theo:continuous}
  The roots of the cubic equation \eqref{cubic_cont} 
possess the following {\em separation property}
\begin{equation}\label{roots}
\Real(\lambda_1(s))<0,\qquad\Real(\lambda_2(s))>0,\qquad\Real(\lambda_3(s))>0.
\end{equation}
\end{theorem}
\noindent The proof of the theorem \ref{theo:continuous} is given in Appendix \ref{proof:theo:cont}.\\

This result is crucial for defining later the TBCs;
the \textsl{separation} property allows to separate the fundamental
solutions into outgoing and incoming waves.\\
}

\begin{remark}\label{cubiceq}
Considering, as in \cite{ZhWeHa2008}, the case $U_1=0$ and $U_2=1$ we have
\begin{equation*}
\lambda_1(s)=-\sqrt[3]{s},\qquad\lambda_2(s)=-\omega\sqrt[3]{s},
\qquad\lambda_3(s)=-\omega^2\sqrt[3]{s}.
\end{equation*}
\end{remark}

Now using the decay condition \eqref{TBC4bis}, the general solution \eqref{TBC5}, 
the separation property \eqref{roots} 
and since solutions of \eqref{TBC4} have to belong to $L^2(\mathbb{R})$, we obtain
\begin{equation}\label{ccc}
c_1(s)=0\quad\text{for}\quad x\le a,\qquad
c_2(s)=c_3(s)=0 \quad\text{for}\quad x\ge b,
\end{equation}
which yields the following TBCs in the Laplace-transformed space
\begin{gather}
\widehat{u}_{xx}(s,a)- \bigl(\lambda_2(s)+\lambda_3(s)\bigr)\,\widehat{u}_x(s,a)
+\lambda_2(s)\lambda_3(s)\,\widehat{u}(s,a)=0
,\label{TBC7}\\
\widehat{u}(s,b)-\frac{1}{\lambda_1^2(s)}\,\widehat{u}_{xx}(s,b)=0,
\qquad 
\widehat{u}_x(s,b)-\frac{1}{\lambda_1(s)}\,\widehat{u}_{xx}(s,b)=0. \label{TBC6} 
\end{gather}
Since $\lambda_1$, $\lambda_2$ and $\lambda_3$ are roots of the cubic equation
\eqref{cubic_cont} we obtain immediately
\begin{equation*}
\lambda_2(s) \lambda_3(s)=-\frac{s}{U_2 \lambda_1(s)}
\qquad\text{and}\qquad\lambda_2(s)+\lambda_3(s)=-\lambda_1(s), 
\end{equation*}
and hence the transformed left TBC \eqref{TBC7} 
can be rewritten solely in terms of $\lambda_1(s)$
\begin{equation}\label{TBC7prim}
\widehat{u}(s,a)-\frac{U_2\lambda_1(s)^2}{s}\,
\widehat{u}_x(s,a)-\frac{U_2\lambda_1(s)}{s}\,\widehat{u}(s,a)=0. 
\end{equation}
Now applying the inverse Laplace transform to equations~\eqref{TBC7prim} and \eqref{TBC6} we get
\begin{gather}
u(t,a)-U_2\,\mathcal{L}^{-1}\left(\frac{\lambda_1(s)^2}{s}\right)*u_x(t,a)-U_2\,
\mathcal{L}^{-1}\left(\frac{\lambda_1(s)}{s}\right)*u_{xx}(t,a)=0, \label{TBC9}\\
u(t,b)-\mathcal{L}^{-1}\left(\frac{1}{\lambda_1(s)^2}\right)* u_{xx}(t,b)=0,
\qquad u_x(t,b)-\mathcal{L}^{-1}\left(\frac{1}{\lambda_1(s)}\right)*u_{xx}(t,b)=0, \label{TBC8} 
\end{gather}
where $\mathcal{L}^{-1}(f(s))$ stands for the inverse Laplace transform of $f$
and $*$ denotes the convolution operator. 
We emphasize that those
boundary conditions strongly depend on $U_1$ and $U_2$ through the root $\lambda_1(s)$.

\begin{remark}{\label{TBCZheng}}\label{r22}
Considering, as in \cite{ZhWeHa2008}, 
the special case $U_1=0$ and $U_2=1$ we easily obtain from \eqref{TBC9}--\eqref{TBC8}
\begin{gather}
u(t,a)-I_t^{1/3}u_x(t,a)+I_t^{2/3}u_{xx}(t,a)=0, \label{TBCZheng9}\\
u(t,b)-I_t^{2/3}u_{xx}(t,b)=0,
\qquad u_x(t,b)+I_t^{1/3}u_{xx}(t,b)=0, \label{TBCZheng8}
\end{gather}
where $I_t^p$ with $p >0$ is the nonlocal-in-time fractional integral operator 
given by the {\it Riemann-Liouville formula}
\begin{equation*}
I_t^p f(t) =\frac{1}{\Gamma(p)}\int_0^t (t-\tau)^{p-1}f(\tau) \,d\tau,
\end{equation*}
where $\Gamma(z)=\int_0^{+\infty} e^{-t}t^{z-1}\,dt$ is the Gamma function. 
We refer to \cite{ZhWeHa2008} for more details.
\end{remark}

To summarize our findings so far, the derived initial boundary value problem reads
\begin{gather}
u_t+U_1u_x+U_2u_{xxx}=0, \qquad t \in \R^+, ~ x \in [a,b], \label{TBC10} \\
u(0,x)=u_0(x), \qquad x \in [a,b], \label{TBC11} \\
u(t,a)-U_2\,\mathcal{L}^{-1}\biggl(\frac{\lambda_1(s)^2}{s}\biggr)*u_x(t,a)-U_2\,\mathcal{L}^{-1}
\biggl(\frac{\lambda_1(s)}{s}\biggr)*u_{xx}(t,a)=0,\label{TBC12} \\
u(t,b)-\mathcal{L}^{-1}\Bigl(\frac{1}{\lambda_1(s)^2}\Bigr)* u_{xx}(t,b)=0,\label{TBC13} \\
u_x(t,b)-\mathcal{L}^{-1}\Bigl(\frac{1}{\lambda_1(s)}\Bigr)*u_{xx}(t,b)=0. \label{TBC14}
\end{gather}
Note that a solution of \eqref{TBC10}--\eqref{TBC14}
can be regarded as the restriction on $[a,b]$ 
of the solution on the whole space domain.

For the special case $U_1=0$ and $U_2=1$ (cf.\ Remark~\ref{r22}) 
the following stability theorem is shown in \cite{ZhWeHa2008}.

\begin{theorem}
[\cite{ZhWeHa2008}]
The initial boundary value problem \eqref{TBC10}--\eqref{TBC14} 
for $U_1=0$ and $U_2=1$ is $L^2$-stable. 
More precisely, for any $t>0$, there is a constant positive number $c(t)$ such that
\begin{equation}
\label{TBC15}
\int_a^b u^2(t,x)\,dx 
\le c(t) \,\Bigl(\int_a^b u_0^2(x)\,dx
+\int_0^t\int_a^b h^2(t,x)\,dxdt\Bigr).
\end{equation}
\end{theorem}
In the sequel we use the same procedure to obtain the artificial boundary conditions 
for the fully discrete case, i.e.\ for the discretized version of the equation \eqref{eqKdVL}.
These so-called discrete artificial boundary conditions are better adapted to 
the numerical scheme and thus do not alter the stability properties. 
Also, they do not suffer from discretization errors of convolution integrals
and theoretically do not produce any unphysical reflections. 
However, since some steps in the calculation of the convolution coefficients, like the 
root finding and the inverse $\mathcal{Z}$-transformation have to be 
done numerically, this procedure will lead to some small errors.

\section{Discrete transparent boundary conditions} \label{Fullydisc}

In this section we present how to obtain the artificial boundary conditions 
in the fully discrete case for the problem \eqref{eqCauchy1}--\eqref{eqCauchy3}. 
For simplicity we focus here on the case without source term, {\it i.e.}\ 
we assume $h(t,x)=0$ for all $t>0$ and $x\in\R$. 
Moreover we consider the problem restricted to the computational interval $[a,b]$ for
the  finite time $t\in[0,T]$, {\it i.e.}\
\begin{align}
	u_t+U_1u_x+U_2u_{xxx}&=0, \qquad t \in [0,T], \quad x \in [a,b], \label{fulldisc1} \\
	u(0,x)&=u_0(x),\qquad x \in [a,b]. \label{fulldisc2} 
\end{align}

Let us denote by $(t_n)_{0 \le n \le N}$ a uniform subdivision 
of the time interval $[0,T]$
 given by $t_n=n\Delta t$ with the temporal step size $\Delta t=T/N$:
\begin{equation*}
      0=t_0 < t_1 < \cdots < t_{N-1} < t_N=T.
\end{equation*}
We also define $(x_j)_{0 \le j \le J}$ a uniform subdivision of $[a,b]$
 given by $x_j=a+j\Delta x$ with the spatial step size $\Delta x=(b-a)/J$:
\begin{equation*}
          a=x_0 < x_1 < \cdots < x_{J-1} < x_J=b.
\end{equation*}
We emphasize here that the temporal discretization must remain uniform
due to the usage of the $\mathcal{Z}$-transform to derive the discrete TBCs. 
On the other hand,
the space discretization in the interior domain could have been non uniform.
In the following, we denote by $u_j^{(n)}$ the pointwise approximation of 
the solution $u(t_n,x_j)$.

We will consider in the sequel two different numerical schemes based
on trapezoidal rule in time (semi discrete Crank-Nicolson
approximation). 
The first one is the \textsl{Rightside Crank-Nicolson} (proposed 
by Mengzhao \cite{Me1983})
(R-CN) scheme defined for $U_1=0$ and $U_2>0$. It reads
\begin{multline}\label{fulldisc4}
\frac{u_j^{(n+1)}-u_j^{(n)}}{\Delta t} 
+\frac{U_2}{2(\Delta x)^3}\left(u_{j+2}^{(n+1)}-3u_{j+1}^{(n+1)}+3u_j^{(n+1)}-u_{j-1}^{(n+1)}\right)  \\
+\frac{U_2}{2(\Delta x)^3}\left(u_{j+2}^{(n)}-3u_{j+1}^{(n)}+3u_j^{(n)}-u_{j-1}^{(n)}\right)=0.
\end{multline}
The second one is the \textsl{Centered Crank-Nicolson} (C-CN) 
scheme \cite{Me1983} which is used for the generalized linear
Korteweg-de Vries equation 
\eqref{eqKdVL} where {$U_1\in \R$} and $U_2\geq 0$. It reads
\begin{equation}\label{fulldisc4bis}
\begin{split}
\frac{u_j^{(n+1)}-u_j^{(n)}}{\Delta t} &+ \frac{U_1}{4\Delta x}\left(u_{j+1}^{(n+1)}-u_{j-1}^{(n+1)}\right)
+\frac{U_1}{4\Delta x}\left(u_{j+1}^{(n)}-u_{j-1}^{(n)}\right) \\
&+\frac{U_2}{4(\Delta x)^3}\left(u_{j+2}^{(n+1)}-2u_{j+1}^{(n+1)}+2u_{j-1}^{(n+1)}-u_{j-2}^{(n+1)}\right)  \\
&+\frac{U_2}{4(\Delta x)^3}\left(u_{j+2}^{(n)}-2u_{j+1}^{(n)}+2u_{j-1}^{(n)}-u_{j-2}^{(n)}\right)=0.
\end{split}
\end{equation}
Here, the convection term is discretized in a centered way. 
Indeed, using simply an {\it upwind Crank-Nicholson} scheme for the
first order term and (R-CN) scheme for the third order term leads to a
strongly dissipative scheme.  

Both schemes are absolutely stable and their truncation errors are
respectively
\begin{equation}
  \label{eq:order}
  \begin{split}
    E_{\rm R-CN}&=O(\Delta x  + \Delta t^2),\\
    E_{\rm C-CN}&=O(\Delta x^2  + \Delta t^2).\\
  \end{split}
\end{equation}
The stencil of the different scheme involves respectively 4 nodes for (R-CN) 
and 5 nodes for (C-CN) schemes. This structure will have a strong influence
on the computation of the roots for the corresponding equation to
\eqref{cubic_cont} at the discrete level. 
For the (R-CN) scheme, we
will recover as in the continuous case a cubic equation, but a quartic
equation for the (C-CN) scheme. 
The later case will turn out to be more difficult.

\subsection{Discrete artificial boundary conditions for (R-CN) scheme}\label{U_1equal0}
Let us first consider the (R-CN) scheme \eqref{fulldisc4} for the interior problem, 
i.e.\ with a spatial index $j$ such that $1\le j\le J-2$. 
Let us recall that this scheme is only valid in the case $U_{1}=0$ and $U_{2}>0$.
For this scheme, as in the continuous case, we will obtain one
boundary condition at point $x_0=a$ and two boundary conditions at the
right side which will involve the two nodes $x_{J-1}$ and $x_J$, 
cf.\ the continuous ABCs \eqref{TBC12}--\eqref{TBC14}.

In order to derive appropriate artificial boundary conditions, we follow the same 
procedure as in Section~\ref{continuous}, but on a purely discrete level. 
First we apply the $\mathcal{Z}$-transform with respect to the time index $n$, 
which is the discrete analogue of the Laplace transform in time, to the
partial difference equation \eqref{fulldisc4}. 
We refer the reader to the appendix of \cite{review, Eh2001, EA01} for a proper definition 
of the $\mathcal{Z}$-transform and its basic properties. 
The standard definition reads
\begin{equation}
  \label{Ztransf}
  \hat{u}(z)=\ZT{(u^n)_n}(z)=\sum_{k=0}^\infty
  u^kz^{-k},\qquad |z|>R\ge 1,
\end{equation}
where $R$ is the convergence radius of the Laurent series and $z\in\mathbb{C}$.

Denoting by $\widehat{u}_j=\widehat{u}_j(z)$ the $\mathcal{Z}$-transform 
of the sequence $(u_j^{(n)})_{n\in\N_0}$ we obtain from 
\eqref{fulldisc4} the homogeneous {\em third order difference equation} 
\begin{equation}
\label{fulldisc5}
\widehat{u}_{j+2}-3\widehat{u}_{j+1}
+\left(3+\frac{2(\Delta x)^3}{U_2\Delta t}\frac{z-1}{z+1}\right)
\widehat{u}_{j}-\widehat{u}_{j-1}=0, \qquad 1 \le j \le J-2.
\end{equation}
It is well-known that homogeneous difference equations
 with constant coefficients possess solutions
of the power form $\widehat{u}_{j}=\sum_{k}c_k(z)\,r^j_k(z)$, 
where $r=r(z)$ solves the cubic equation
\begin{equation}
\label{fulldisc6}
r^3-3r^2+\left(3+p\right)r-1=0,
\end{equation}
{with
$$
p=\mu\frac{z-1}{z+1}, \quad \mu=\frac{2(\Delta x)^3}{U_2\Delta t}.
$$}
Equation \eqref{fulldisc6} admits three fundamental solutions denoted 
here by $r_1$, $r_2$ and $r_3$ that can be computed analytically or
numerically up to a very high precision.
Thus the general solution of \eqref{fulldisc5} on the exterior domains
is of the form 
\begin{equation*}
\widehat{u}_{j}(z)= c_1(z)\,r_1^j(z) + c_2(z)\,r_2^j(z) + c_3(z)\,r_3^j(z),
\qquad j\le 1 \text{ or } j\ge J-2. 
\end{equation*}
Let 
{
$$
\zeta(z)=\left ( \frac{-p+\sqrt{p^{2}+\frac{4}{27}p^{3}}}{2} \right )^{1/3}.
$$}
The three solutions of \eqref{fulldisc6} are
{
\begin{equation} \label{expl_roots}
r_j(z)=\omega^{j-1} \zeta(z)-\frac{p}{3}\frac{1}{\omega^{j-1}\zeta(z)}+1,\quad j=1,2,3.
\end{equation}}
{If we consider $z=\rho e^{i\theta}$, the roots $r_j(z)$ may be discontinuous due to branch changes with respect to $\theta$ in the complex plane. One example of this phenomenom can be seen on left part of figure \ref{fig:branch} where we plot the evolution of roots $r_i$ for $\rho=(1+10^{-2})$ and $\mu=1$. On one side, we clearly see that branch changes occur. On the other side, we always have simultaneously one root inside the unit disk and two outside. Instead of considering roots $r_k$, it is more convenient to identify roots by continuity. We refer to these continuous roots as $\ell_k$ (see right part of figure \ref{fig:branch}). We therefore have the following theorem.
  \begin{figure}[htbp]
    \begin{tabular}{cc}
    \includegraphics[width=0.48\textwidth]{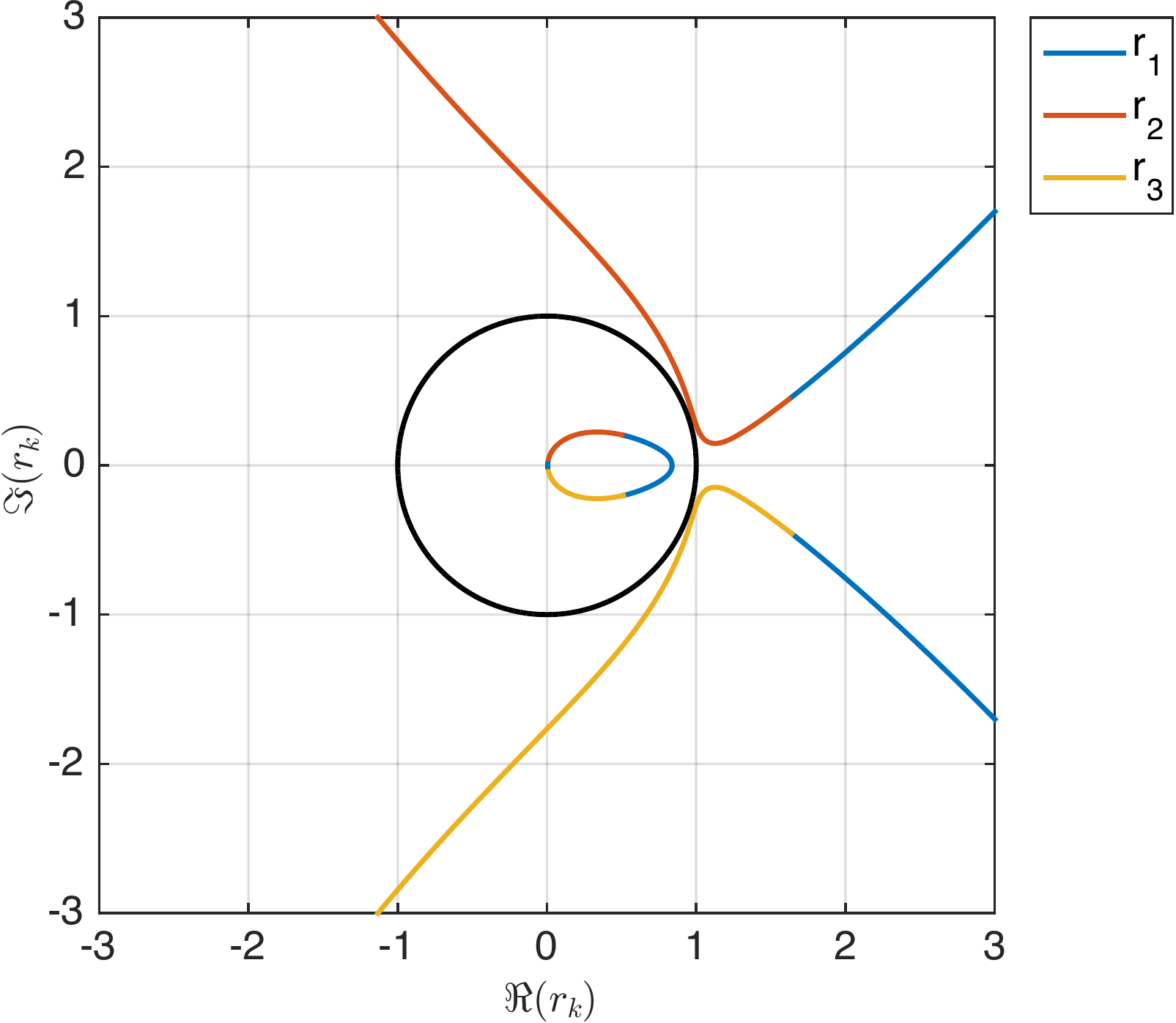}  &
    \includegraphics[width=0.48\textwidth]{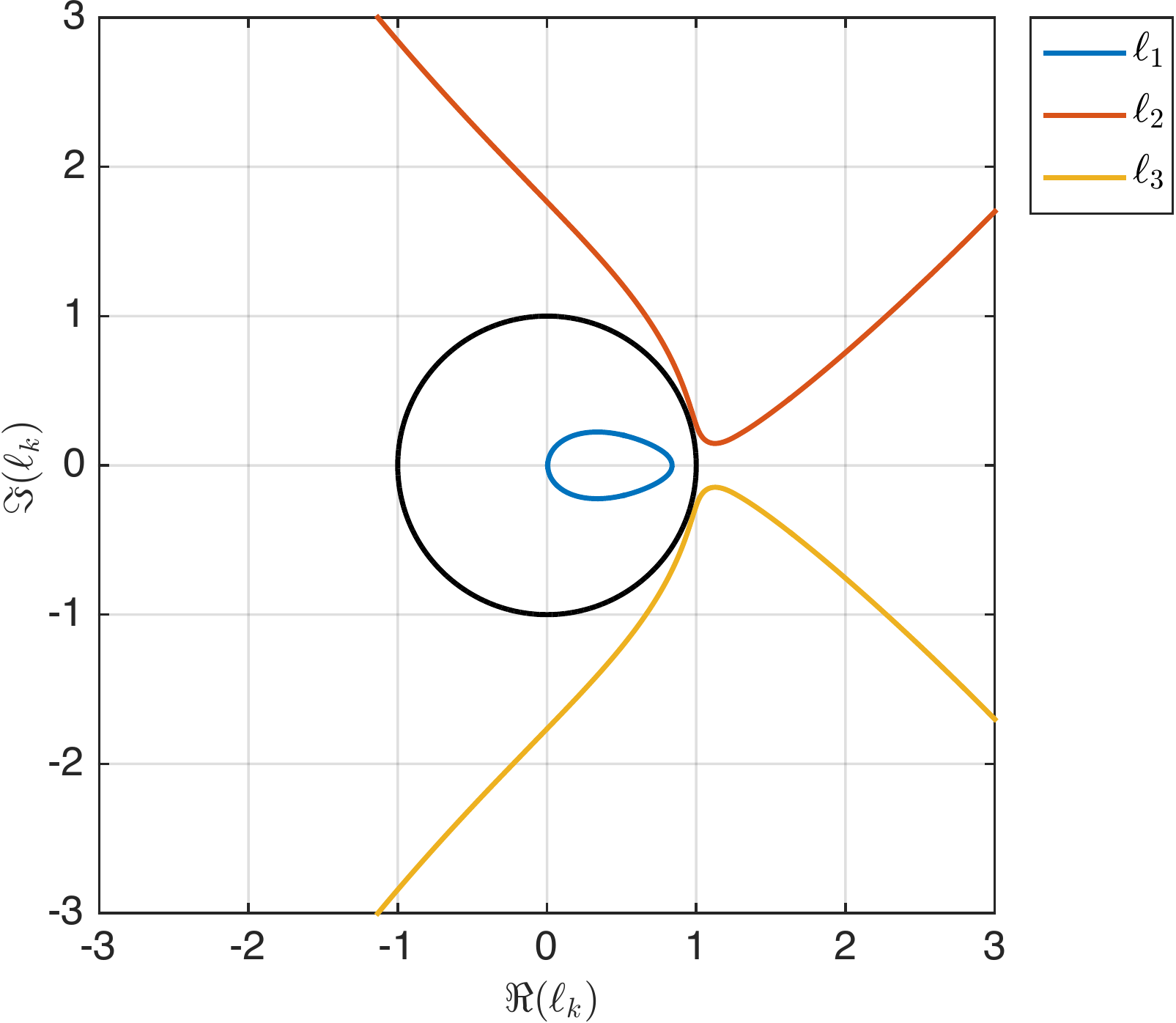} \\
    Branch change of roots $r_k$ & Continuous roots $\ell_k$ 
    \end{tabular}
    \caption{Discontinuous and continuous roots $r_k$ and $\ell_k$}
    \label{fig:branch}
  \end{figure}
\begin{theorem}\label{theo:discr3}
For any $\Delta x>0$, $\Delta t>0$ and $|z|>1$, the continuous roots of the cubic algebraic equation \eqref{fulldisc6} are well separated according to
\begin{equation}\label{cubic1}
|\ell_1(z)|<1,\qquad|\ell_2(z)|>1,\qquad|\ell_3(z)|>1,\qquad\text{for all }z,
\end{equation}
which defines the {\em discrete separation property}.
\end{theorem}
\noindent The proof of the theorem \ref{theo:discr3} is given in Appendix \ref{proof:theo:discr3}.}

{\noindent Like in the continuous case \eqref{ccc}, using the decay condition we obtain,}
\begin{itemize}
\item for the left exterior domain $c_1(z)=0$, $j\le 1$ and thus 
       $\widehat{u}_{j}(z)=c_2(z)\,\ell_2^j(z)+c_3(z)\,\ell_3^j(z)$, $j\le 1$
\item for the right exterior domain $c_2(z)=c_3(z)=0$, $j\ge J-2$ and thus
$\widehat{u}_{j}(z)=c_1(z)\,\ell_1^j(z)$, $j\ge J-2$
\end{itemize}

Let us now derive the boundary conditions for (R-CN) scheme.

\noindent{\bf Left boundary.} 
On the left boundary we only need one relation. It is easy to see that
\begin{equation}\label{fulldisc8}
\widehat{u}_{j+1}(z)-\bigl(\ell_2(z)+\ell_3(z)\bigr)\,\widehat{u}_{j}(z)
+\ell_2(z)\ell_3(z)\,\widehat{u}_{j-1}(z)=0.
\end{equation}
Applying it for $j=1$ and denoting by $\IZT{f(z)}$ the inverse 
$\mathcal{Z}$-transform of $f(z)$, we obtain
\begin{equation}
\label{fulldisc10}
\IZT{\ell_2(z)\ell_3(z))}*_d u_0^{(n)} -\IZT{\ell_2(z)+\ell_3(z))}*_d u_{1}^{(n)}+u_{2}^{(n)}=0, 
\quad n=0,1,2\dots,
\end{equation}
where $*_d$ stands for the discrete convolution with respect to the temporal
index $n$:
\begin{equation*}
P *_d u_i^{(n)}=\sum_{k=0}^{n}P^{(k)} u_i^{(n-k)},
\end{equation*}
for $P=(P^{k})$ a sequence and $i$ an integer. Let us denote the 
convolution kernels by 
$k_{1,R}(z)=\ell_2(z)+\ell_3(z)$ and $k_{2,R}(z)=\ell_2(z)\ell_3(z)$ 
and by $Y_{i,R}$ the
sequences of the inverse $\mathcal{Z}$-transform of kernel $k_{i,R}$
{\it i.e.} $Y_{i,R}=\IZT{k_{i,R}(z)}$. 
Then \eqref{fulldisc10} can be written
\begin{equation}\label{fulldisc10bis}
Y_{2,R}*_d u_0^{(n)}-Y_{1,R}*_d u_{1}^{(n)}+u_{2}^{(n)}=0, \quad n=0,1,2\dots.
\end{equation}

\noindent{\bf Right boundary.} 
On the right boundary we need two relations
 since the fully discrete scheme involved four grid points. 
It is easy to see that
\begin{equation}
\label{fulldisc7}
\widehat{u}_{j+2}(z)=\ell_1(z)^2\,\widehat{u}_{j}(z),\qquad\text{and}
\qquad\widehat{u}_{j+1}(z)=\ell_1(z)\,\widehat{u}_{j}(z),
\end{equation}
Applying them for $j=J-2$ and using inverse $\mathcal{Z}$-transformation, we obtain
\begin{equation}\label{fulldisc9}
u_J^{(n)}-Y_{4,R} *_d u_{J-2}^{(n)}=0,
\quad u_{J-1}^{(n)}-Y_{3,R} *_d u_{J-2}^{(n)}=0,\quad n=0,1,2\dots,
\end{equation}
where $k_{3,R}(z)=\ell_1(z)$ and $k_{4,R}(z)=\ell_1^2(z)$.

\begin{remark}\label{obtaincoeff}
We can easily obtain continuous roots $\ell_k$. We compute roots $r_k(z)$ where $z=\rho e^{i\theta}$ for a fixed radius $\rho$ and various angle $\theta$ thanks to \eqref{expl_roots}. For each values of $\theta$, we sort roots thanks to the relation $r_1(\theta)< r_2(\theta)<r_3(\theta)$.
\end{remark}

\begin{figure}[ht]
\begin{center}
\begin{tabular}{cc}
\includegraphics[width=0.4\textwidth]{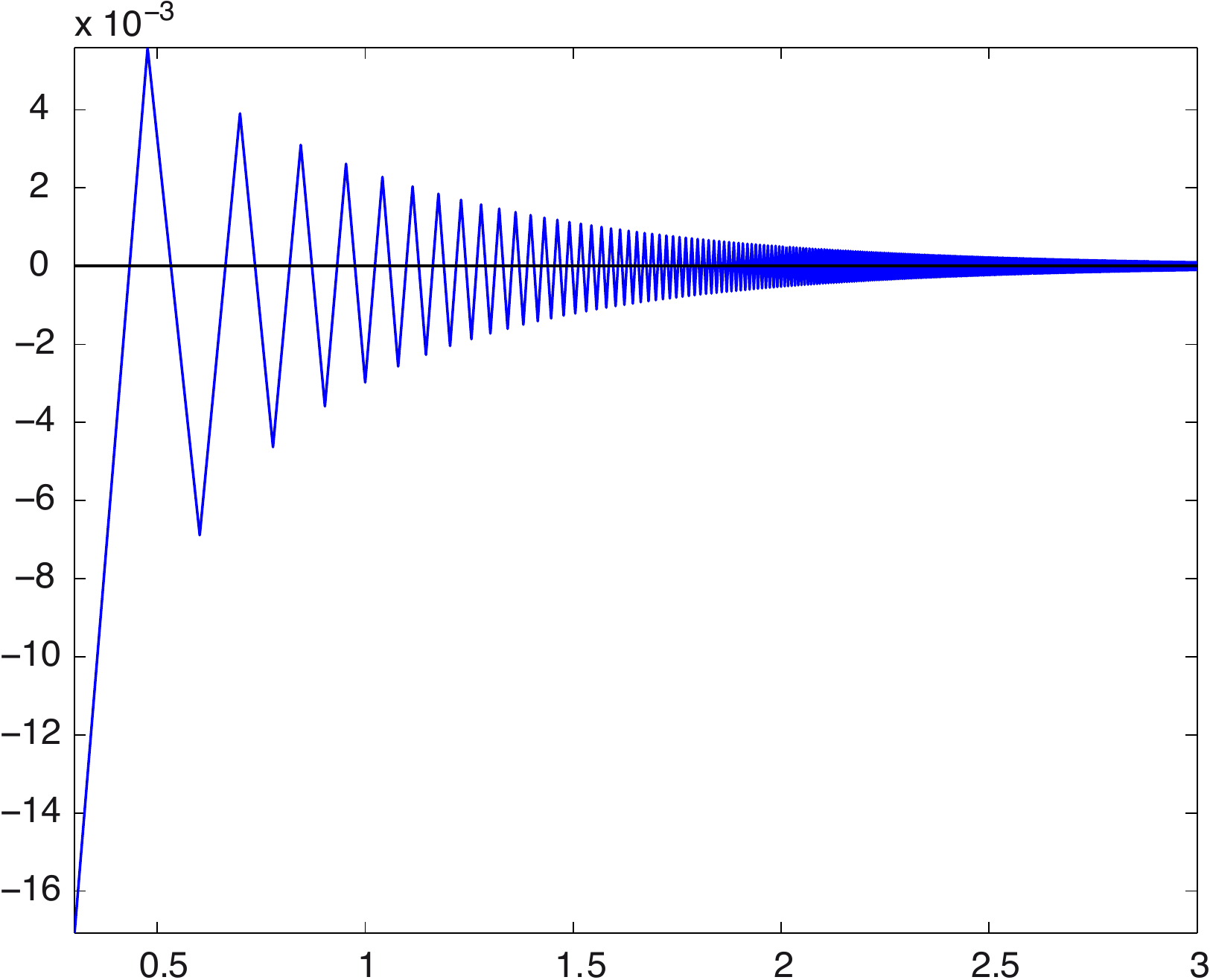}
&
\includegraphics[width=0.4\textwidth]{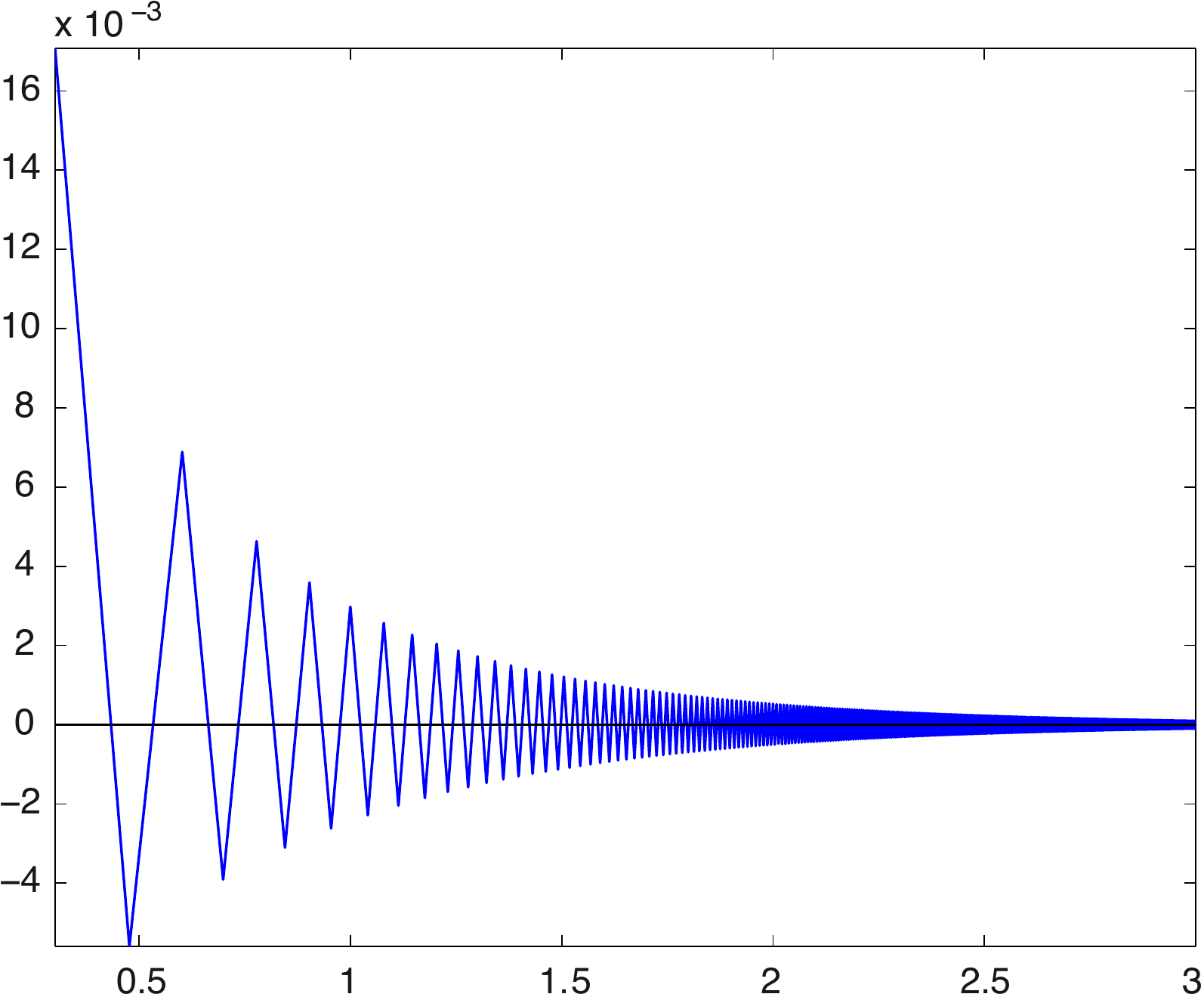} \\
$\mathcal{Z}^{-1}(k_{1,R}(z))$ & $\mathcal{Z}^{-1}(k_{3,R}(z))$ \\
\includegraphics[width=0.4\textwidth]{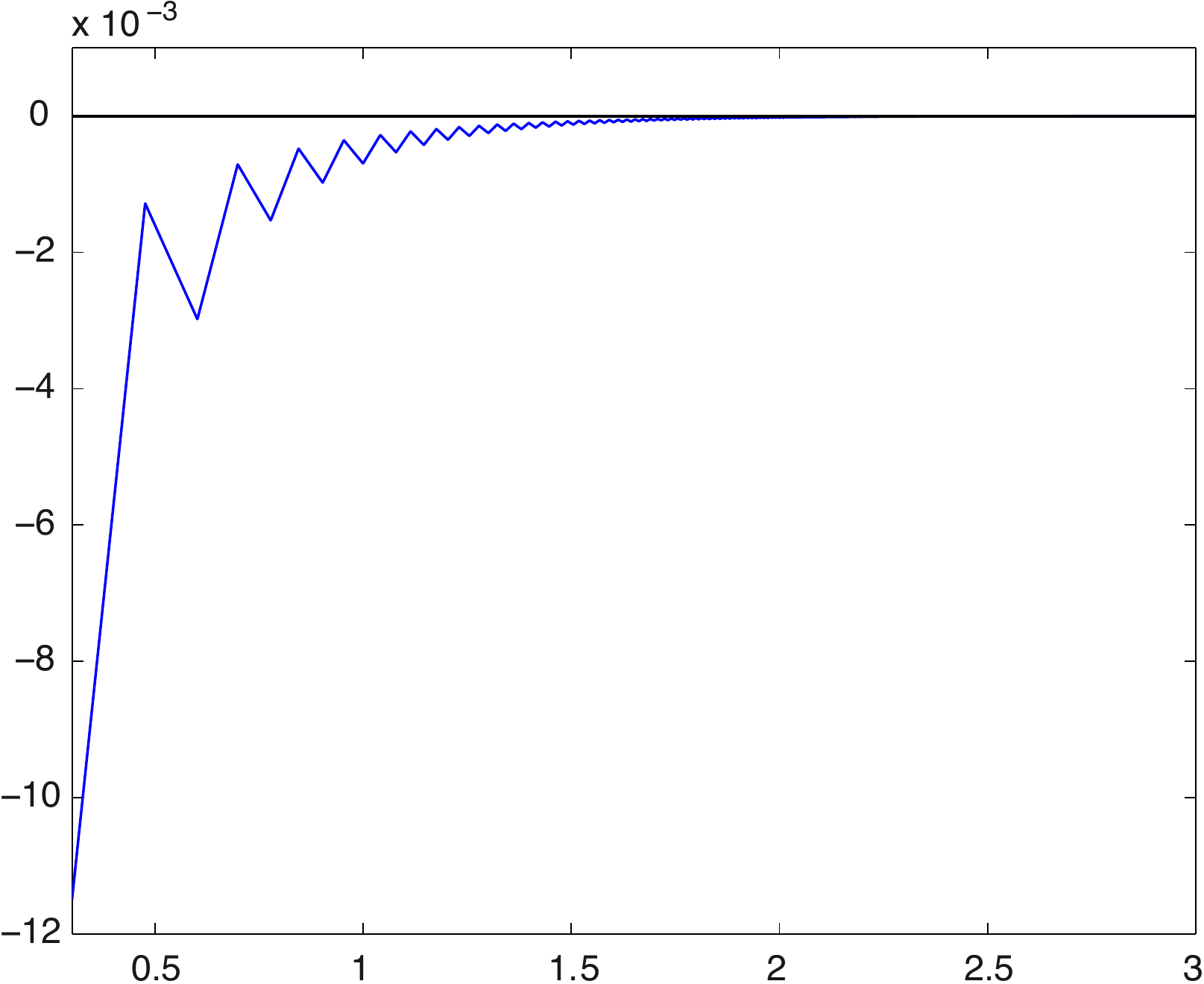}
&
\includegraphics[width=0.4\textwidth]{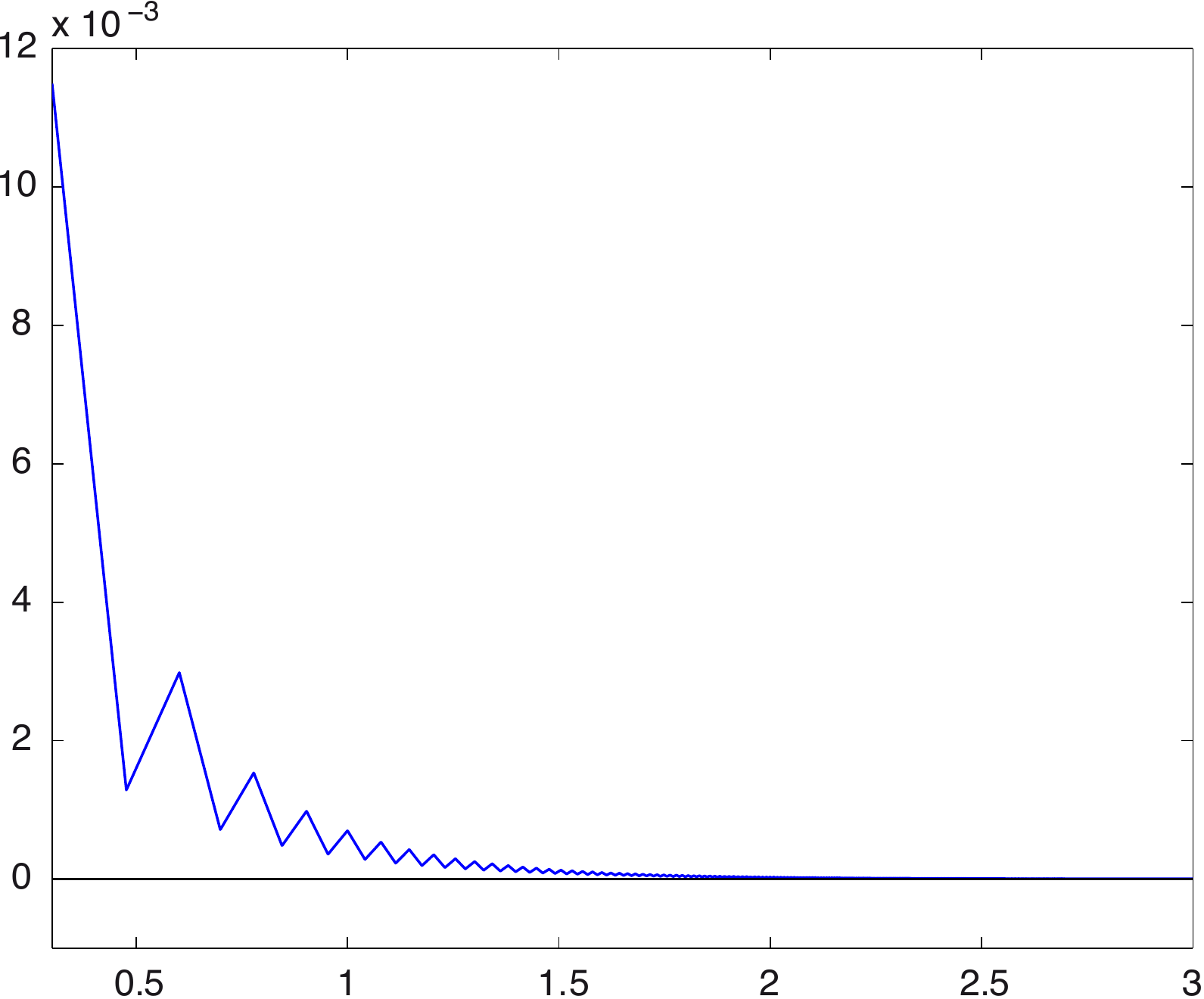}\\
$\mathcal{Z}^{-1}(\xi(z)k_{1,R}(z))$ & $\mathcal{Z}^{-1}(\xi(z)k_{3,R}(z))$ \\
\end{tabular}
\end{center}
\caption{Coefficients of the inverse $\mathcal{Z}$-transform for the 
kernels $k_{i,R}(z)$, $i=1,3$ and kernels $\xi(z)k_{i,R}(z)$, $i=1,3$ 
with $\Delta t=4/2560$, 
$\Delta x=12/5000$ and $r=1.001$.}
\label{kernelcoeff}
\end{figure}

\begin{remark}\label{remcoeff}
We draw on Figure \ref{kernelcoeff} (top figures) the behaviour of the inverse
$\mathcal{Z}$-transform only for the two kernels $k_{1,R}(z)$ and $k_{3,R}(z)$ (the behavior of $k_{2,R}(z)$ and $k_{4,R}(z)$ being analogous). We clearly see that
the signs of the coefficients alternate. This will possibly create
subtractive cancellation errors when we will use them in the boundary
convolutions. Following \cite{AES03,AESS12}, we modify the
convolution kernels. The idea is to multiply a
$\mathcal{Z}$-transformed kernel $k_{i,R}(z)$ by $\xi(z)=1+z^{-1}$ which
corresponds to add two neighboured values in the series $\mathcal{Z}^{-1}\{k_{i,R}(z)\}$. We draw on Figure \ref{kernelcoeff} (bottom figures) the behavior of the inverse
$\mathcal{Z}$-transform for $\xi(z)k_{1,R}(z)$ and $\xi(z)k_{3,R}(z)$. We can clearly see that the signs of the coefficients do not alternate anymore.
In the sequel we introduce the notations $Y_{i,R}^\xi=\IZT{\xi(z)k_{i,R}(z)}$. We refer to section \ref{num} for more details on the numerical procedure used to compute the inverse $\mathcal{Z}$-transform for a kernel.
\end{remark}

Remark~\ref{remcoeff} yields the algorithm used in Section~\ref{num} 
to solve numerically the problem. 
Assuming that the solution on the previous time level 
$(u_j^{(n)})_{0\le j \le J}$ is known, 
$(u_j^{(n+1)})_{0\le j \le J}$ is given for $n \ge0$ by
\begin{equation}\label{algo1}
\begin{cases}
Y_{2,R}^\xi *_d u_0^{(n+1)}-Y_{1,R}^\xi*_d u_1^{(n+1)}+ u_2^{(n+1)}=-u_2^n, \\[2mm]
-\alpha u_{j-1}^{(n+1)}+(3\alpha+1)u_{j}^{(n+1)}-3\alpha u_{j+1}^{(n+1)}+\alpha u_{j+2}^{n+1} \\[2mm]
\hskip2cm=\alpha u_{j-1}^{(n)}-(3\alpha-1)u_{j}^{(n)}+3\alpha u_{j+1}^{(n)}-\alpha u_{j+2}^{(n)},
\qquad 1\le j\le J-2, \\[2mm]
u_{J-1}^{(n+1)}-Y_{3,R}^\xi *_d  u_{J-2}^{(n+1)} =-u_{J-1}^n, \\[2mm]
u_{J}^{(n+1)}-Y_{4,R}^\xi *_d u_{J-2}^{(n+1)} =-u_J^n, 
\end{cases}
\end{equation}
with the mesh ratio $\alpha=U_2\Delta t/(2({\Delta x})^3)$.

\subsection{Discrete artificial boundary conditions for (C-CN) scheme}\label{scheme2}
We treat here the case of the (C-CN) scheme for the interior 
nodes with a spatial index $j$ such that $2 \le j \le J-2$. 
Since we consider a difference scheme with a five points stencil, 
we need two artificial boundary conditions on each side of the
computational interval $[a,b]$. 

In order to derive suitable artificial boundary conditions for the (C-CN) scheme \eqref{fulldisc4bis} 
we follow the same procedure as in Section~\ref{U_1equal0} for (R-CN) scheme. 
First we apply the $\mathcal{Z}$-transform with respect to the time index $n$, 
denoting by $\widehat{u}_j$ the $\mathcal{Z}$-transform 
of the sequence $(u_j^{(n)})_{n\in\N_0}$ we obtain from \eqref{fulldisc4bis}
the homogeneous {\em fourth order difference equation}:
\begin{equation}
\label{fulldisc5bis}
\widehat{u}_{j+2}-\left(2-\frac{U_1 (\Delta x)^2}{U_2}\right)\widehat{u}_{j+1}
+\frac{4(\Delta x)^3}{U_2\Delta t}\frac{z-1}{z+1}\widehat{u}_{j}+\left(2-\frac{U_1 (\Delta x)^2}{U_2}\right)
\widehat{u}_{j-1}-\widehat{u}_{j-2}=0, 
\end{equation}
for the spatial index range $2 \le j \le J-2$. 

The solutions of this difference equation are again
of the power form $\widehat{u}_{j}(z)=\sum_kc_k(z)\,\ell^j_k(z)$, 
where $\ell=\ell(z)$ solves now the \textsl{quartic} equation
{
\begin{equation}\label{fulldisc6bis}
\ell^4-\left(2-a\right)\ell^3+2p\ell^2+\left(2-a\right)\ell-1=0,
\end{equation}
with $a=U_{1}(\Delta x)^{2}/U_{2}$ and $p=2\lambda (z-1)/(z+1)$.}

{
Equation \eqref{fulldisc6bis} admits four roots which can be computed numerically or analytically by the well-known Ferrari's solution formula. As in the case of the cubic equation we identify these roots by continuity and refer to them as $\ell_{k}$ for $k=1,2,3,4$.
Thus the general solution of \eqref{fulldisc5bis} is of the form}
\begin{equation*}
\widehat{u}_{j}(z)=c_1(z)\,\ell_1^j(z)+c_2(z)\,\ell_2^j(z)+c_3(z)\,\ell_3^j(z)
+c_4(z)\,\ell_4^j(z). 
\end{equation*}
{Like in the previous sections we can show the following theorem
\begin{theorem}
\label{theo:discr4}
For any $U_{1} \in \R$, $U_{2}>0$, $\Delta x>0$, $\Delta t>0$ and $|z|>1$, the continuous roots of the quartic algebraic equation \eqref{fulldisc6bis} are well separated according to
\begin{equation}\label{cubic1}
|\ell_1(z)|<1,\qquad|\ell_2(z)|<1,\qquad|\ell_3(z)|>1, \qquad|\ell_4(z)|>1, \qquad\text{for all }z,
\end{equation}
which defines the {\em discrete separation property}.
\end{theorem}
\noindent The proof of the theorem is given in Appendix \ref{proof:theo:discr4}.
}

{As previously, using theorem \ref{theo:discr4} and the decay condition  of the solution we obtain}
\begin{itemize}
\item for the left exterior domain $c_1(z)=c_2(z)=0$, $j\le 2$ and thus 
       $\widehat{u}_{j}(z)=c_3(z)\,\ell_3^j(z)+c_4(z)\,\ell_4^j(z)$, $j\le 2$
       
	\item for the right exterior domain $c_3(z)=c_4(z)=0$, $j\ge J-2$ and thus
	    $\widehat{u}_{j}(z)=c_1(z)\,\ell_1^j(z)+c_2(z)\,\ell_2^j(z)$, $j\ge J-2$
\end{itemize}

\noindent{\bf Right boundary.} 
On the right boundary we need two relations. 
It is easy to see that
\begin{equation}
\label{fulldisc7bis}
\widehat{u}_{j+2}(z)-
\bigl(\ell_1(z)+\ell_2(z)\bigr)\widehat{u}_{j+1}(z)+\ell_1(z)\ell_2(z)\widehat{u}_{j}(z)=0,
\end{equation}
\begin{equation}
\label{fulldisc7bisbis}
\widehat{u}_{j+2}(z)-2(\ell_1(z)+\ell_2(z))\widehat{u}_{j+1}(z)+(\ell_1(z)+\ell_2(z))^2
\widehat{u}_{j}(z)-(\ell_1(z)\ell_2(z))^2\widehat{u}_{j-2}(z)=0,
\end{equation}
which will give a link between $u_J^{(n)}$, $u_{J-1}^{(n)}$, $u_{J-2}^{(n)}$ 
and $u_{J-4}^{(n)}$ with $j=J-2$. 

For brevity of the notation we introduce  
$Y_{i,C}=\IZT{k_{i,C}(z)}, i=1, 2, 3, 4$ with
\begin{gather*}
k_{1,C}(z)=\ell_1(z)+\ell_2(z),\qquad
k_{2,C}(z)=\bigl(\ell_1(z)+\ell_2(z)\bigr)^2, \\
k_{3,C}(z)=\ell_1(z)\ell_2(z),\qquad 
k_{4,C}(z)=\bigl(\ell_1(z)\ell_2(z)\bigr)^2,
\end{gather*}
and we obtain from \eqref{fulldisc8bis1}--\eqref{fulldisc8bisbis} 
\begin{equation}\label{fulldisc9bis}
\begin{split}
u_J^{(n)}-Y_{1,C} *_d u_{J-1}^{(n)}+Y_{3,C} *_d u_{J-2}^{(n)}&=0,\\
u_{J}^{(n)}-2Y_{1,C} *_d u_{J-1}^{(n)}+Y_{2,C} *_d u_{J-2}^{(n)}-Y_{4,C} *_d
u_{J-4}^{(n)}&=0.
\end{split}
\end{equation}

\noindent{\bf Left boundary.} 
On the left boundary we need now two relations. We can easily verify
\begin{equation}\label{fulldisc8bis1}
\widehat{u}_{j}(z)-\bigl(\ell_3(z)+\ell_4(z)\bigr)\,\widehat{u}_{j-1}(z)
+\ell_3(z)\ell_4(z)\,\widehat{u}_{j-2}(z)=0,
\end{equation}
\begin{equation}\label{fulldisc8bisbis}
\widehat{u}_{j+2}(z)-2(\ell_3(z)+\ell_4(z))\widehat{u}_{j+1}(z)+(\ell_3(z)
+\ell_4(z))^2\widehat{u}_{j}(z)-(\ell_3(z)\ell_4(z))^2\widehat{u}_{j-2}(z)=0,
\end{equation}
which give a link between $u_0^{(n)}$, $u_1^{(n)}$, $u_2^{(n)}$, 
$u_3^{(n)}$ and $u_4^{(n)}$ setting $j=2$. 
Indeed, denoting by $Y_{i,C}=\IZT{k_{i,C}(z)}$, $i=5, 6, 7, 8$ 
with \begin{gather*}
k_{5,C}(z)=\ell_3(z)+\ell_4(z),\qquad 
k_{6,C}(z)=(\ell_3(z)+\ell_4(z))^2, \\
k_{7,C}(z)=\ell_3(z)\ell_4(z), \qquad
k_{8,C}(z)=(\ell_3(z)\ell_4(z))^2, 
     \end{gather*}
we obtain from \eqref{fulldisc8bis1}-\eqref{fulldisc8bisbis} 
\begin{equation}\label{fulldisc10bis1}
\begin{split}
Y_{7,C}*_d u_0^{(n)}-Y_{5,C} *_d u_{1}^{(n)}+ u_{2}^{(n)}&=0,\\
-Y_{8,C}*_du_{0}^{(n)}+ Y_{6,C} *_d u_{2}^{(n)}-2Y_{5,C} *_d u_{3}^{(n)}+ u_{4}^{(n)}&=0. 
\end{split}
\end{equation}

Following Remark~\ref{remcoeff}, we finally obtain the algorithm used
in Section~\ref{num} to solve numerically the problem. 
Assuming that the solution $(u_j^{(n)})_{0\le j\le J}$ on the previous time level 
is known, then
$(u_j^{(n+1)})_{0\le j\le J}$ is given for $n \ge 0$ by
\begin{equation}\label{algo2}
\begin{cases}
Y_{7,C}^\xi*_d u_0^{(n)}-Y_{5,C}^\xi *_d u_{1}^{(n)}+ u_{2}^{(n)}=-u_2^{(n-1)}, \\[2mm]
-Y_{8,C}*_du_{0}^{(n)}+ Y_{6,C} *_d u_{2}^{(n)}-2Y_{5,C} *_d u_{3}^{(n)}+ u_{4}^{(n)}=-u_4^{(n-1)}, \\[2mm]
\displaystyle{-\frac{\alpha}{2} u_{j-2}^{(n+1)}+(\alpha-\beta)u_{j-1}^{(n+1)}+u_{j}^{(n+1)}+(-\alpha+\beta)u_{j+1}^{(n+1)}+\frac{\alpha}{2}u_{j+2}^{(n+1)}} \\[2mm]
\quad=\displaystyle{\frac{\alpha}{2} u_{j-2}^{(n)}+(-\alpha+\beta)u_{j-1}^{(n)}+u_{j}^{(n)}+(\alpha-\beta)u_{j+1}^{(n+1)}-\frac{\alpha}{2}u_{j+2}^{(n+1)}},
\quad 1\le j\le J-2, \\[2mm]
u_J^{(n)}-Y_{1,C}^\xi *_d u_{J-1}^{(n)}+Y_{3,C}^\xi *_d u_{J-2}^{(n)}=-u_J^{(n-1)},\\[2mm]
u_{J}^{(n)}-2Y_{1,C}^\xi *_d u_{J-1}^{(n)}+Y_{2,C}^\xi *_d u_{J-2}^{(n)}-Y_{4,C}^\xi *_d
u_{J-4}^{(n)}=-u_{J}^{(n-1)}.
\end{cases}
\end{equation}
with the mesh ratios $\alpha=U_2\Delta t/(2({\Delta x})^3)$ and $\beta=U_1\Delta t/(4\Delta x)$. 

\begin{remark}
Concerning the implementation, it is usual to define the midpoint
unknown $v_j^{(n+1/2)}=(u_j^{(n+1)}+u_j^{(n)})/2$, 
with $v_j^{(-1/2)}=u_j^{(0)}$. 
In this case, \eqref{algo1} for the (R-CN) scheme reads
\begin{equation}\label{fulldisc4midpoint}
\begin{cases}
\displaystyle
Y_{2,R}^{\xi,(0)}v_0^{(n+1/2)}-Y_{1,R}^{\xi,(0)}v_1^{(n+1/2)}+v_2^{(n+1/2)}=-v_2^{(n-1/2)}
-\sum_{k=1}^nY_{2,R}^{\xi,(k)}v_0^{(n+1/2-k)}-\\
\hspace*{4cm}\displaystyle Y_{2,R}^{\xi,(n+1)}u_0^{(0)}/2+\sum_{k=1}^nY_{1,R}^{\xi,(k)}
v_1^{(n+1/2-k)}+Y_{1,R}^{\xi,(n+1)}u_1^{(0)}/2, \\[2mm]
   v_j^{(n+1/2)} +\alpha \left (
    v_{j+2}^{(n+1/2)}-3v_{j+1}^{(n+1/2)}+3v_j^{(n+1/2)}-v_{j-1}^{(n+1/2)}\right
  ) = u_j^n, \ 1\le j \le J-2.    \\[2mm]
\displaystyle -Y_{3,R}^{\xi,(0)}v_{J-2}^{(n+1/2)}+v_{J-1}^{(n+1/2)}=-v_{J-1}^{(n-1/2)}
+\sum_{k=1}^nY_{3,R}^{\xi,(k)}v_{J-2}^{(n+1/2-k)}+Y_{3,R}^{\xi,(n+1)}u_{J-2}^{(0)}/2, \\[2mm]
\displaystyle -Y_{4,R}^{\xi,(0)}v_{J-2}^{(n+1/2)}+v_{J}^{(n+1/2)}=-v_{J}^{(n-1/2)}
+\sum_{k=1}^nY_{4,R}^{\xi,(k)}v_{J-2}^{(n+1/2-k)}+Y_{4,R}^{\xi,(n+1)}u_{J-2}^{(0)}/2,
  \end{cases}
\end{equation}
For the (C-CN) scheme, \eqref{algo2} reads for $ \ 1\le j \le J-2$.
\begin{equation}\label{fulldisc4midpoint2}
\begin{cases}
\displaystyle
Y_{7,C}^{\xi,(0)}v_0^{(n+1/2)}-Y_{5,C}^{\xi,(0)}v_1^{(n+1/2)}+v_2^{(n+1/2)}=-v_2^{(n-1/2)}
-\sum_{k=1}^n Y_{7,C}^{\xi,(k)}v_0^{(n+1/2-k)}-\\
\hspace*{4cm}\displaystyle Y_{7,C}^{\xi,(n+1)}u_0^{(0)}/2+\sum_{k=1}^n Y_{5,C}^{\xi,(k)}
v_1^{(n+1/2-k)}+Y_{5,C}^{\xi,(n+1)}u_1^{(0)}/2, \\[2mm]
\displaystyle -Y_{8,C}^{\xi,(0)}v_0^{(n+1/2)}+Y_{6,C}^{\xi,(0)}v_2^{(n+1/2)}-2Y_{5,C}^{\xi,(0)}v_3^{(n+1/2)}+v_4^{(n+1/2)}\\
\displaystyle \quad=-v_4^{(n-1/2)}+\sum_{k=1}^n Y_{8,C}^{\xi,(k)}v_0^{(n+1/2-k)}+Y_{8,C}^{\xi,(n+1)}u_0^{(0)}/2-\sum_{k=1}^n Y_{6,C}^{\xi,(k)}v_2^{(n+1/2-k)}\\[2mm]
\displaystyle \qquad-Y_{6,C}^{\xi,(n+1)}u_2^{(0)}/2+2\sum_{k=1}^n Y_{5,C}^{\xi,(k)}v_3^{(n+1/2-k)}+Z_5^{(n+1)}u_3^{(0)}, \\[2mm]
\displaystyle  -\frac{\alpha}{2} v_{j-2}^{(n+1/2)} +(\alpha-\beta)v_{j-1}^{(n+1/2)}+v_{j}^{(n+1/2)}
+(-\alpha+\beta)v_{j+1}^{(n+1/2)}+\frac{\alpha}{2}v_{j+2}^{(n+1/2)}= u_j^n,\\[2mm]
\displaystyle
Y_{3,C}^{\xi,(0)}v_{J-2}^{(n+1/2)}-Y_{1,C}^{\xi,(0)}v_{J-1}^{(n+1/2)}+v_J^{(n+1/2)}=-v_J^{(n-1/2)}-\sum_{k=1}^n Y_{3,C}^{\xi,(k)}v_{J-2}^{(n+1/2-k)}-\\
\hspace*{4cm}\displaystyle Y_{3,C}^{\xi,(n+1)}u_{J-2}^{(0)}/2+\sum_{k=1}^n Y_{1,C}^{\xi,(k)}v_{J-1}^{(n+1/2-k)}+Y_{1,C}^{\xi,(n+1)}u_{J-1}^{(0)}/2, \\[2mm]
\displaystyle -Y_{4,C}^{\xi,(0)}v_{J-4}^{(n+1/2)}+Y_{2,C}^{\xi,(0)}v_{J-2}^{(n+1/2)}-2Y_{1,C}^{\xi,(0)}v_{J-1}^{(n+1/2)}+v_J^{(n+1/2)}\\
\displaystyle \quad=-v_J^{(n-1/2)}+\sum_{k=1}^n Y_{4,C}^{\xi,(k)}v_{J-4}^{(n+1/2-k)}+Y_{4,C}^{\xi,(n+1)}u_{J-4}^{(0)}/2-\sum_{k=1}^n Y_{2,C}^{\xi,(k)}v_{J-2}^{(n+1/2-k)}\\
\displaystyle \qquad-Y_{2,C}^{\xi,(n+1)}u_{J-2}^{(0)}/2+2\sum_{k=1}^n Y_{1,C}^{\xi,(k)}v_{J-1}^{(n+1/2-k)}+Y_{1,C}^{\xi,(n+1)}u_{J-1}^{(0)}.
  \end{cases}
\end{equation}
Solving \eqref{fulldisc4midpoint} or  \eqref{fulldisc4midpoint2}, 
we recover $u_j^{(n+1)}$ by $u_j^{(n+1)}=2v_{j}^{(n+1/2)}-u_j^{(n)}$.
\end{remark}

\section{The Sum-of-Exponentials Approach}\label{s:expo}
An ad-hoc implementation of the discrete convolutions of the form
\begin{equation*}
\sum_{k=1}^{n}X_m^{(k)}u_j^{(n-k)}
\end{equation*}
with convolution coefficients $X_m^{(n)}$
has still one disadvantage. 
The boundary conditions are non--local in time
(and space for higher dimensions) and therefore computations are too expensive. 
As a remedy, to get rid of the time non-locality, we use as in \cite{AES03} the sum of
exponentials ansatz, i.e.\ to approximate the convolution coefficients $X_m^{(n)}$
by a finite sum (say $L_m$ terms) of exponentials that \emph{decay} 
with respect to time.
This approach allows for a fast 
(approximate) evaluation of the discrete convolution 
since the convolution can now be evaluated with a simple recurrence
formula for $L_m$ auxiliary terms and the numerical effort per time step now stays constant.


\subsection{The Exponential Approximation}
To do so we will follow the ideas of \cite{AES03} and approximate the
coefficients of a sequence $X_m^{(n)}$ 
by the following {\em sum-of-exponentials ansatz}
\begin{equation}\label{cc8}
 X_m^{(n)}\approx \widetilde{X}_m^{(n)}:=
 \begin{cases}
 X_m^{(n)}, &n=0,\dots,\nu_m-1,\\
 \sum_{l=1}^{L_m} b_{m,l} q_{m,l}^{-n},&n=\nu_m, \nu_m+1,\dots\,,
 \end{cases}
\end{equation}
where $L_m\in\Z$, $\nu_m\ge 0$ are given integer parameters, e.g.\ $L_m=20$, $\nu_m=2$,
that have to be chosen appropriately to guarantee good
approximation properties of $\widetilde{X}_m^{(n)}$. In the
following, $X_m^{(n)}$ has to be seen as $Y_{m,R}^{\xi,(n)}$ or $Y_{m,C}^{\xi,(n)}$
respectively for (R-CN) and (C-CN) schemes.

In \cite{AES03} the authors presented a deterministic method 
of choosing such an optimal approximation, i.e.\
finding  the set $\{b_{m,l}, q_{m,l}\}$ for fixed $L_m$ and $\nu_m$.

The ``split'' definition of $\widetilde{X}_m^{(n)}$ in \eqref{cc8} is motivated
by the fact that the implementation of the discrete TBCs 
involves a convolution sum with $k$ ranging only from 1 to $k=n$.
Since the first coefficient $X_m^{(0)}$ does not appear in this convolution,
it makes no sense to include it in our sum-of-exponential approximation,
which aims at simplifying the evaluation of the convolution.
Hence, one may choose $\nu_m=1$ in \eqref{cc8}.
We observe numerically that the two first coefficients have a larger magnitude
compared to the other ones. This suggests even to exclude $Z_m^{(1)}$ from this approximation
and to choose $\nu_m=2$ in \eqref{cc8}.
We use this choice in our numerical implementation.


Let us fix $L_m$ and consider the formal power series:
\begin{equation}\label{cc9}
  g_m(x):= s^{(\nu_m)} + s^{(\nu_m+1)}x + s^{(\nu_m+2)}x^2 + \dots,\quad |x|\le1.
\end{equation}
If there exists the $[L_m-1|L_m]$ \emph{Pad\'{e} approximation}
\begin{equation*}
 \tilde{g}_m(x):=\frac{P_{L_m-1}(x)}{Q_{L_m}(x)}
\end{equation*}
of \eqref{cc9}, then its Taylor series
\begin{equation*}
 \tilde{g}_m(x)= \widetilde{X}_m^{(\nu_m)} + \widetilde{X}_m^{(\nu_m+1)}x 
 + \widetilde{X}_m^{(\nu_m+2)}x^2 + \dots
\end{equation*}
satisfies the conditions
\begin{equation}\label{cc11}
 \widetilde{X}_m^{(n)} = X_m^{(n)}, \qquad n=\nu_m,\nu_m+1,\dots,2L_m+\nu_m-1,
\end{equation}
due to the definition of the Pad\'{e} approximation rule.

\begin{theorem}[\cite{AES03}]\label{tcc1}
Let $Q_{L_m}(x)$ have $L_m$ simple roots
$q_{m,l}$ with $|q_{m,l}| > 1, \quad l=1,\dots,L_m$. Then
\begin{equation}\label{cc12}
  \widetilde{X}_m^{(n)} = \sum_{l=1}^{L_m} b_{m,l}^{\phantom{n}} q_{m,l}^{-n}, 
  \qquad n=\nu_m,\nu_m+1,\dots\,,
\end{equation}
where
\begin{equation}\label{cc13}
 b_{m,l} := -\frac{P_{L_m-1}(q_{m,l})}{Q_{L_m}'(q_{m,l})}\,q_{m,l}^{\phantom{i}} \neq 0, \qquad
 l=1,\dots,L_m.
\end{equation}
\end{theorem}

Evidently, the approximation of the convolution coefficients $X_m^{(n)}$ by the
representation \eqref{cc8} using a $[L_m-1|L_m]$ Pad\'{e} approximant to
\eqref{cc9} behaves as follows: the first $2L_m$ coefficients are reproduced
exactly, see \eqref{cc11}. 
However, the asymptotics of $X_m^{(n)}$ and
$\widetilde{X}_n^{(n)}$ (as $n\to\infty$) differ strongly -- 
algebraic versus exponential decay. 


\subsection{Fast Evaluation of the Discrete Convolution}
Let us consider the approximation \eqref{cc8} with $\nu_m=2$ for the discrete
convolution kernel appearing in the discrete TBCs. 
With these ``exponential'' coefficients the \emph{approximated convolution}
\begin{equation}\label{cc18}
 \tilde{C}_{m,j}^{(n)} := \sum_{k=2}^{n}\widetilde{X}_m^{(k)}u_{j}^{(n-k)},
 \qquad\widetilde{X}_m^{(n)} = \sum_{l=1}^{L_m} b_{m,l}^{\phantom{n}} q_{m,l}^{-n},
 \qquad|q_l|>1,
\end{equation}
of the discrete function $u_{j}^{(n-k)}$, $k=1,2,\dots$ with the
coefficients $\widetilde{X}_m^{(n)}$ can be calculated by recurrence
formulas, and this will reduce the numerical effort significantly.

A straightforward calculation (\cite{AES03}) yields (for $\nu_m=2$):
\begin{equation}\label{cc20}
 \tilde{C}_{m,j}^{(n)}(\{u_j^{(n)})_n\}=\sum_{l=1}^{L_m} \tilde{C}_{m,j,l}^{(n)},\qquad n\ge2,
\end{equation}
where
\begin{equation*}
\tilde{C}_{m,j,l}^{(2)} \equiv 0,
\end{equation*}
\begin{equation}\label{cc21}
 \tilde{C}_{m,j,l}^{(n)}=q_{m,l}^{-1} \tilde{C}_{m,j,l}^{(n-1)}
                         +b_{m,l}^{\phantom{n}} q_{m,l}^{-1} u_{j}^{(n-2)},
\end{equation}
$n=2,3,\dots$, $l=1,\dots,L_m$.  

In order to use this fast evaluation procedure in our implementation
point of view, we must transform it before to use it for midpoint
$v_j^{(n+1/2)}$ unknown. It is easy to see that the second relation of
\eqref{cc21} can be transformed as
\begin{equation}
  \label{cc22}
  \tilde{C}_{m,j,l}^{(n)}\{(u_j^{(n)})_n\}=b_{m,l}\sum_{k=2}^{n-1} q_{m,l}^{-k}u_j^{(n-k)}.
\end{equation}
For example, let us consider the last TBC of \eqref{algo1}
\begin{equation*}
u_{J}^{(n+1)}-Y_{4,R}^{\xi} *_d u_{J-2}^{(n+1)} =-u_J^n.
\end{equation*}
We have to see $Y_{4,R}^{\xi} *_d u_{J-2}^{(n+1)}$ as
\begin{equation*}
Y_{4,R}^{\xi,(0)}u_{J-2}^{(n+1)}+Y_{4,R}^{\xi,(1)}u_{J-2}^{(n)}+\tilde{C}_{4,J-2}^{(n+1)}\{(u_{J-2}^{(n+1)})_n\}. 
\end{equation*}
In
order to get an equation for $v_{J-2}^{(n+1/2)}$, we write the previous
relation at discrete time level $n$ and average the equations. 
We therefore get
\begin{multline}
  \label{cc23}
v_{J}^{(n+1/2)}-Y_{4,R}^{\xi,(0)}v_{J-2}^{(n+1/2)}-Y_{4,R}^{\xi,(1)}v_{J-2}^{(n-1/2)}\\
=-v_J^{(n-1/2)}
+\frac{1}{2}\left ( \tilde{C}_{4,J-2}^{(n+1)}\{(u_{J-2}^{(n+1)})_n\} +
  \tilde{C}_{4,J-2}^{(n)}\{(u_{J-2}^{(n)})_n\}\right).
\end{multline}
Thanks to \eqref{cc20} and \eqref{cc22}, we obtain
\begin{multline}
  \label{cc24}
v_{J}^{(n+1/2)}-Y_{4,R}^{\xi,(0)}v_{J-2}^{(n+1/2)}-Y_{4,R}^{\xi,(1)}v_{J-2}^{(n-1/2)}\\
=-v_J^{(n-1/2)}
+  \tilde{C}_{4,J-2}^{(n)}\{(v_{J-2}^{(n+1/2)})_n\}+\sum_{l=1}^{L_m}b_{4,l}q_{4,l}^{-n}v_0^{(1/2)}.
\end{multline}
These equation then replaces the last one of \eqref{fulldisc4midpoint}.

These computations can be easily transferred to other convolutions
appearing in other TBCs.

\section{Numerical Results}\label{num}

In this section we first present the numerical procedure used to compute 
the inverse $\mathcal{Z}$-transforms required for the discrete absorbing 
boundary conditions. Then we consider two different examples for which 
we give some numerical results. 
The first example can be considered for either (R-CN) or (C-CN) scheme since 
$U_1=0$, $U_2=1$.
The second one can only be considered for the (C-CN) scheme since in this 
ones $U_1=U_2=1$.

\subsection{Numerical procedure for the inverse $\mathcal{Z}$-transform}
In this section, we recall for a self-contained presentation the numerical
procedure presented in \cite{Zi02003} to compute the inverse $\mathcal{Z}$-transform. 

Let us recall that if we consider the sequence $(u_n)_{n\in\N}$, then its $\mathcal{Z}$-transform reads
$$
U(z)=\sum_{k=0}^\infty u_k z^{-k},\quad |z|>R.
$$
Assuming we know $U(z)$, we can recover the value of the sequence $u_n$ thanks to the relation
$$
u_n=\frac{1}{2i\pi} \oint_{S_r} U(z) z^{n-1} dz, \quad r>R,
$$
where $S_r$ denotes any circle of radius $r>R$. Performing the change of variable $z=r e^{i\varphi}$, we obtain
$$u_n=\frac{r^n}{2\pi} \int_0^{2\pi} U(r e^{i\varphi}) e^{in\varphi}\, d\varphi.
$$
Discretizing the angular variable $\varphi$ with $N$ nodes, we obtain the approximation
$$
u_n\approx \frac{r^n}{N} \sum_{k=0}^{N-1} U\left (r e^{i\frac{2\pi}{N}k}\right ) e^{i\frac{2\pi}{N}kn}.
$$
Let us now consider the finite sequence $\displaystyle\big(f_n\big)_{n=0}^{N-1}$. Its discrete Fourier transform is given by
$$
F_k=\mathscr{F}\{f_n\}(k)=\sum_{n=0}^{N-1}f_n\omega_N^{-nk}
$$
and we recover the value of $f_n$ through
$$
f_n=\mathscr{F}^{-1}\{F_k\}(n)=\frac{1}{N}\sum_{k=0}^{N-1} F_k\omega_N^{nk}
$$
where $\omega_N=e^{i\frac{2\pi}{N}}$. Therefore, if we define $U_k=U(r \omega_N^k)$, we have
$$
u_n\approx \frac{r^n}{N} \sum_{k=0}^{N-1} U_k \omega_N^k = r^n\mathscr{F}^{-1}\{U_k\}(n), \quad 0\leq n <N.
$$
Thus, in order to obtain approximation of $u_n$, we have to multiply the inverse discrete Fourier transform of $U(z)$ evaluated at nodes $z_k=r e^{2i\pi k/N}$ by $r^k$. The inverse discrete Fourier transform if easily obtain thanks for inverse fast Fourier transform.

Note that the choice of the inversion radius $r$ is crucial to guarantee the 
good approximation of the inverse $\mathcal{Z}$-transform and then the convergence of the 
numerical scheme as presented in Figure~\ref{error_wrt_r}. 
Note that in \cite{Zi02003}
the best choice of $r$ seems to be $1.02$ 
while in our case it seems to be $1.001$. For a concise discussion on the choice of $r$ we refer
the reader to \cite{Zi02003}.

\subsection{Numerical Example 1}\label{num1}
Let us first consider the example from Zheng, Wen and Han \cite{ZhWeHa2008} 
which is concerned with the following equation ($U_1=0,~U_2=1$):
\begin{align}
u_t+u_{xxx}&=0, \qquad x \in \R, \label{ex1_1} \\
u(0,x)&=e^{-x^2}, \qquad x \in \R, \label{ex1_2} \\
u &\to 0,\qquad |x| \to \infty.  \label{ex1_3}
\end{align}
The fundamental solution of equation \eqref{ex1_1} is \cite{ZhWeHa2008}
\begin{equation*}
E(t,x)=\frac{1}{\sqrt[3]{3t}}\Ai\left(\frac{x}{\sqrt[3]{3t}}\right),
\end{equation*}
where $\Ai(\cdot)$ is the Airy function. 
The exact solution of \eqref{ex1_1}-\eqref{ex1_3} can be written in terms of $E(t,x)$ as
\begin{equation*}
u_{\rm exact}(t,x)=E(t,x) * e^{-x^2}, 
\end{equation*}
where $*$ denotes the convolution product on the whole real axis.

\begin{figure}[ht]
\begin{center}
\begin{tabular}{cc}
\includegraphics[width=0.48\textwidth]{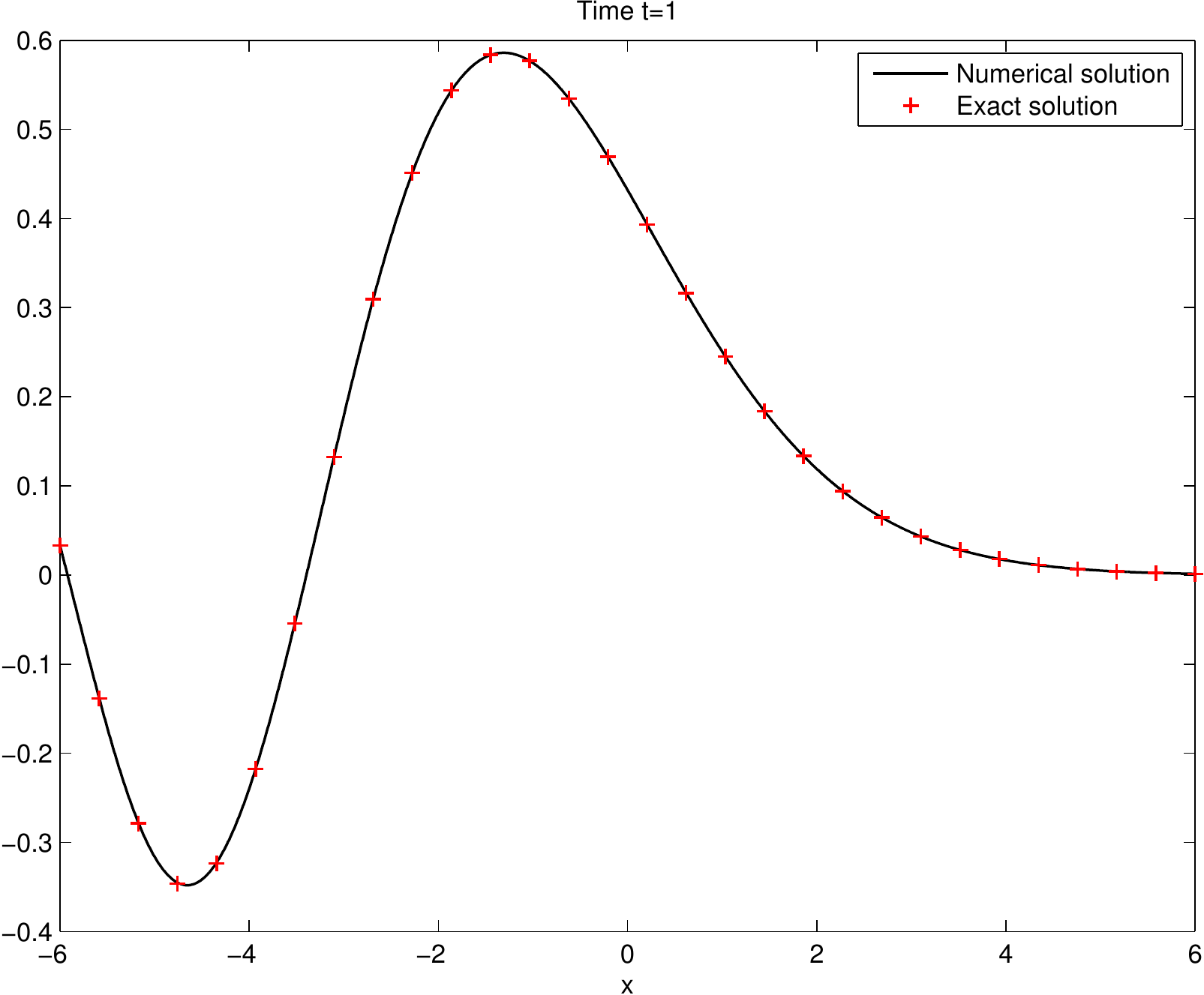} & 
\includegraphics[width=0.48\textwidth]{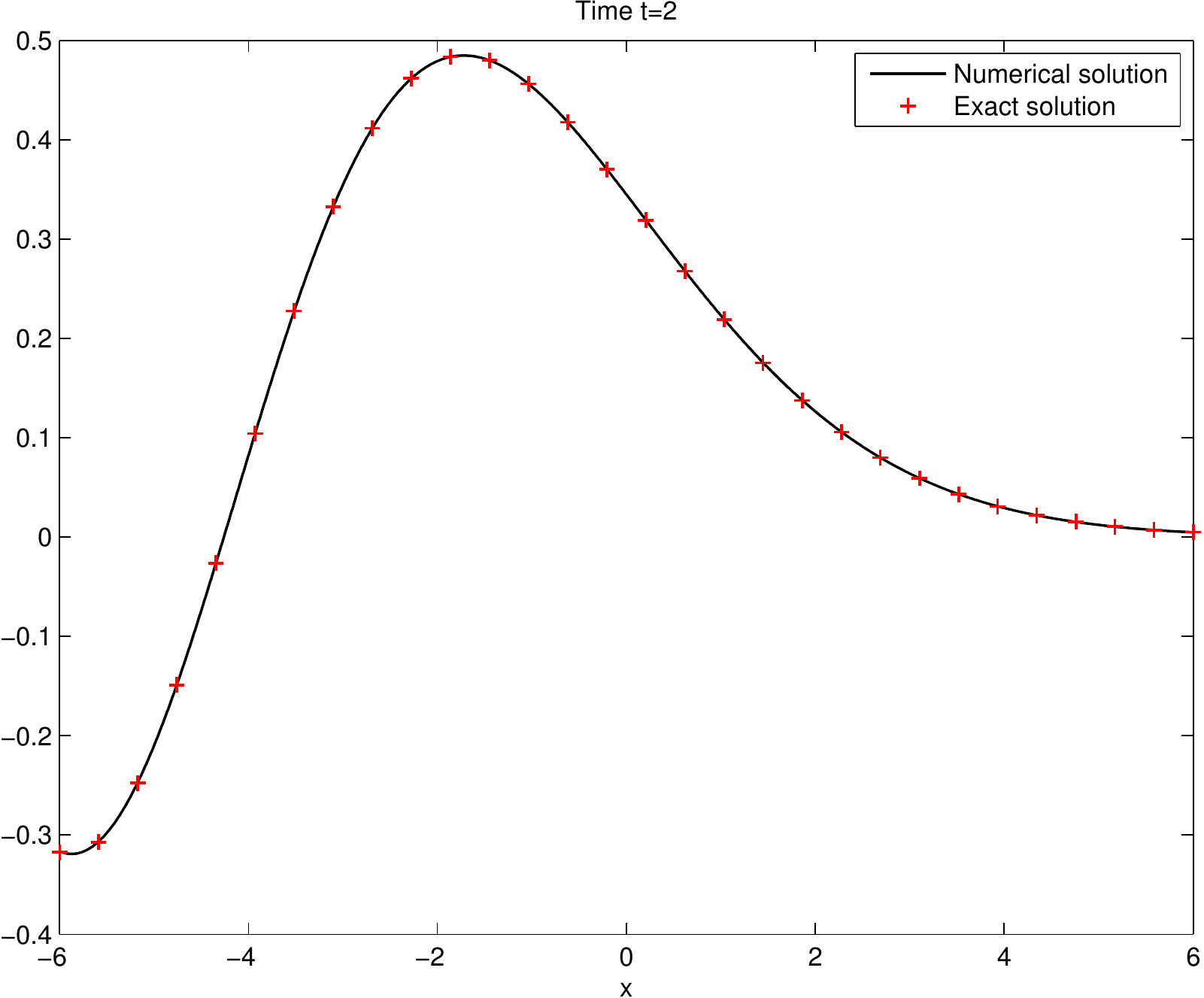} \\
\includegraphics[width=0.48\textwidth]{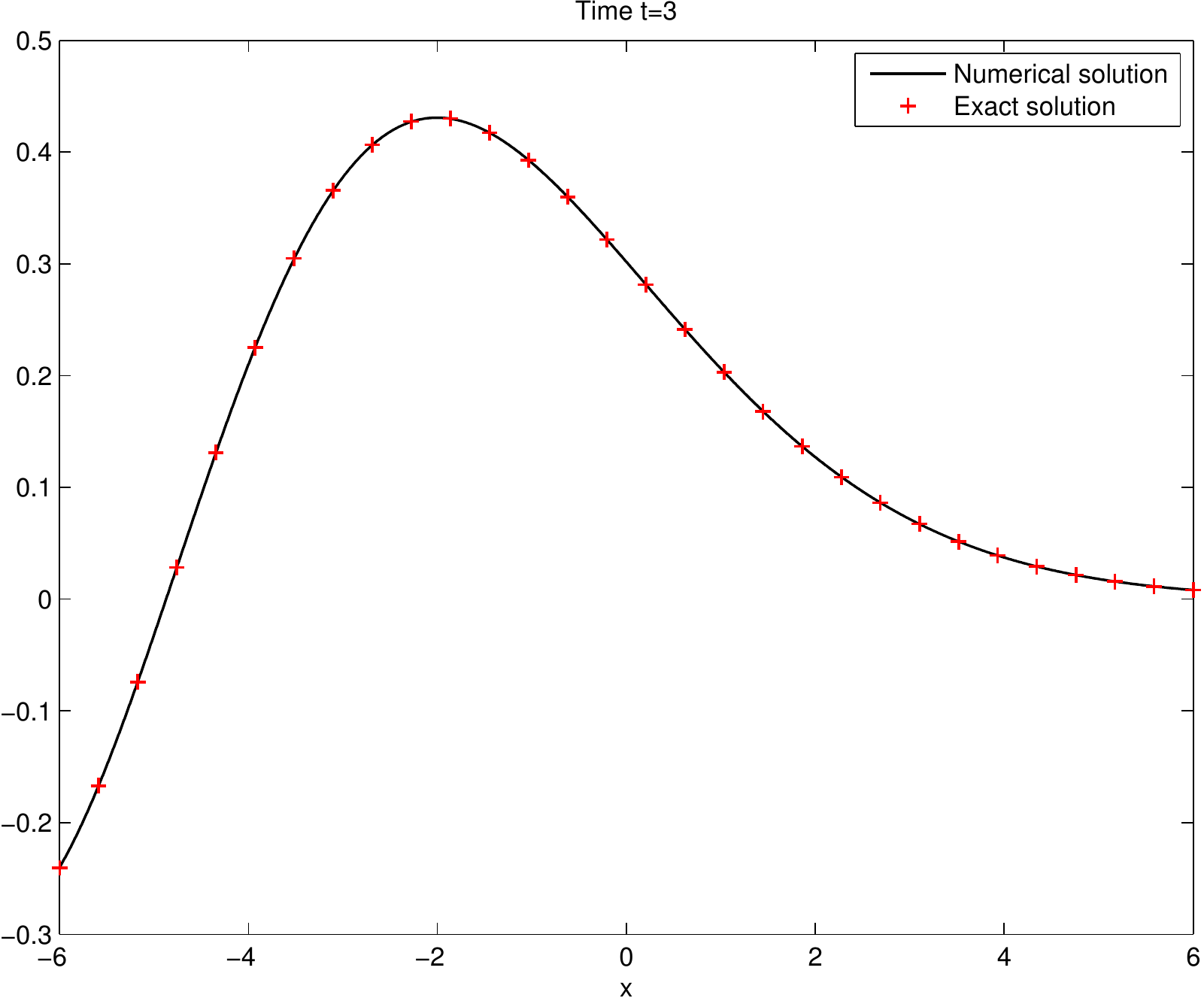} & 
\includegraphics[width=0.48\textwidth]{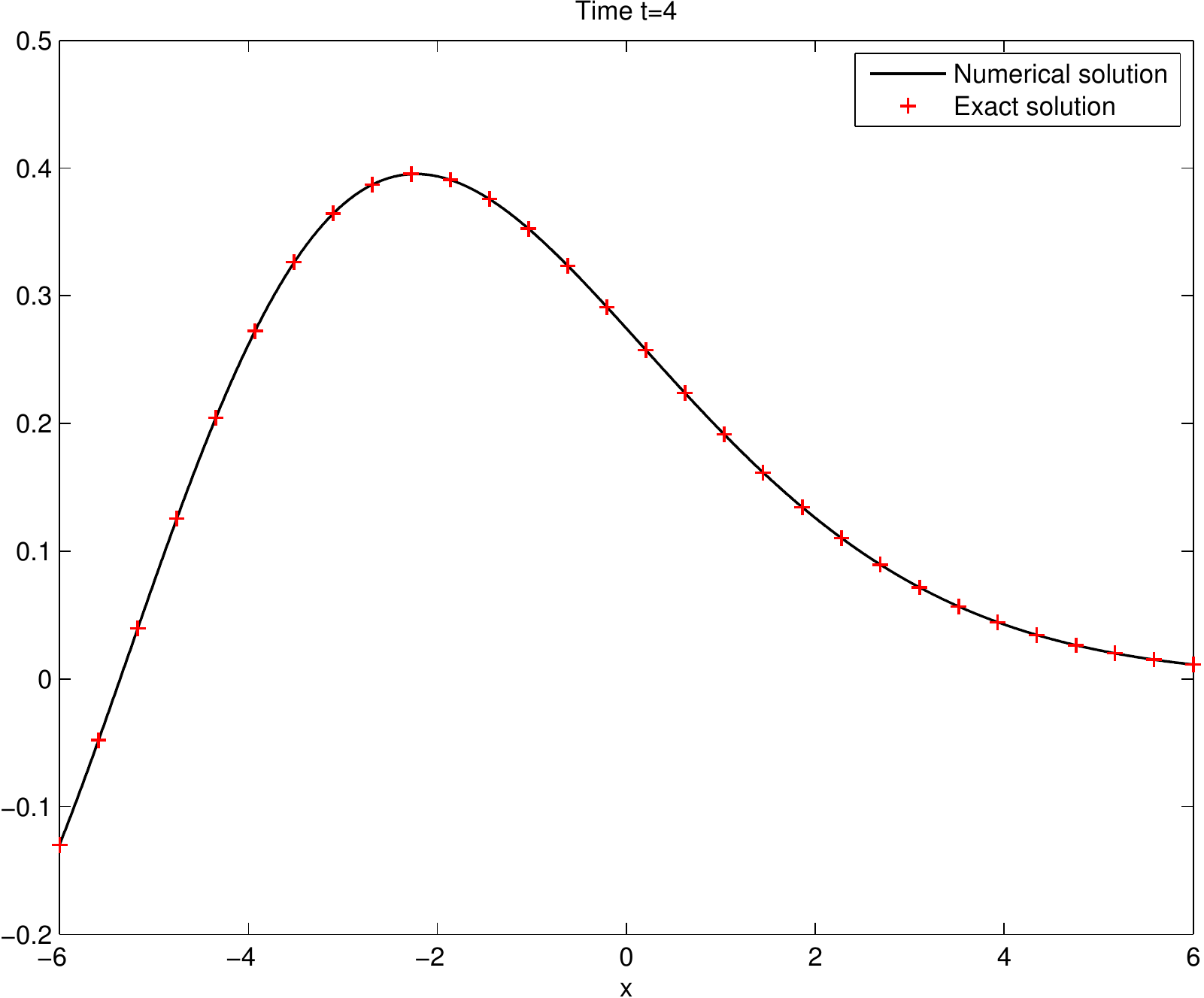}
\end{tabular}
\end{center}
\caption{Numerical and exact solutions at times $t=1$, $t=2$, 
$t=3$ and final time $T=4$ for the first example with $\Delta t=4/2560$, 
$\Delta x=12/5000$ and $r=1.001$.}
\label{compZhengRCN}
\end{figure}

We present in Figure~\ref{compZhengRCN} the exact solution and the approximate
solution obtained with (R-CN) scheme for $\Delta t=4/2560$, $\Delta x=12/5000$ and $r=1.001$ 
at different times $t=1,2,3,4$. We see that the (R-CN) solution is a very good approximation of the 
exact solution all along the time. No unphysical reflections can be seen at the boundaries. 
The same can be obtained by using the (C-CN) scheme. 

\begin{figure}[ht]
\begin{center}
\begin{tabular}{cc}
\includegraphics[width=0.48\textwidth]{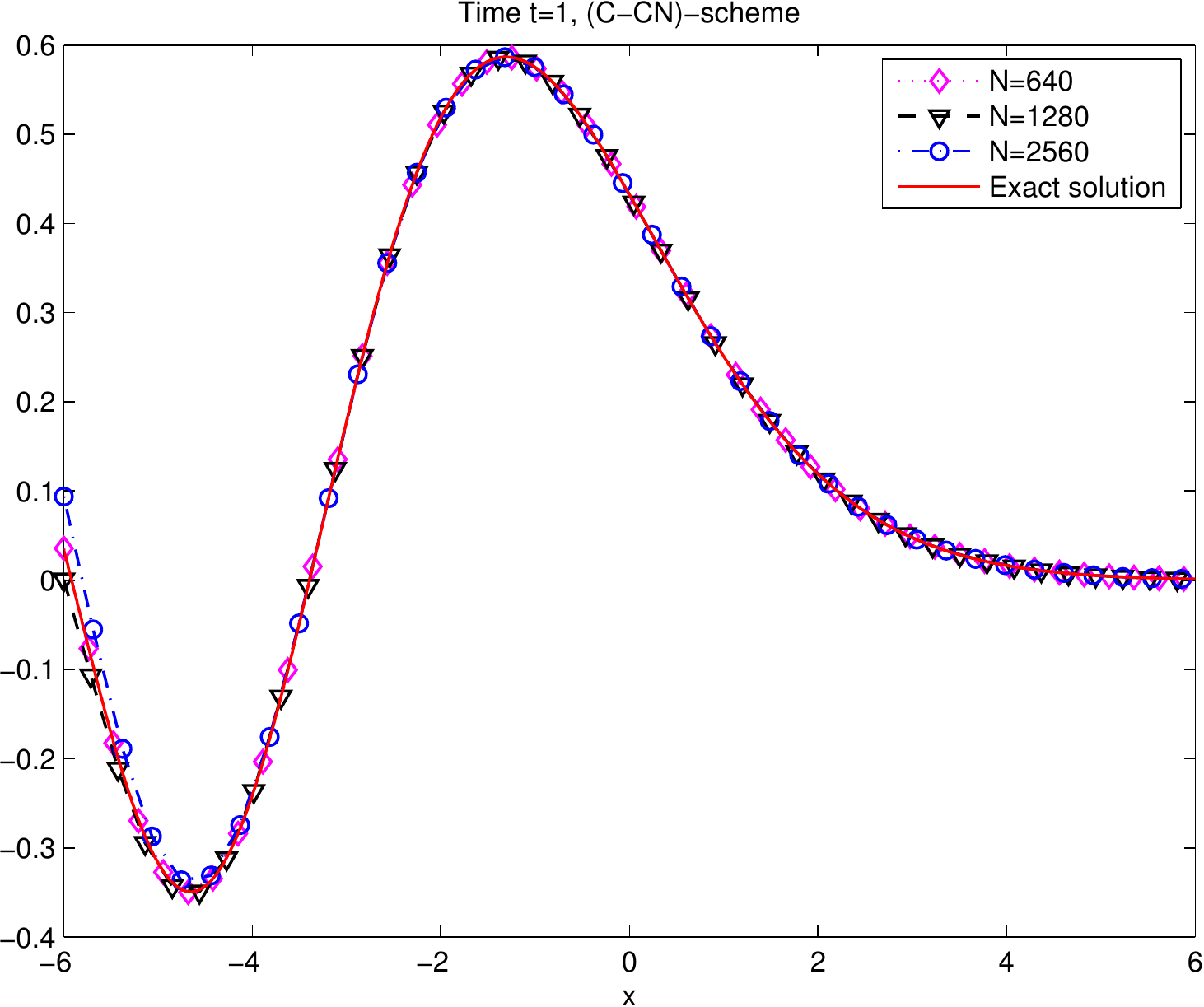} & 
\includegraphics[width=0.48\textwidth]{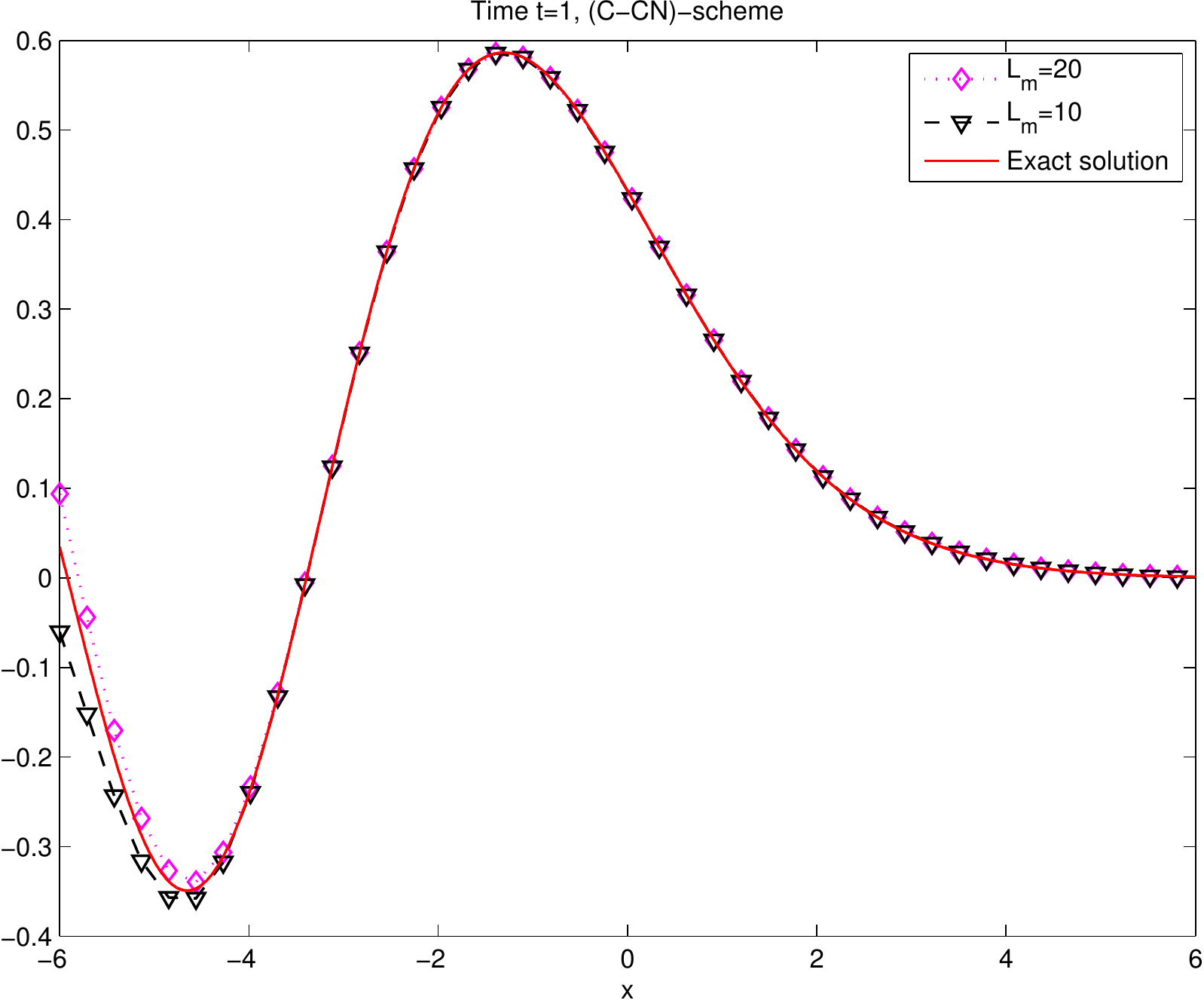} \\
\includegraphics[width=0.48\textwidth]{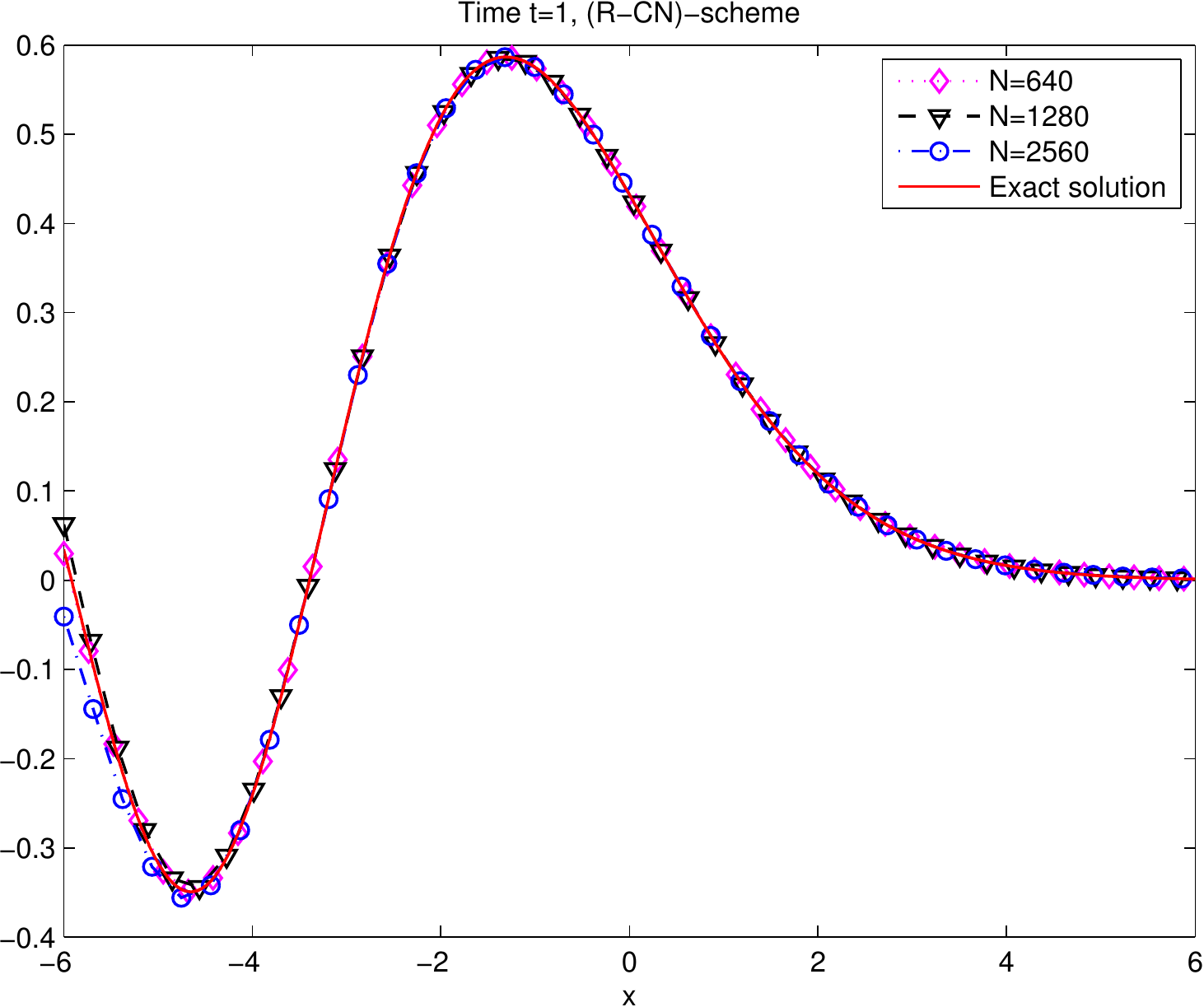} & 
\includegraphics[width=0.48\textwidth]{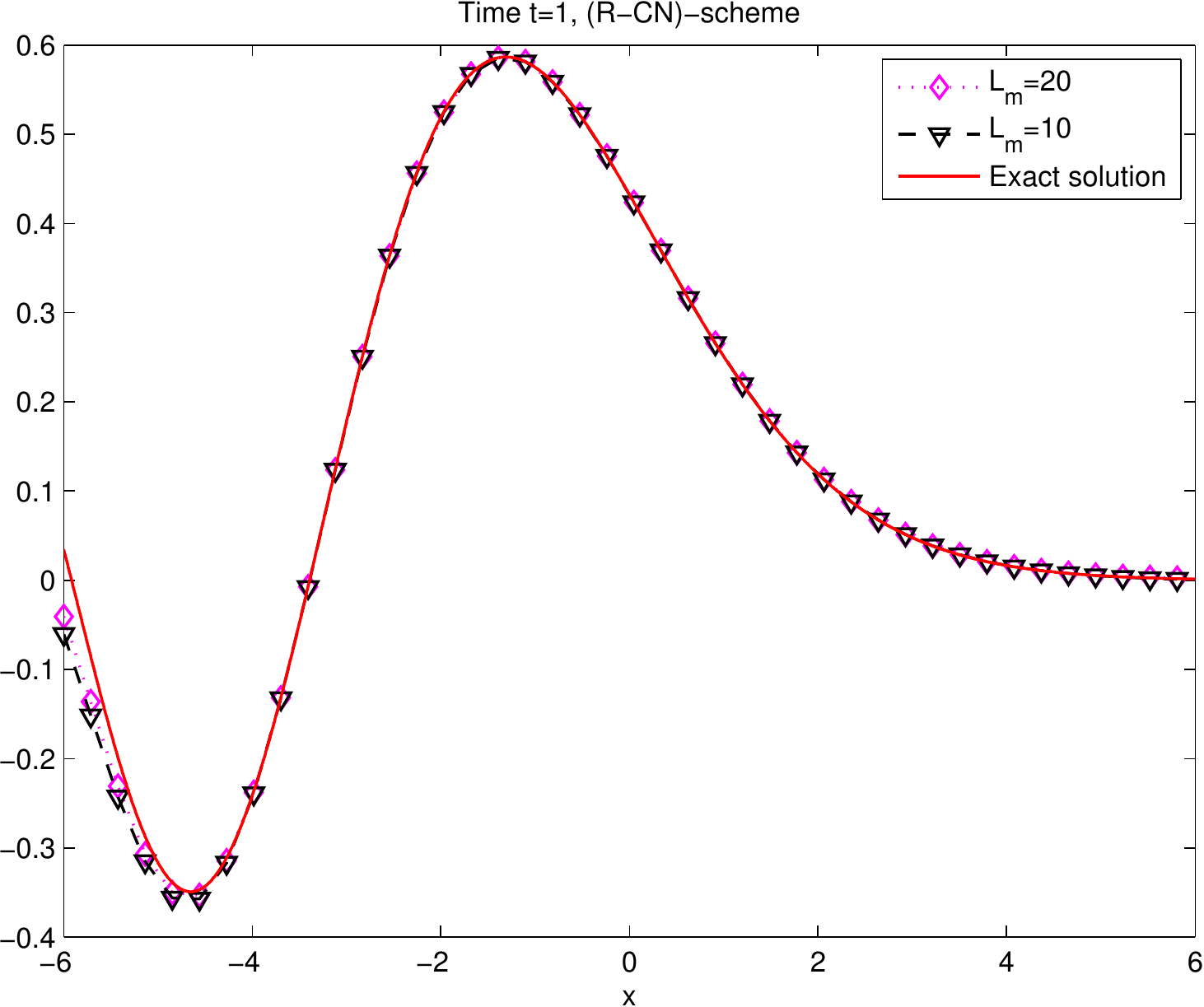}
\end{tabular}
\end{center}
\caption{Comparison at time $T=1$ between the exact solution and the approximate solution obtained with 
the sum of exponential approach for the (C-CN) scheme (top figures) and the (R-CN)-scheme
(bottom figures). The left figures are obtained with a fixed $L_m=20$ and various $N=640, 
1280, 2560$ and the right figures with a fixed $N=2560$ and various $L_m=10, 20$ for 
(C-CN) and (R-CN)-schemes.}
\label{sumexp}
\end{figure}

We present in Figure~\ref{sumexp} a comparison at time $T=1$ between the exact solution 
and the approximate solution obtained with the sum of exponential approach either for 
various values of $N$ ($N=640, 1280, 2560$) and a fixed value of $L_m$ ($L_m=20$) or 
for a fixed value  of $N$ ($N=2560$) and various values of $L_m$ ($L_m=10, 20$). In each case $\Delta x=12/5000$ and $r=1.001$.
We observe that the accuracy of the approximate solution depend on $N$ for a fixed $L_m$ and on 
$L_m$ for a fixed $N$.

Let us define as $e^{(n)}$ the {\em relative $\ell^2$-error} at time $t=n\Delta t$ given by:
\begin{equation*}
e^{(n)}=\left\|u^{(n)}_{\rm exact}-u^{(n)}_{\rm num} \right\|_2/\left\| u^{(n)}_{\rm exact}\right\|_2,
\end{equation*}
where we use trapezoidal rule to compute the $\ell^2$-norm. 
Note that here $u_{\rm num}$ stands for the numerical solution computed with either (R-CN) or (C-CN) scheme. 
We decided to compute from  $e^{(n)}$ two error functions; first  
the maximum in time and secondly the $\ell^2$-error in time: 
\begin{equation*}
rel.ErrTm=\max_{0<n<N} \left(e^{(n)}\right),\qquad
rel.ErrL2=\left(\Delta t \sum_{n=1}^N (e^{(n)})^2\right)^{1/2}.
\end{equation*}

\begin{figure}[ht]
\begin{center}
\begin{tabular}{cc}
\includegraphics[width=0.48\textwidth]{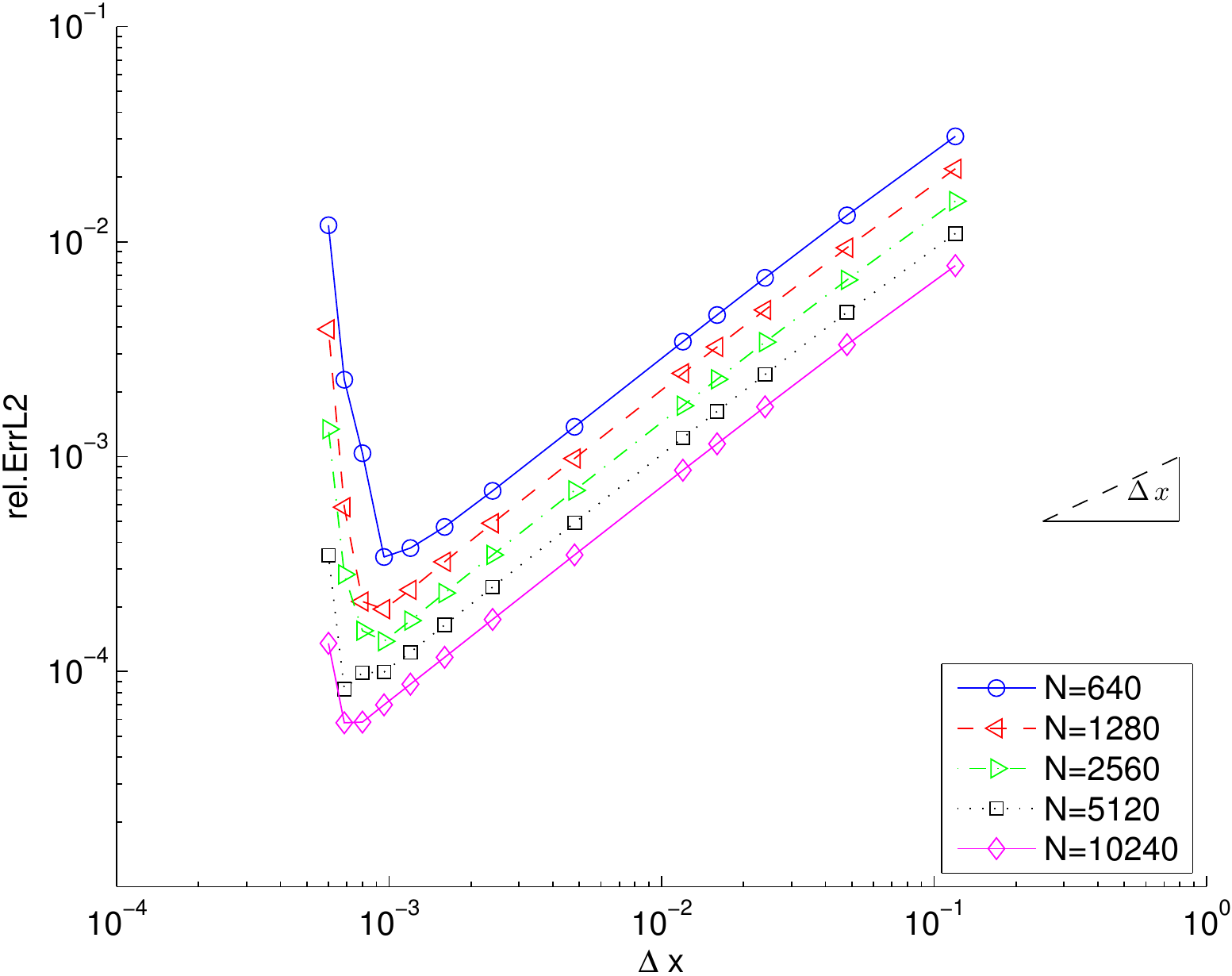} & 
\includegraphics[width=0.48\textwidth]{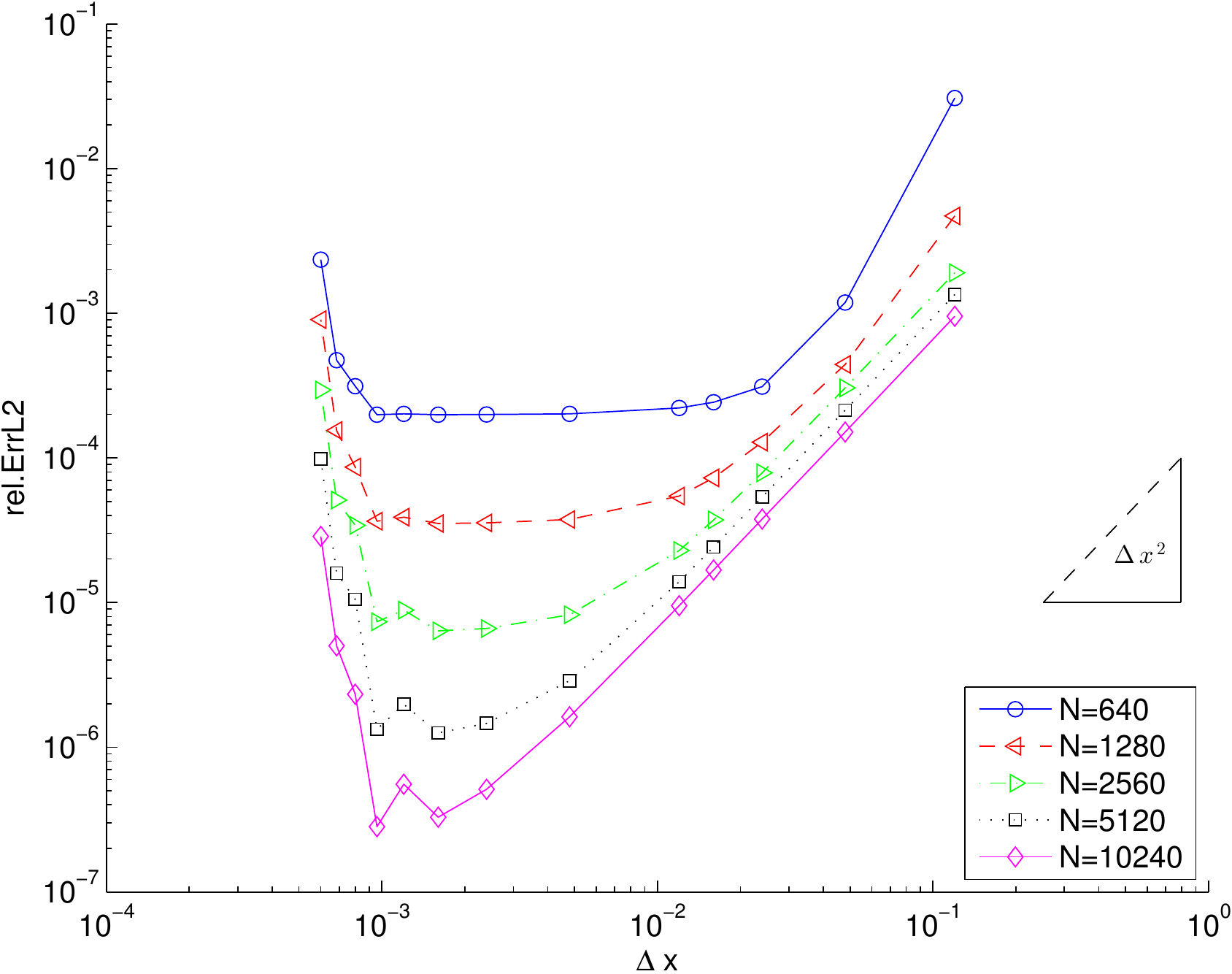} \\
\includegraphics[width=0.48\textwidth]{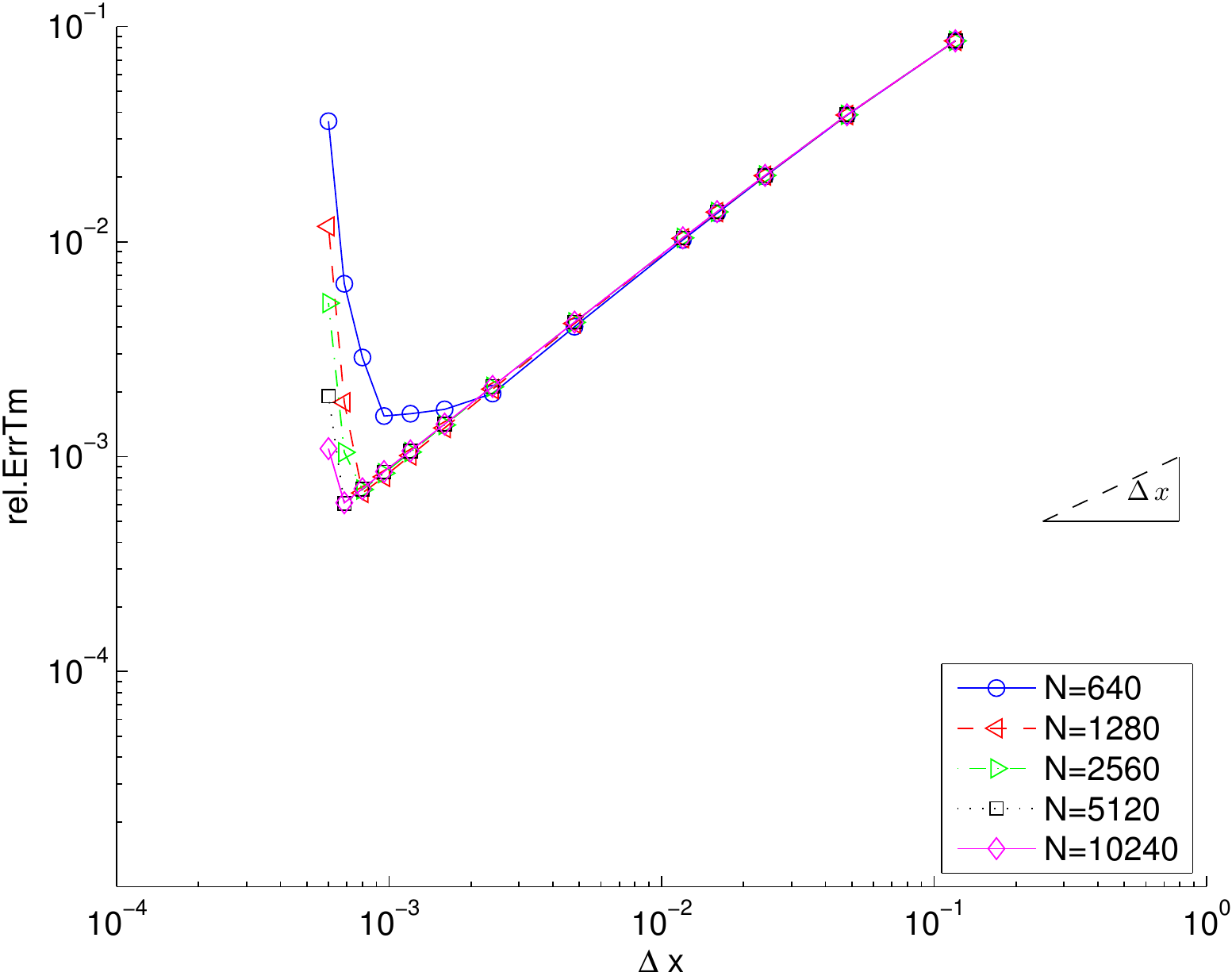} & 
\includegraphics[width=0.48\textwidth]{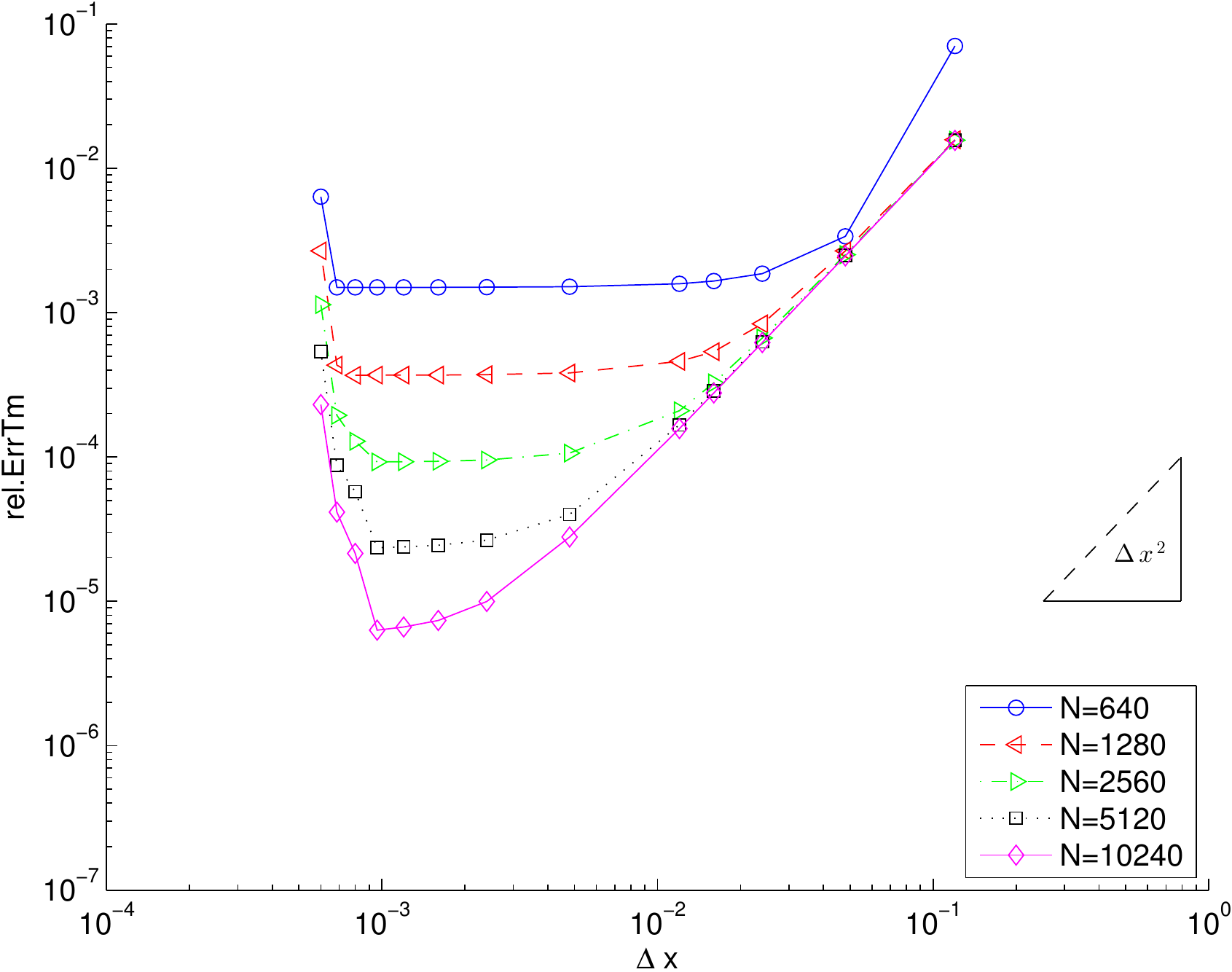}
\end{tabular}
\end{center}
\caption{ Relative errors with respect to $\Delta x$ at time $T=4$ for the (R-CN) scheme 
(figures on left) and the (C-CN)-scheme (figures on right) and for different values of $N$.}
\label{example1_error_dx}
\end{figure}

\begin{figure}[ht]
\begin{center}
\begin{tabular}{cc}
\includegraphics[width=0.48\textwidth]{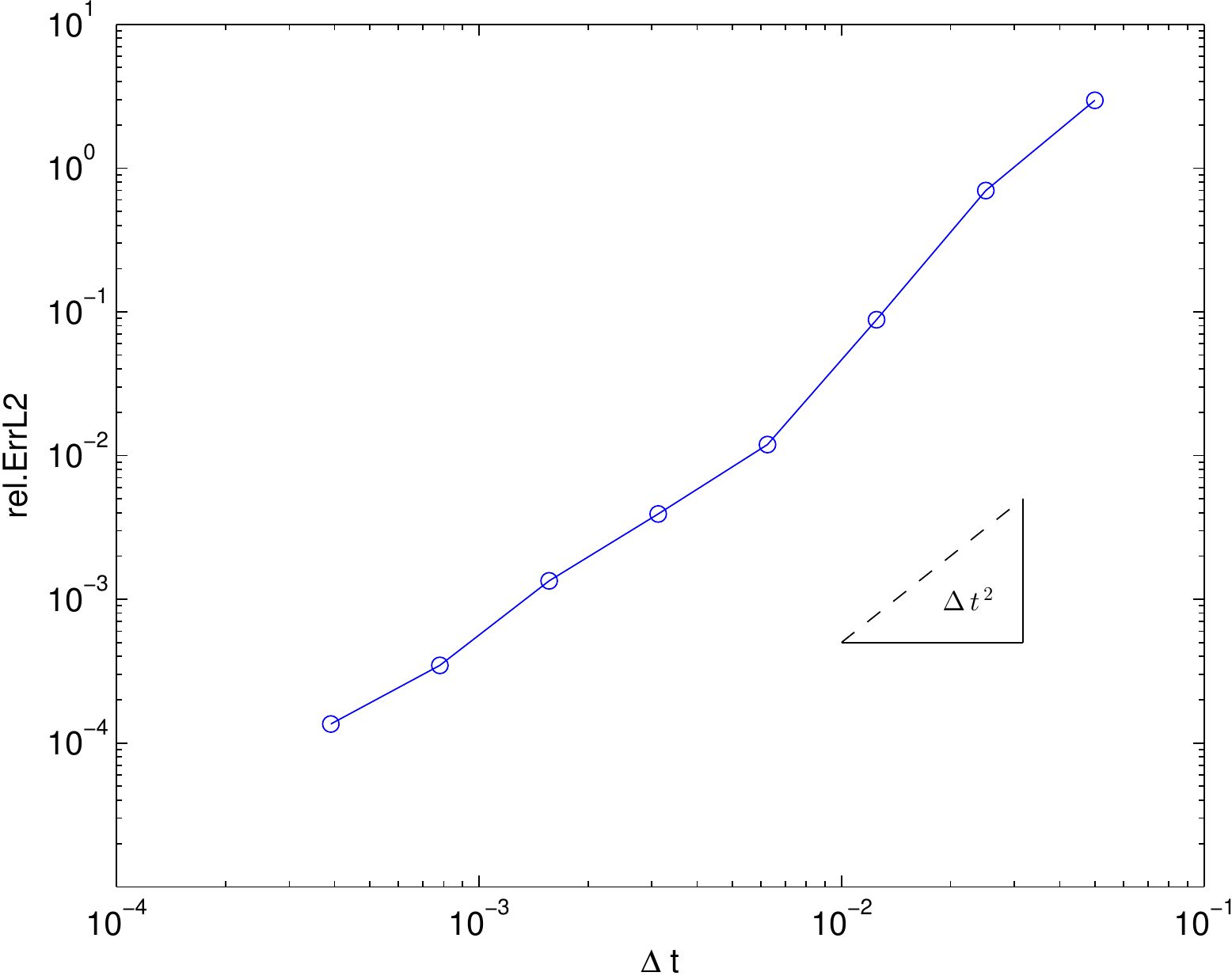} & 
\includegraphics[width=0.48\textwidth]{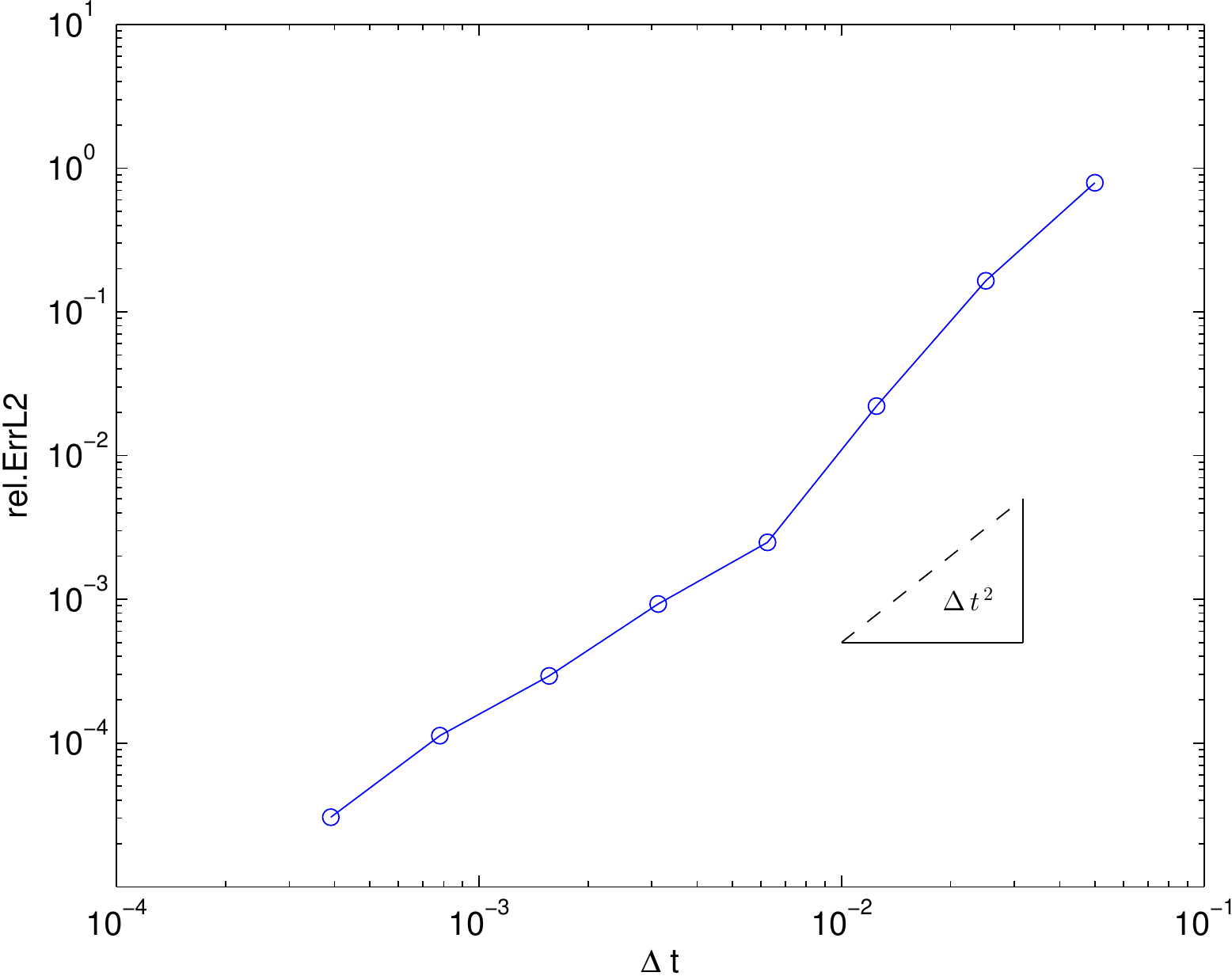} \\
\includegraphics[width=0.48\textwidth]{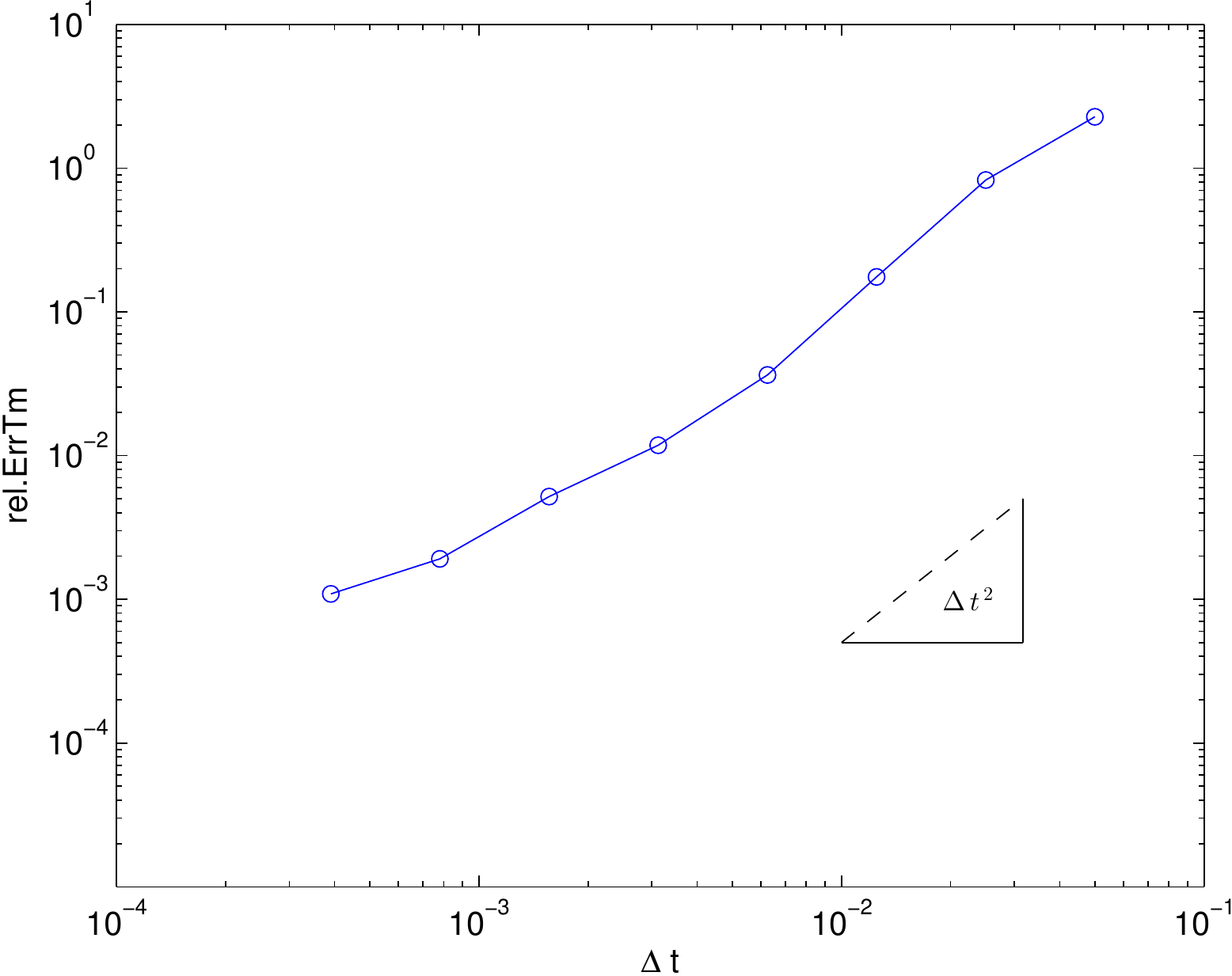} & 
\includegraphics[width=0.48\textwidth]{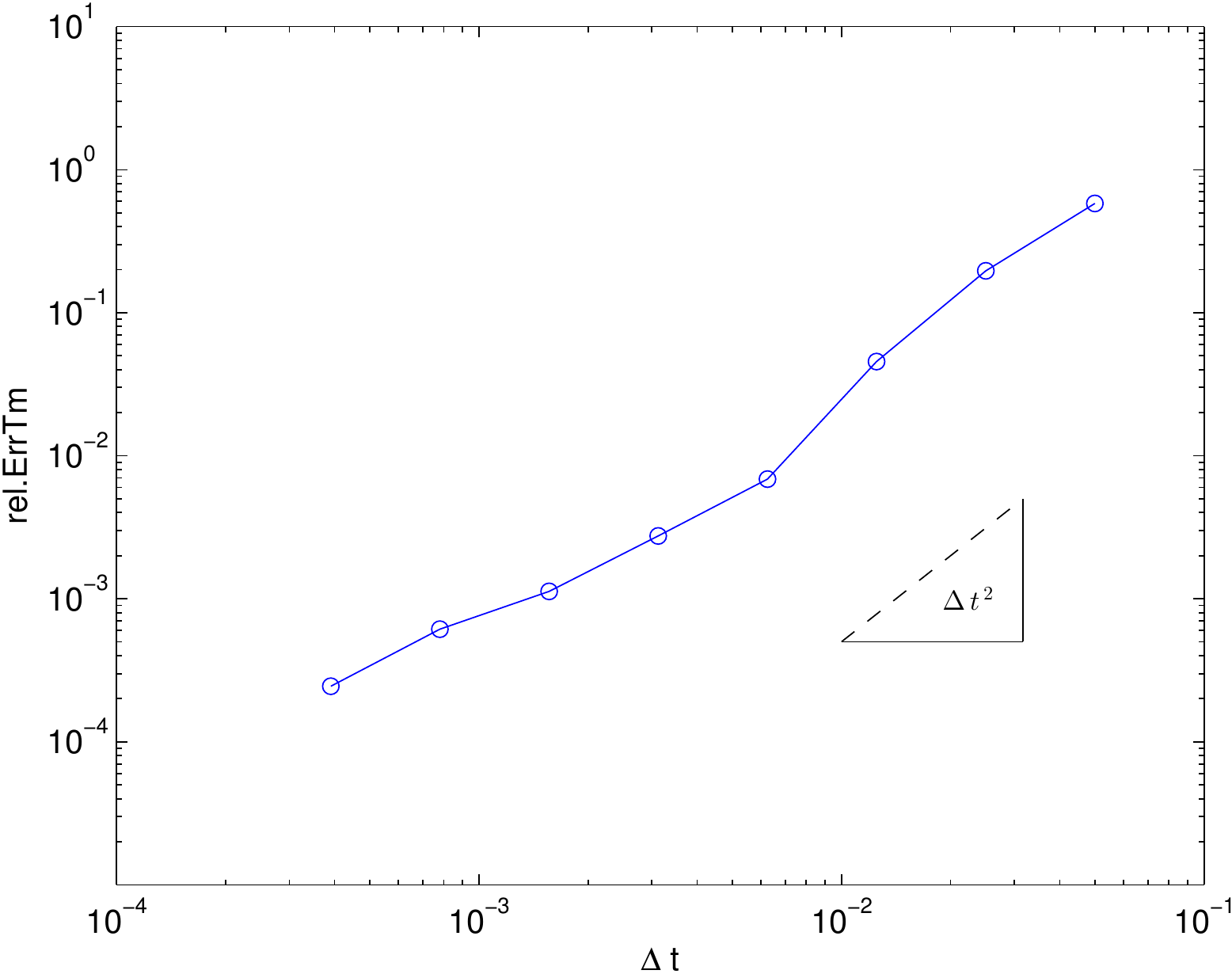}
\end{tabular}
\end{center}
\caption{ Relative errors with respect to $\Delta t$ at time $T=4$ for the (R-CN) scheme (figures on left)
and the (C-CN)-scheme (figures on right) and for $J=20 000$.}
\label{example1_error_dt}
\end{figure}

The behaviour of these two errors with respect to $\Delta x$ are
presented in Figure~\ref{example1_error_dx}. 
We observe that we obtained numerically the expected order of accuracy
for each scheme: the (R-CN) scheme has a convergence order of one 
and the (C-CN) scheme is of order two. 
We can also see that for each value of $N$ there is a saturation phenomena for the error,
for very small $\Delta x$ the round-off errors balances with the errors in the solution.
Also, changing $\Delta x$ also modifies the roots of the cubic/quartic equations 
needed in the boundary convolution and the numerical inverse $\mathcal{Z}$-transform of the 
convolution kernels may
degrade the overall accuracy
(at least for the selected inversion radius).

We present in Figure~\ref{example1_error_dt} the $rel.ErrTm$ and 
$rel.ErrL2$ with respect to $\Delta t$ for $J=20 000$ and $r=1.001$. 
Again we obtain for each scheme a numerical rate of convergence of order two in 
time. Surprisingly, the saturation effects from the previous figure do not
show up, although with smaller $\Delta t$ the size of the boundary convolutions 
is increasing and this often leads to additional errors. 


\begin{figure}[ht]
\begin{center}
\begin{tabular}{cc}
\includegraphics[width=0.48\textwidth]{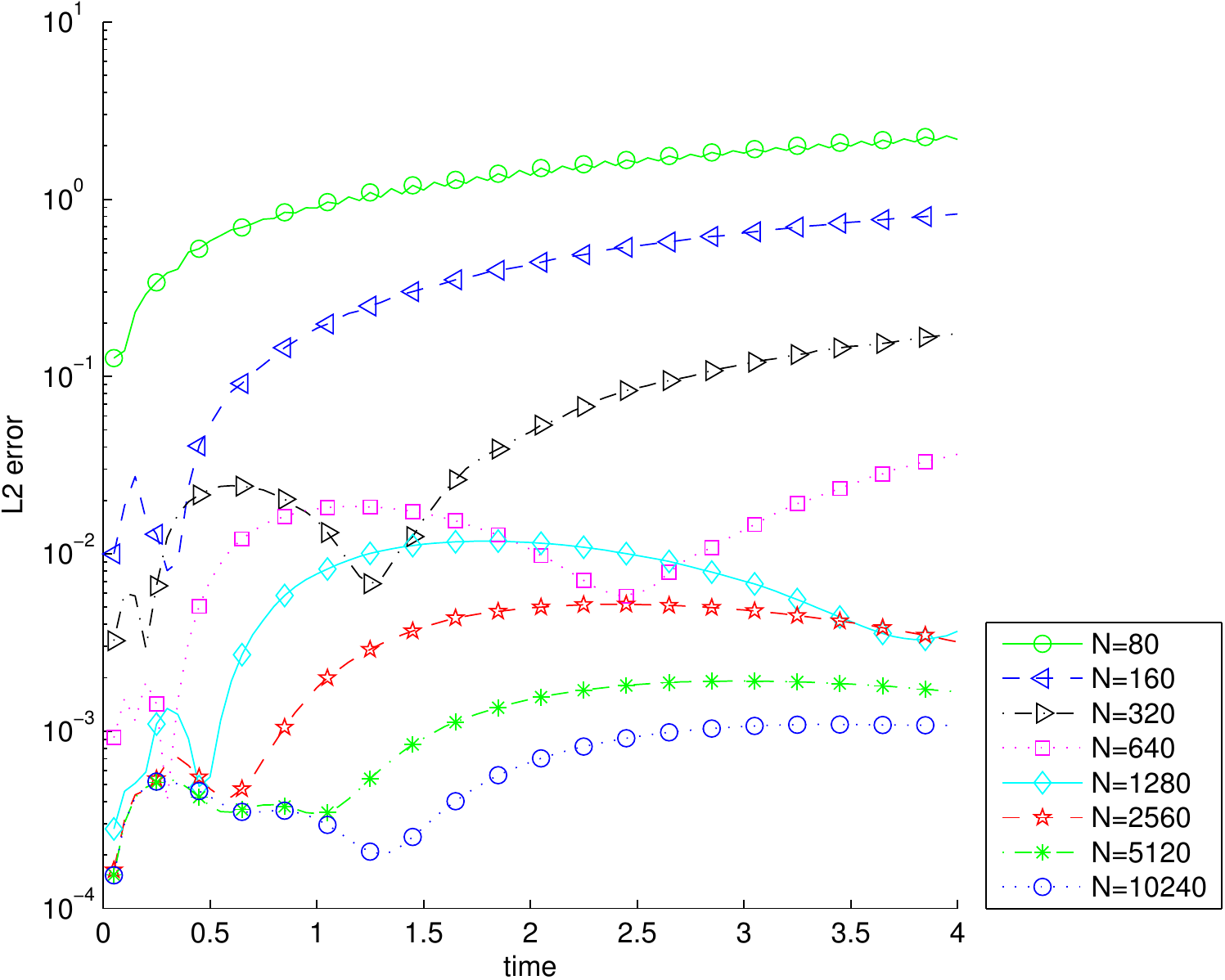} 
& \includegraphics[width=0.48\textwidth]{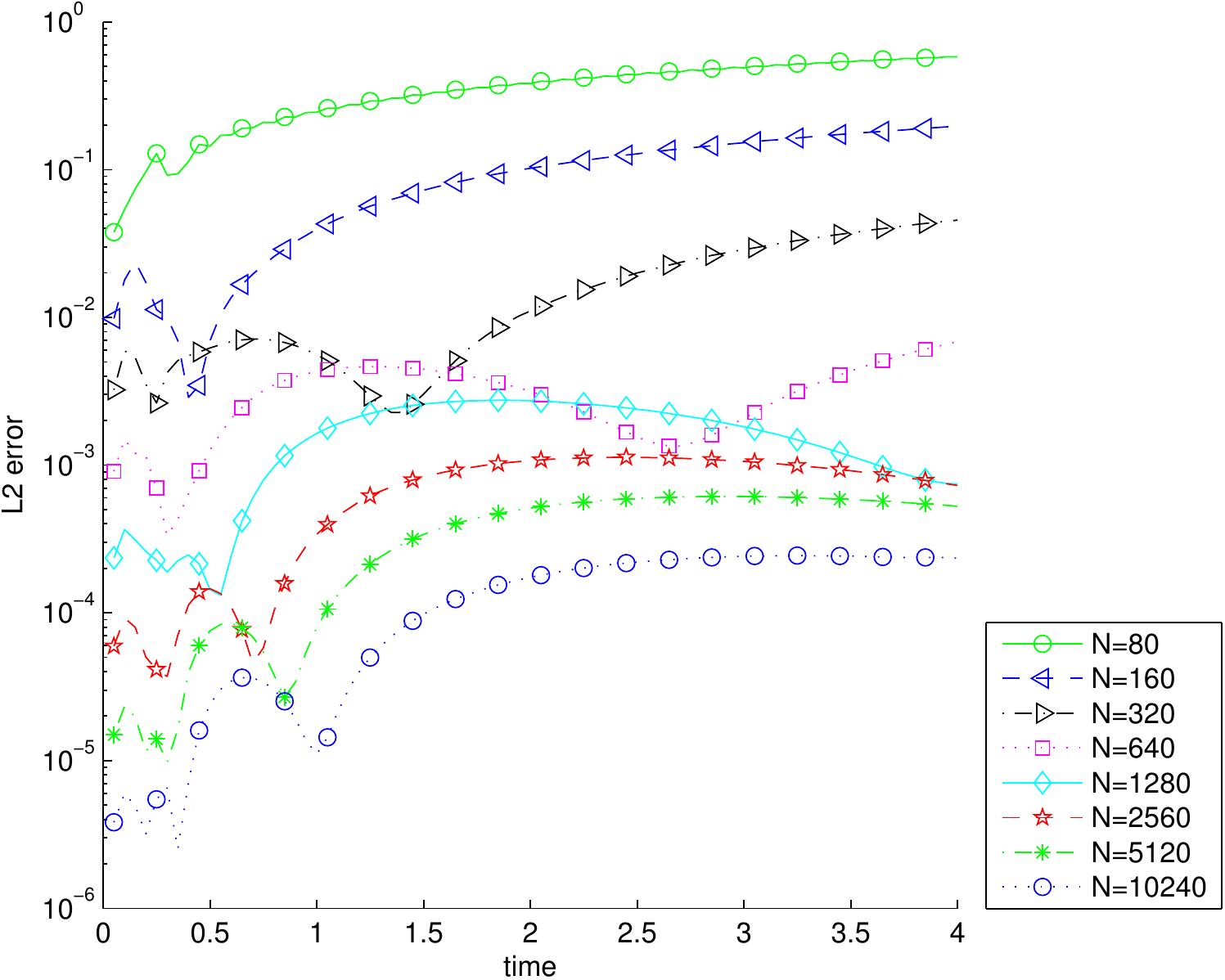} 
\end{tabular}
\end{center}
\caption{ Evolution of the $\ell^2$ error between $T=0$ and $T=4$ for the (R-CN) scheme 
(figure on left) and the (C-CN)-scheme (figure on right) and for $J=20000$ and various values of $N$.}
\label{example1_error_evol}
\end{figure}

We present in Figure~\ref{example1_error_evol} the evolution of the $\ell^2$-error with respect
to time for various values of $N$, $J=20 000$ and $r=1.001$. 
As expected, the error decreases for increasing $N$, i.e.\ finer mesh size. In any case,
the error remains moderately bounded over the whole simulation time which shows the
usefulness of the proposed method.
At the beginning the first increase is due to the interaction with the artificial boundaries
and the second long term growth is due to an accumulation effect of errors, e.g.\
due to the increasing time convolution at the boundaries.

\begin{figure}[ht]
\begin{center}
\begin{tabular}{cc}
\includegraphics[width=0.48\textwidth]{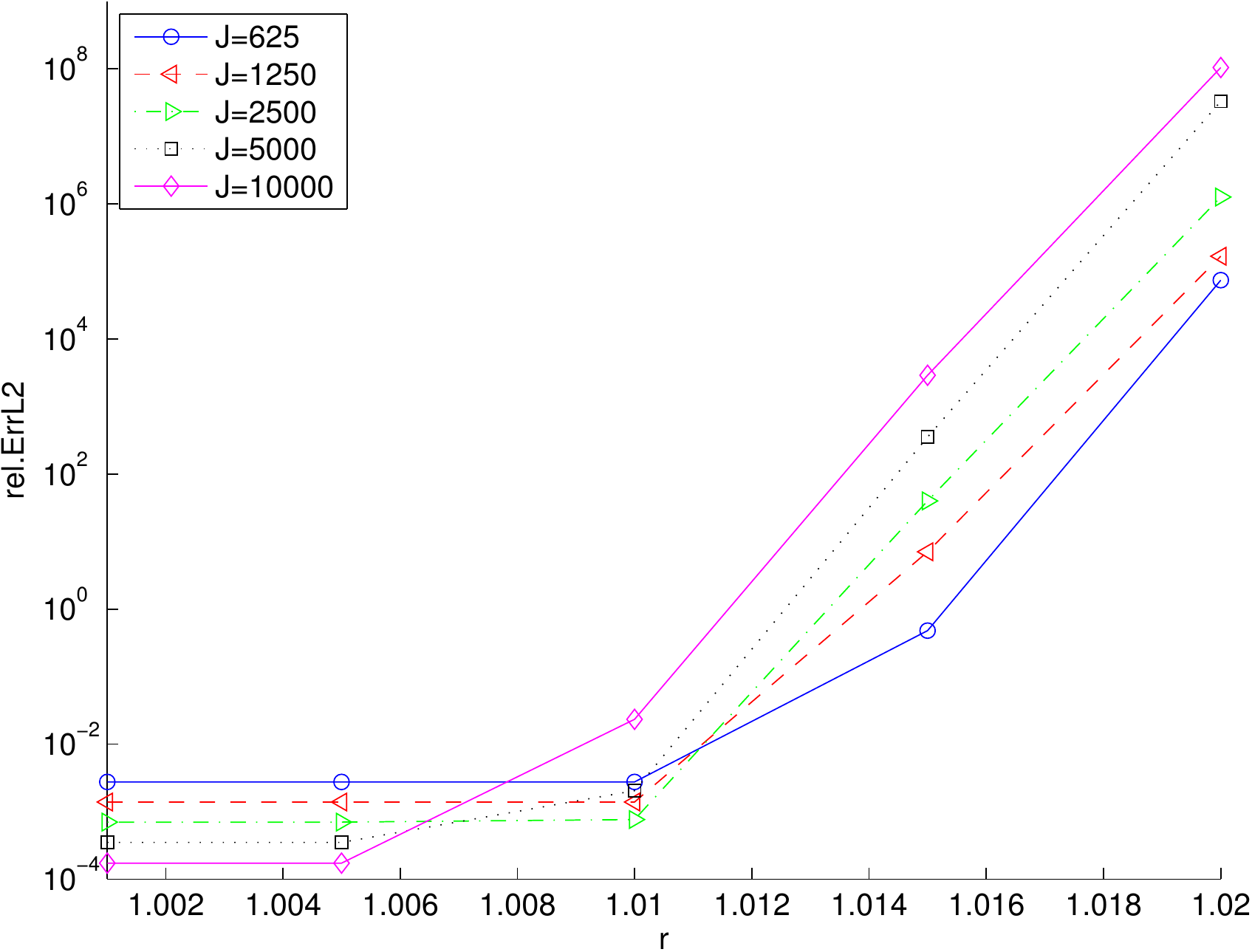} & 
\includegraphics[width=0.48\textwidth]{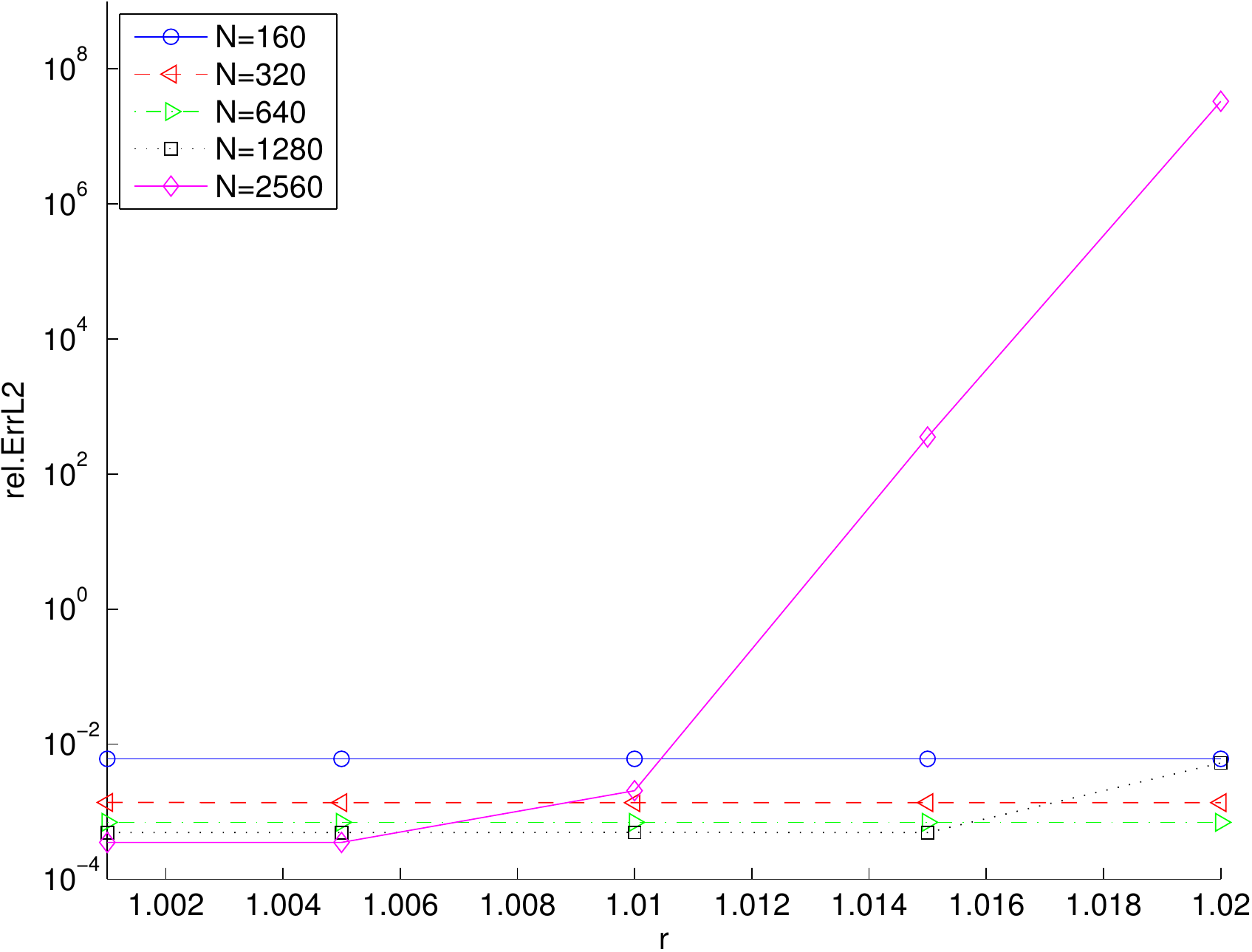}\\
\includegraphics[width=0.48\textwidth]{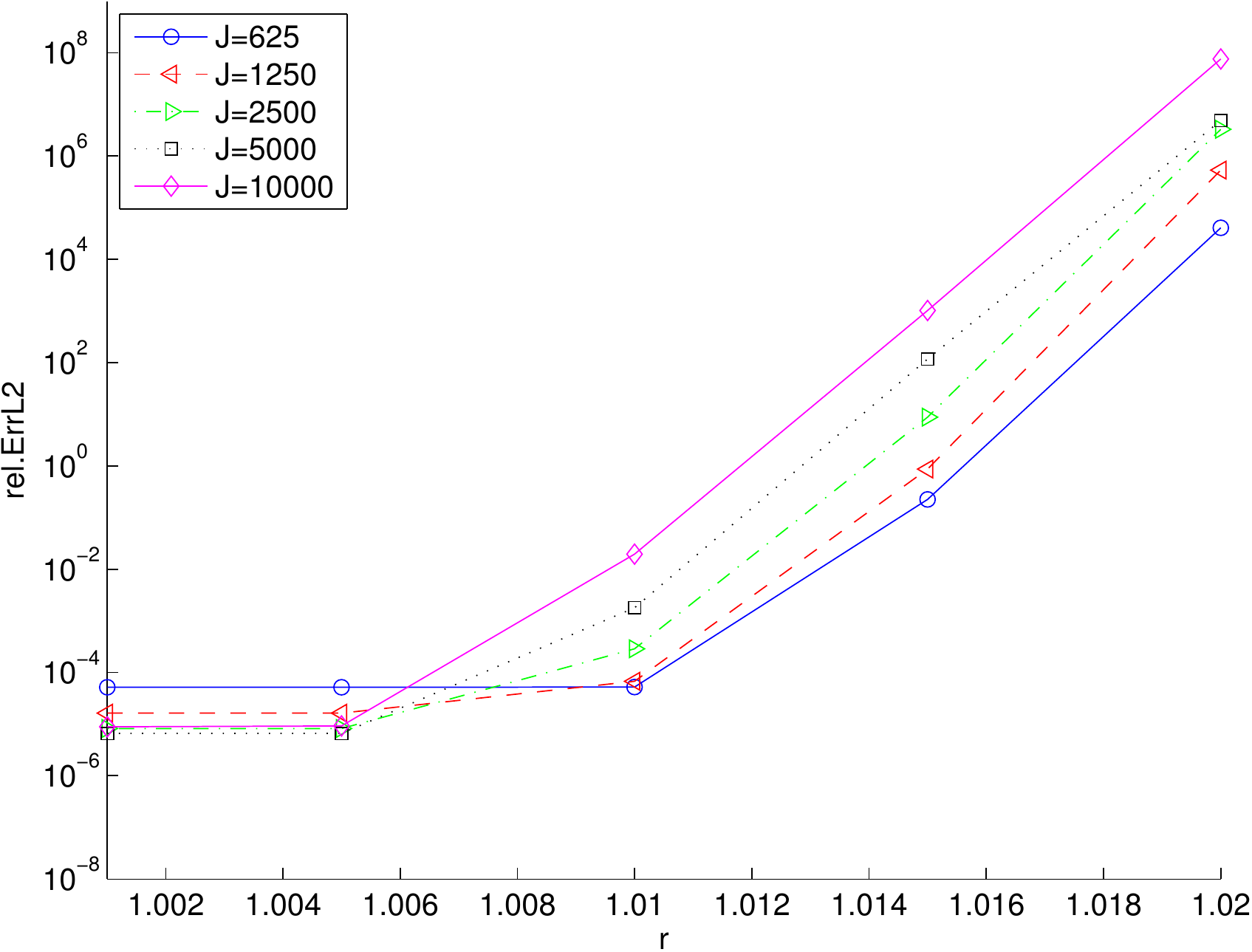} & 
\includegraphics[width=0.48\textwidth]{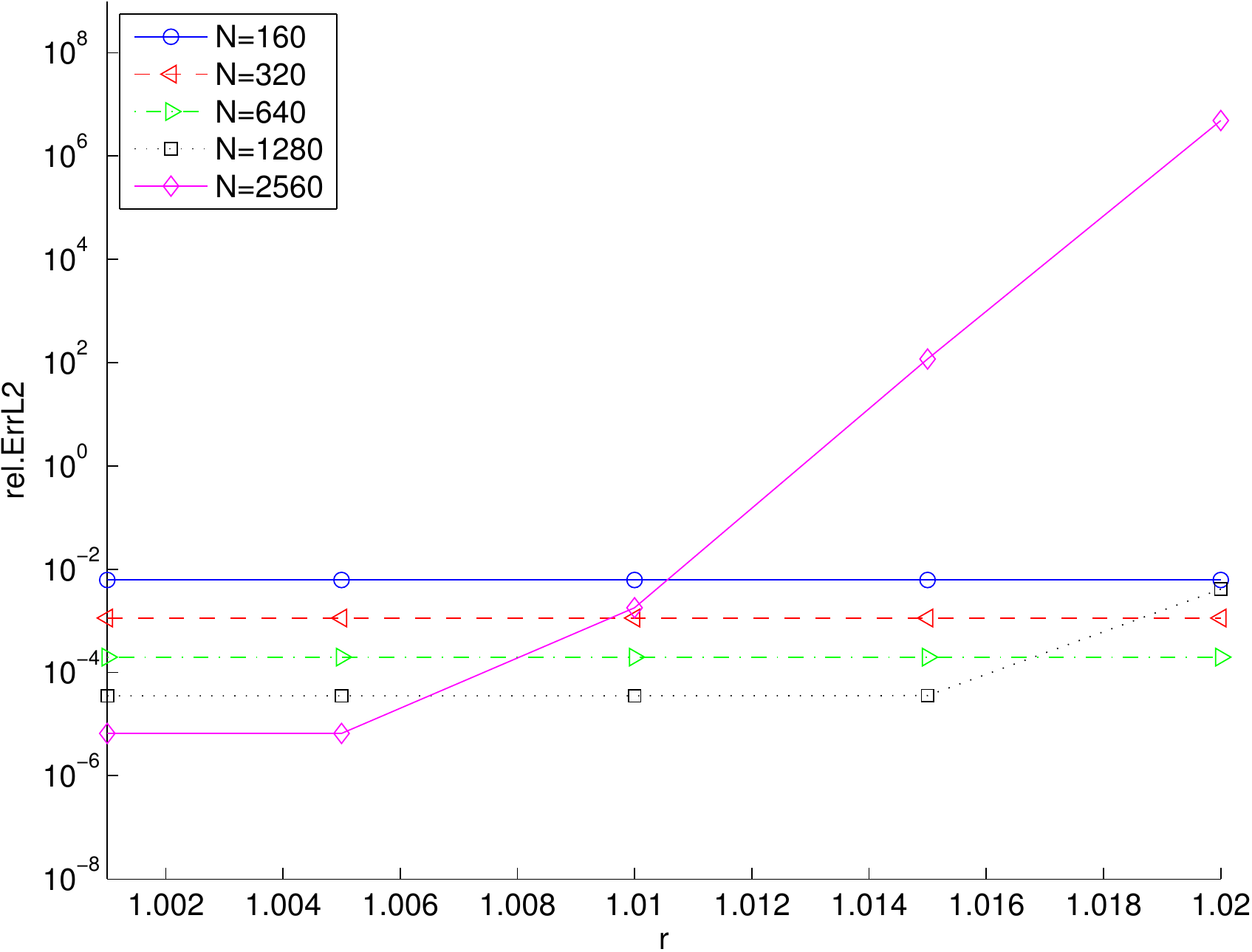}
\end{tabular}
\end{center}
\caption{ Error with respect to $r$ at time $T=4$ for the (R-CN) scheme (top figures) and
the (C-CN)-scheme (bottom figures) with either $N=2560$ and various $J$ (figures on the left) 
or $J=5000$ and various $N$ (figures on the right).}
\label{error_wrt_r}
\end{figure}

We present in Figure~\ref{error_wrt_r} the $rel.ErrL2$ with respect to $r$ 
for each scheme and either with $N=2560$ and various $J$ or $J=5000$ and various $N$. 
The choice of $r$ in the inverse $\mathcal{Z}$-transform procedure is clearly impacting
the error and depend on the values of $J$ and $N$. 
It seems that our choice, $r=1.001$, 
is a good choice for a large set of values of $N$ and $J$.

\subsection{Numerical Example 2}\label{num2}
\begin{figure}[ht]
\begin{center}
\begin{tabular}{cc}
\includegraphics[width=0.48\textwidth]{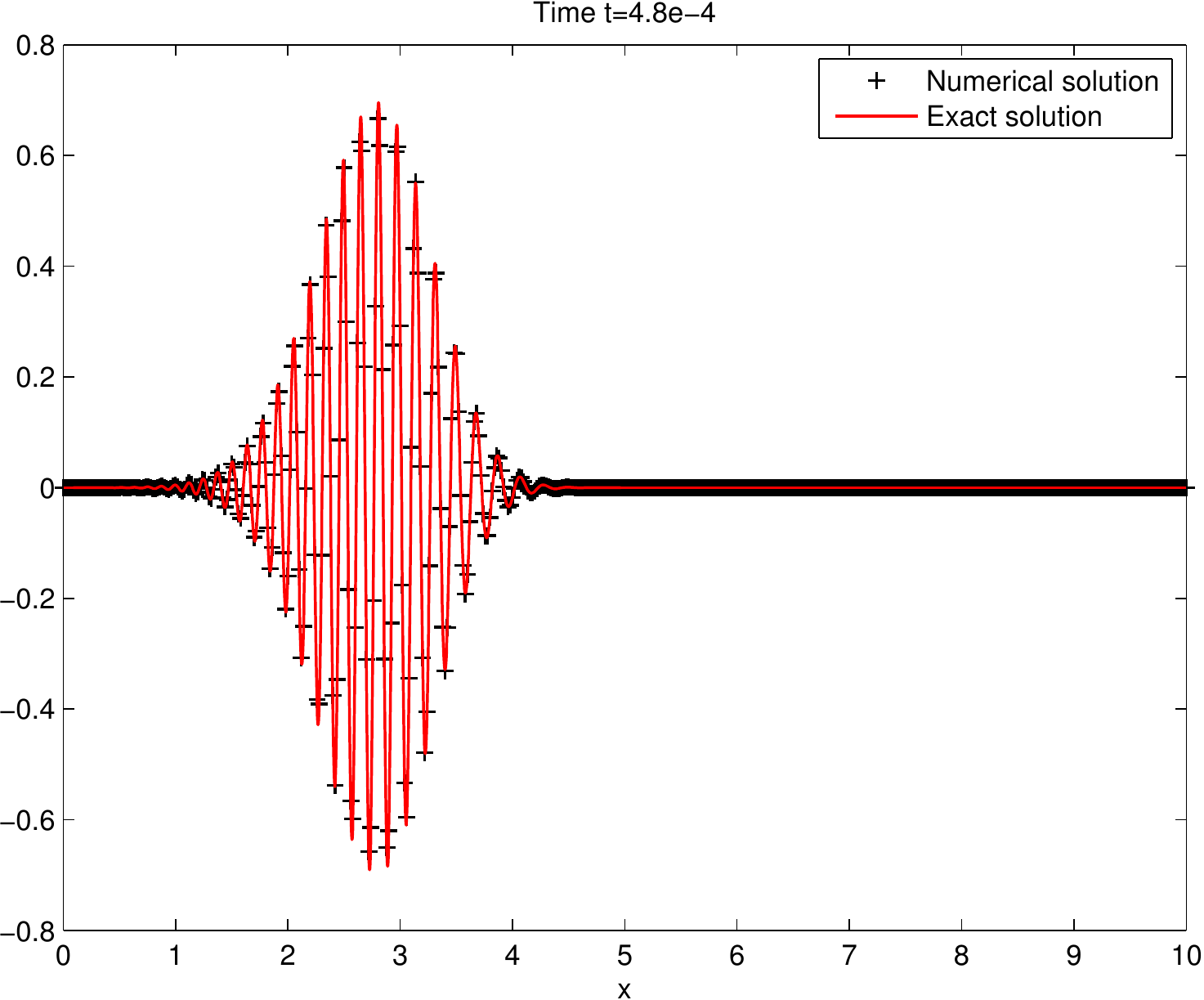} 
\end{tabular}
\end{center}
\caption{Numerical and exact solutions at final time for the second example with $\Delta t=T/2560$, 
$\Delta x=10/5000$ and $r=1.001$.}
\label{compBriggsCCN}
\end{figure}

\begin{figure}[ht]
\begin{center}
\begin{tabular}{cc}
\includegraphics[width=0.48\textwidth]{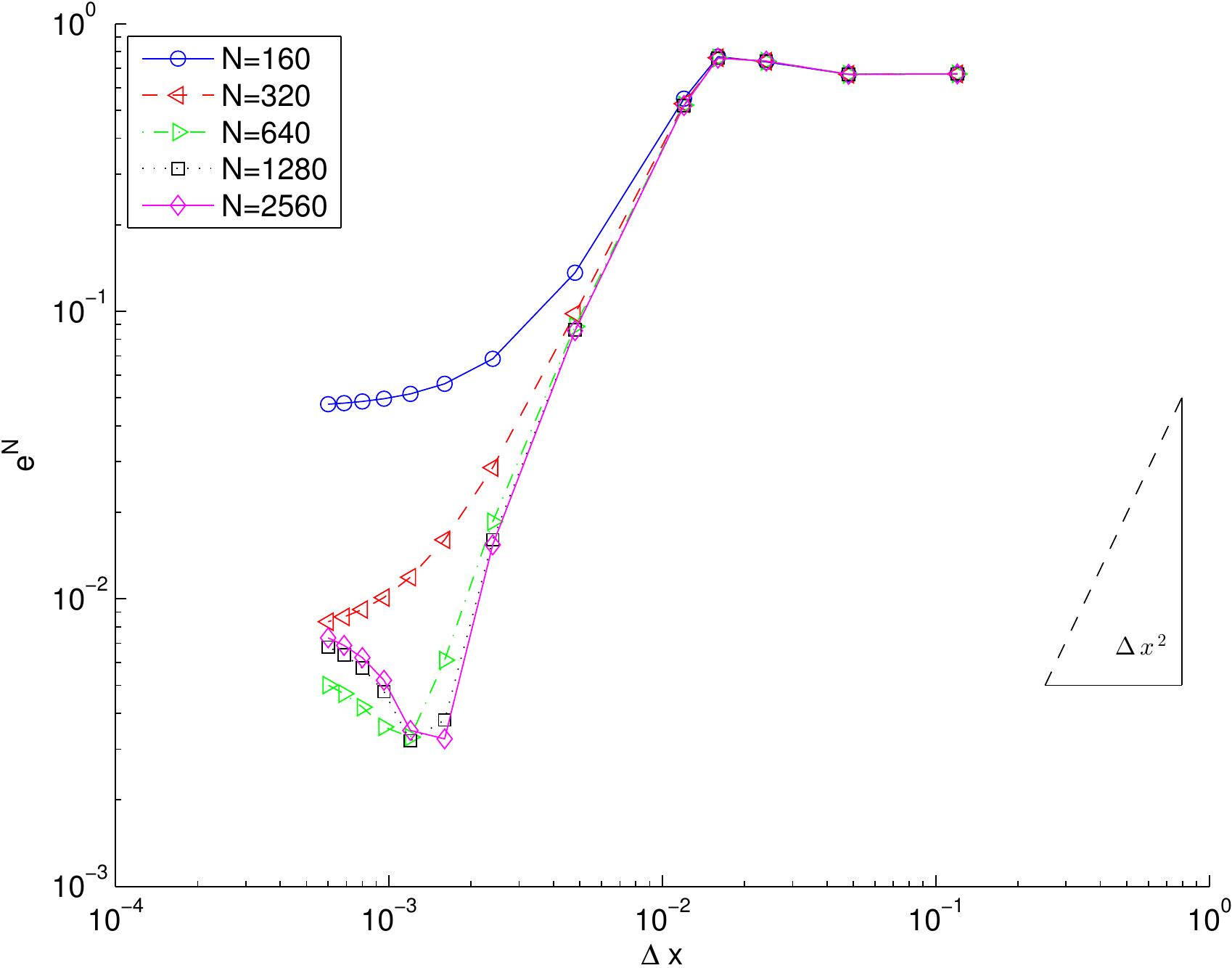} 
\end{tabular}
\end{center}
\caption{ Relative error $e^{(N)}$ with respect to $\Delta x$ for various values of $N$ using 
(C-CN) scheme for the second example}
\label{error_wrt_dx_briggs}
\end{figure}

Let us now consider a second example. 
We consider the dispersive equation \eqref{eqKdVL} with $U_1=U_2=1$ and we choose as initial condition 
\begin{equation*}
u_0(x)=\exp(-8(x-5)^2)\sin\left(\frac{50\pi}{4}\right), 
\end{equation*}
for $0 \le x \le 10$ and for a final time $T=4.8 \times10^{-4}$. 
This example was already considered in \cite{BrSa83}.
 Note that using the Fourier transform, the problem being a linear and periodic problem,
 we can compute the exact solution ${u}_{\rm exact}(t,x)$. 
 Indeed, applying the Fourier transformation in the 
 space variable to the equation \eqref{eqKdVL} we obtain
\begin{equation*}
\widehat{u}_{t}+i\xi\widehat{u}-i\xi^3\widehat{u}=0, 
\end{equation*}
where $\xi$ stands for the Fourier variable. 
Then it is easy to see that the transformed exact solution reads
\begin{equation*}
\widehat{u}_{\rm exact}(t,\xi)
=\widehat{u}_0\exp\left(-\left(i\xi-i\xi^3\right)t\right).
\end{equation*}
Using the inverse Fourier transform we have the exact solution of the problem. 
A reference solution is computed using 50000 points in space and 2560 iterations 
in time and used for Figures~\ref{compBriggsCCN} and \ref{error_wrt_dx_briggs}.

We present in Figure~\ref{compBriggsCCN} the exact solution and the approximate 
solution obtained with (C-CN) scheme for $\Delta t=T/2560$, $\Delta x=10/5000$ 
and $r=1.001$ at final time. We see that the (C-CN) solution is a very good approximation of 
the exact oscillatory solution, the two solutions are nearly indistinguishable. 

We present in Figure~\ref{error_wrt_dx_briggs} the relative $\ell^{2}$-error $e^{(n)}$
computing at final time ({\it i.e.} for n=N) with respect to $\Delta x$ and for various 
values of $N$ and $r=1.001$. 
Again, we see that we obtain the order two as predicted. 
As for the first example there is a saturation phenomena. 

\clearpage 

\section*{Conclusion and Outlook}
In this work we presented some new 
discrete absorbing boundary conditions adapted to two different numerical 
schemes for the linearized KdV equation \eqref{eqKdVL}. 
The orders of each scheme in time and space are shown numerically and given evidence 
that they are not 
perturbed by the discrete absorbing boundary conditions. 
To speed up the calculations of the costly boundary convolutions, 
especially in higher-dimensional cases,
we proposed to use the sum-of-exponentials ansatz.
We gave finally two numerical examples that supported our theoretical findings.

Future work will be to design an automatic good choice of the inversion radius,
establish a transformation rule in the spirit of \cite{AES03} for the KdV equation,
treat the 2D and the nonlinear case.

\section*{Acknowledgements}
This work was partially supported by the French ANR grant MicroWave NT09 460489 (``Programme
Blanc`` call) and Universit\'e Paul Sabatier Toulouse 3. The first author also acknowledges support from the French ANR grant BonD ANR-13-BS01-0009-01. The third author acknowledges support from the team INRIA/RAPSODI and the Labex CEMPI (ANR-11-LABX-0007-01).

\begin{appendices}
  \section{Proof of theorem \ref{theo:continuous} \label{proof:theo:cont}}
\begin{figure}[htbp]
\begin{center}
\begin{tabular}{ll}
\includegraphics[width=0.5\textwidth]{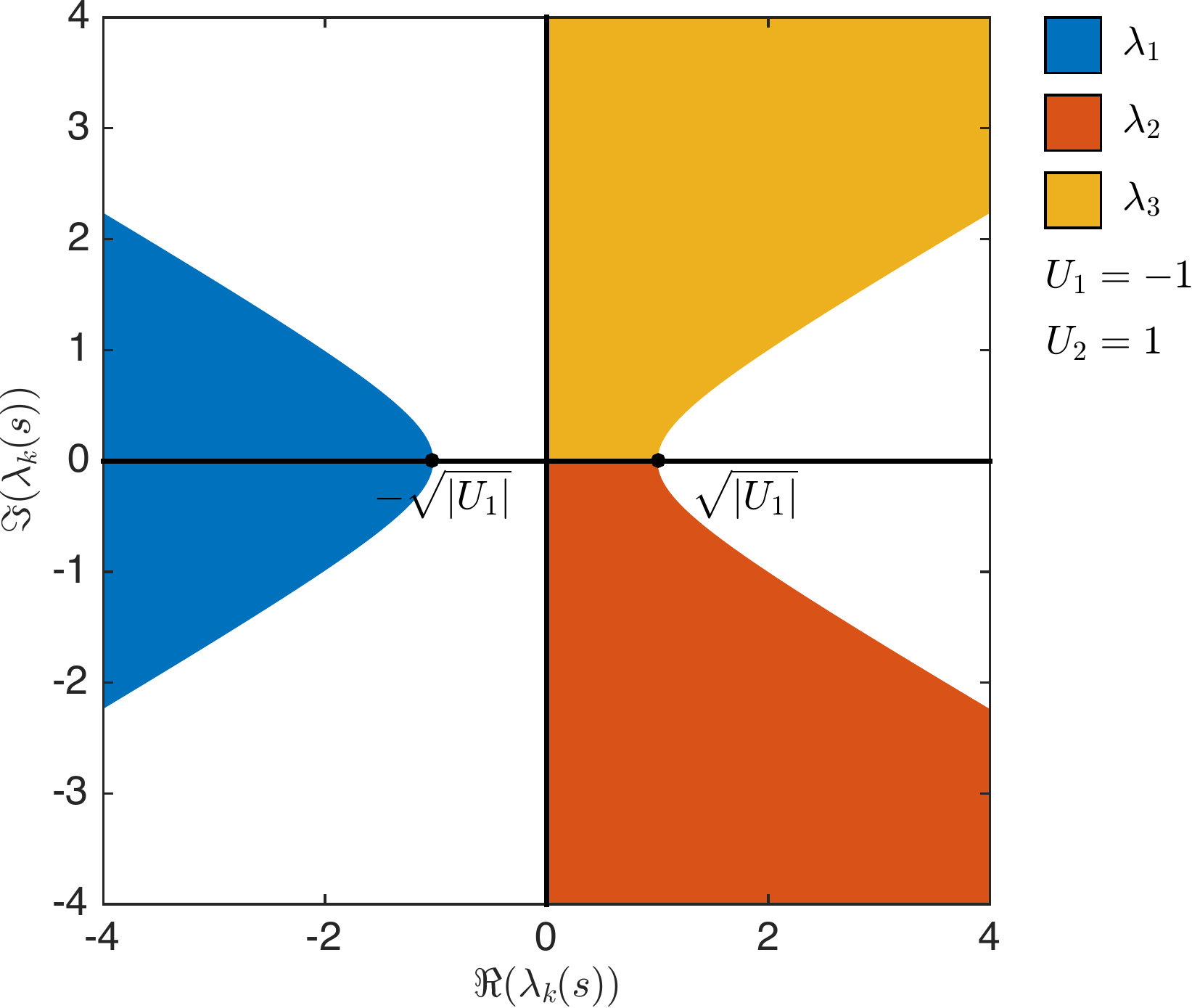}
&
\includegraphics[width=0.5\textwidth]{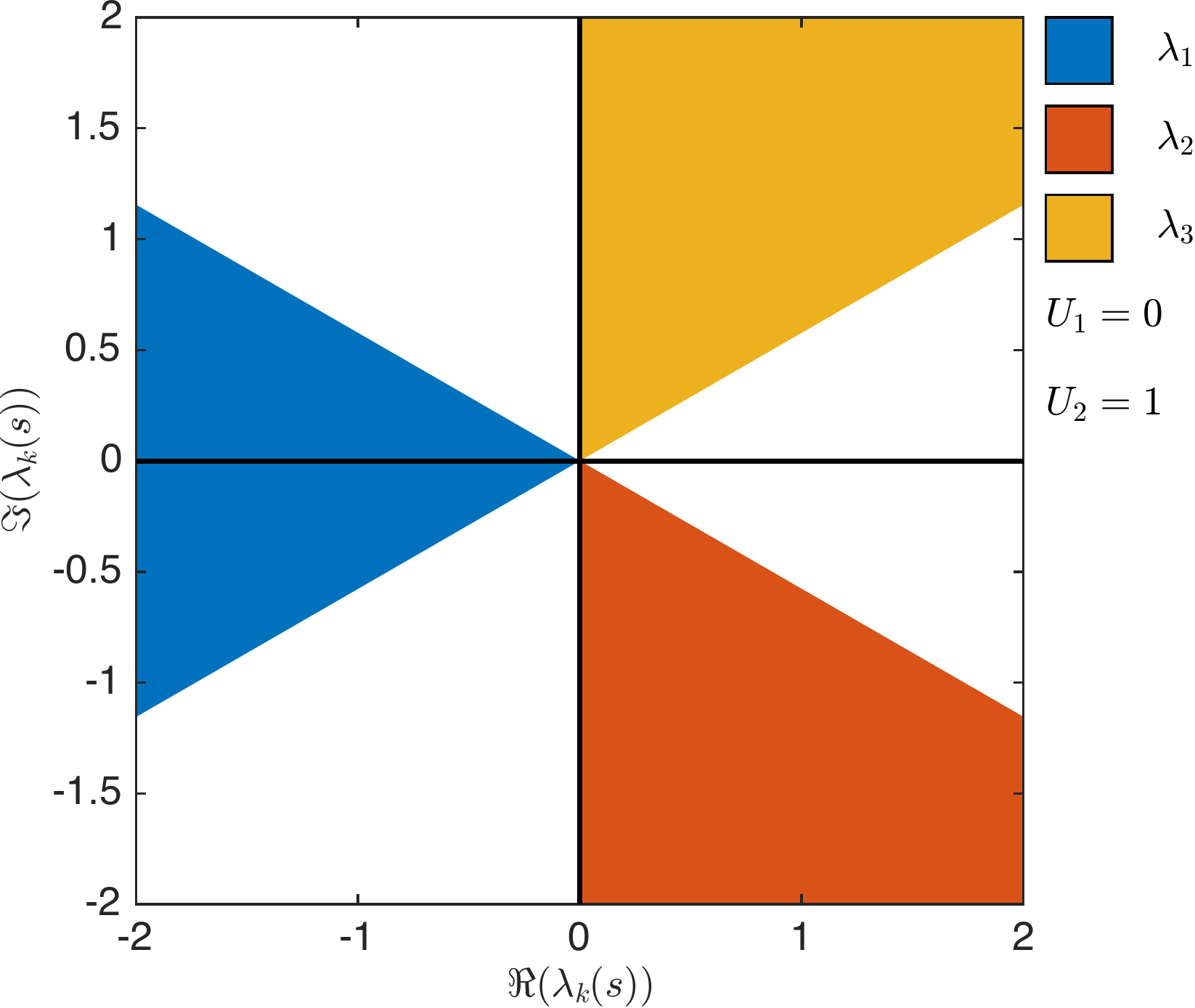}
\end{tabular}
\centerline{\includegraphics[width=0.5\textwidth]{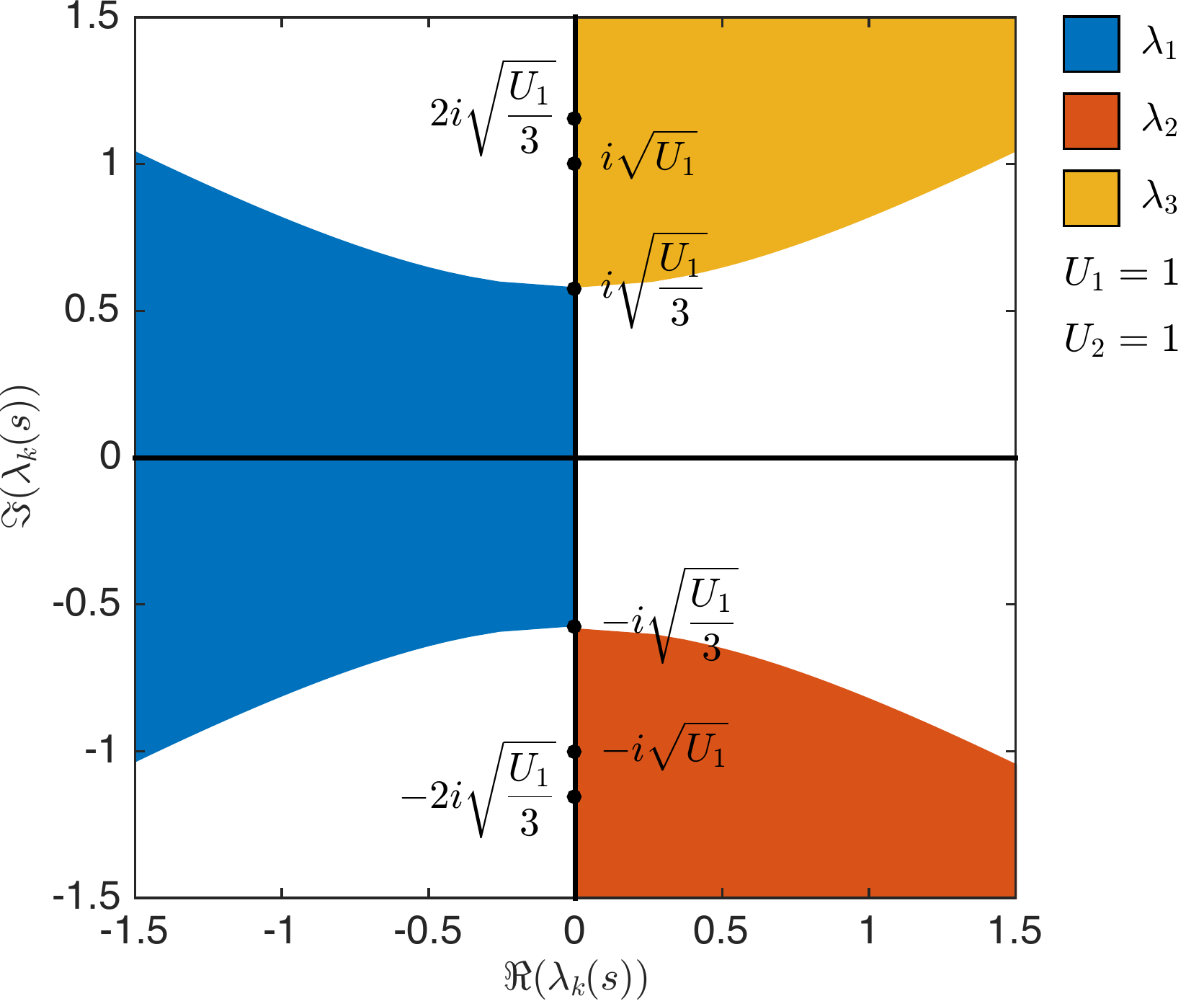}}
\end{center}
\caption{Roots $\lambda_1(s)$, $\lambda_2(s)$, 
$\lambda_3(s)$ to the equation~\eqref{cubic_cont} 
for $U_1=-1,\ 0,\ 1$ and $U_2=1$. }
\label{fig:cubiceqcont}
\end{figure}

We saw in Remark \ref{rmk:U2} that the dispersive constant $U_2$ can be forced to be equal to $1$. Without loss of generality, we therefore consider for this proof $U_2=1$. The behavior of the roots is given on Figure \ref{fig:cubiceqcont}.
We prove here the result for stricly positive or negative velocity $U_1$. In both cases, we want to determine the domains $\mathbb{D}_k$, image of $\mathbb{C}^+=\{z\in\mathbb{C},\ \text{Re } z>0\}$, by roots $\lambda_k$ defined by \eqref{eq:roots}
$$
\begin{array}{ccl}
  \C^+ & \longrightarrow & \mathbb{D}_k\\
  s & \longmapsto & \lambda_k(s).
\end{array}
$$
Since $\C^+$ is simply connected, we just have to identify the boundary of each domain $\mathbb{D}_k$ and a single point inside it to completely determined them. The boundaries are given by the image of $s=\varepsilon+i\xi$, $\varepsilon>0$, $\varepsilon\ll 1$ and $\xi \in \R$.
\subsection{Case $U_1>0$}
We perform an asymptotic expansion with respect to $\varepsilon$. Let us define $A=\sqrt{|\xi^2-(4/27)U_1^3|}$, $B=((\xi+A)/2)^{1/3}$, $C=B-U_1/(3B)$ and $D=B+U_1/(3B)$. Let us consider $\xi>0$.

For $\xi^2>4U_1^3/27$, the  general expression for roots is for $k=1,2,3$
$$
\lambda_k=-\left (\omega^{k-1}Be^{i\pi/6}-\frac{U_1}{3}\frac{1}{\omega^{k-1}Be^{i\pi/6}}\right) + \left (\omega^{k-1}Be^{2i\pi/3}-\frac{U_1}{3}\frac{1}{\omega^{k-1}Be^{2i\pi/3}}\right) \frac{\varepsilon}{3A} + O(\varepsilon^2),
$$
and we have
$$
\begin{array}{l}
 \displaystyle \lambda_1(s)=-\left (\frac{\sqrt{3}}{2} C+\frac{i}{2} D\right ) - \left (\frac{1}{2}C-i\frac{\sqrt{3}}{2} D\right ) \frac{\varepsilon}{3A}+O(\varepsilon^2), \\[2mm]
 \displaystyle \lambda_2(s)=-\left (-\frac{\sqrt{3}}{2} C+\frac{i}{2} D\right ) - \left (\frac{1}{2}C+i\frac{\sqrt{3}}{2} D\right ) \frac{\varepsilon}{3A}+O(\varepsilon^2), \\[2mm]
 \displaystyle \lambda_3(s)=iD+C\frac{\varepsilon}{3A}+O(\varepsilon^2).
\end{array}
$$
It is easy to show that $A$, $B$, $C$ and $D$ are positive for any $\xi^2>4U_1^3/27$. Clearly, we have $\Real{\lambda_1}<0$, $\Real{\lambda_2}>0$ and $\Real{\lambda_3}>0$. Concerning $\lambda_1$, since $-C/2<0$ and $D/2>0$, domain $\mathbb{D}_1$ is located on the left and above the complex curve $-(\sqrt{3}C+iD)/2$. The conclusion for $\lambda_2$ is similar and we conclude that domain $\mathbb{D}_2$ is located on the left and below the complex curve $(\sqrt{3}C-iD)/2$. The domain $\mathbb{D}_3$ is located on the right of the complex curve $iD$.

For $\xi^2<4U_1^3/27$, we define $E=((i\xi+A)/2)^{1/3}$ and have $|E|=\sqrt{U_1/3}$. We therefore can write $E=\sqrt{U_1/3}e^{i\theta(\xi)}$ with $0<\theta(\xi)<\pi/6$. We obtain
$$
\lambda_k(s)=-2\sqrt{\frac{U_1}{3}}i \sin{\left(\theta+(k-1)\frac{2\pi}{3}\right)}-2\sqrt{\frac{U_1}{3}}\cos{\left(\theta+(k-1)\frac{2\pi}{3}\right)}\frac{\varepsilon}{3A}+O(\varepsilon^2).
$$
We can conclude that $\lambda_1$ is located on the left of the complex curve $if(\theta)$, $-\sqrt{U_1/3} < f(\theta)<0$, $\lambda_2$ is on the right of the complex curve $ig(\theta)$, $-\sqrt{U_1}<g(\theta)<-\sqrt{U_1/3}$ and $\lambda_3$ is on the right of $ih(\theta)$, $\sqrt{U_1}<h(\theta)<2\sqrt{U_1/3}$.

The situation is completely symetric if we consider $\xi<0$. We therefore recover the figure \ref{fig:cubiceqcont} and we have shown the {\em separation property}
\begin{equation*}
\Real(\lambda_1(s))<0,\qquad\Real(\lambda_2(s))>0,\qquad\Real(\lambda_3(s))>0.
\end{equation*}
\subsection{Case $U_1<0$}
Like in the previous case, we perform an asymptotic expansion and consider $\xi>0$. Let us define $\tilde{A}=\sqrt{\xi^2+4|U_1|^3/27}$, $\tilde{B}=((\xi+\tilde{A})/2)^{1/3}$, $\tilde{C}=\tilde{B}-U_1/(3\tilde{B})$ and $\tilde{D}=\tilde{B}+U_1/(3\tilde{B})$. For any $\xi \in \R$, we have $\tilde{A}>0$, $\tilde{B}>0$, $\tilde{C}>0$ and $\tilde{D}>0$. The general expression for roots is for $k=1,2,3$
$$
\lambda_k=-\left (\omega^{k-1}\tilde{B}e^{i\pi/6}-\frac{U_1}{3}\frac{1}{\omega^{k-1}\tilde{B}e^{i\pi/6}}\right) + \left (\omega^{k-1}\tilde{B}e^{2i\pi/3}-\frac{U_1}{3}\frac{1}{\omega^{k-1}\tilde{B}e^{2i\pi/3}}\right) \frac{\varepsilon}{3A} + O(\varepsilon^2).
$$
Thus, the asymptotics of $\lambda_k$ are given by
$$
\begin{array}{l}
 \displaystyle \lambda_1(s)=-\left (\frac{\sqrt{3}}{2}\tilde{C}+\frac{i}{2} \tilde{D}\right ) - \left (\frac{1}{2}\tilde{C}-i\frac{\sqrt{3}}{2} \tilde{D}\right ) \frac{\varepsilon}{3A}+O(\varepsilon^2), \\[2mm]
 \displaystyle \lambda_2(s)=-\left (-\frac{\sqrt{3}}{2} \tilde{C}+\frac{i}{2} \tilde{D}\right ) - \left (\frac{1}{2}\tilde{C}+i\frac{\sqrt{3}}{2} \tilde{D}\right ) \frac{\varepsilon}{3A}+O(\varepsilon^2), \\[2mm]
 \displaystyle \lambda_3(s)=i\tilde{D}+\tilde{C}\frac{\varepsilon}{3A}+O(\varepsilon^2).
\end{array}
$$
Thanks to the signs of $\tilde{A}$, $\tilde{B}$, $\tilde{C}$ and $\tilde{D}$, the conclusions can be drawn by a simliar study to the $U_1>0$ case and we have $\Real(\lambda_1(s))<0$, $\Real(\lambda_2(s))>0$ and $\Real(\lambda_3(s))>0$.
  \section{Proof of theorem \ref{theo:discr3} \label{proof:theo:discr3}}
We consider here the continuous roots $\ell_k$ of \eqref{fulldisc6} and ordered them thanks to the relation $|\ell_1(z)| \leq |\ell_2(z)| \leq |\ell_3(z)|$. We defined $p=\mu(z-1)/(z+1)$. If we note $z=x+iy$, $(x,y)\in \R^2$, then 
$$
p=\lambda \frac{(x^2+y^2-1) + 2iy}{(x+1)^2+y^2}.
$$
Since $|z|>1$, we have $x^2+y^2-1>0$ and
$$\Real{p}>0.$$
Thus, instead of studying $\ell_k$ as functions of $z$, we now consider them as functions of variable $p$ and try to identify the domains $\mathbb{D}_k$ 
$$
\begin{array}{llcl}
  \ell_k:& \C^+ &\longrightarrow &\mathbb{D}_k\\
& p &\longmapsto &\ell_k(p)
\end{array}
$$
Let us first show that $\ell_k$ can not belong to the unit circle. Let us assume that there exists $\theta \in [0,2\pi)$ such that $\ell_2=e^{i\theta}$. Since $\ell_k$ is a root of \eqref{fulldisc6}, then
$$
e^{3i\theta}-3 e^{2i\theta}+(3+p)e^{i\theta}-1=0.
$$
It leads to $p=e^{i\theta}-3+e^{-i\theta}-e^{2i\theta}$. We therefore obtain
$$
\Real{p}=-2(\cos\theta -1)^2\leq 0.
$$
We have a contradiction since we had assume that $\Real{p}>0$. Thus, $\ell_k$ can not belong to the unit circle.
Since we consider the continuous roots $\ell_k$ and $\C^+$ is simply connected, then $\mathbb{D}_k\in B(0,1)$ or $\mathbb{D}_k \in \bar{B}(0,1)$, the complementary of $B(0,1)$ in $\C$. In order to find if $\mathbb{D}_k$ lie inside or outside the unit circle, we just have to determine the domains for a single value of $p$.

Next, we know that if we consider a third order algebraic equation
$$
x^3+bx^2+cx+d=0
$$
then its roots satisfies
$$
\begin{array}{l}
  x_1+x_2+x_3=-b\\
  x_1x_2+x_2x_3+x_3x_1=c\\
  x_1x_2x_3=-d.
\end{array}
$$
Therefore, $\ell_1$, $\ell_2$ and $\ell_3$ satisfy
$$
\begin{array}{l}
  \ell_1\ell_2\ell_3=1\\
\ell_1+\ell_2+\ell_3=3\\
\ell_1\ell_2+\ell_2\ell_3+\ell_3\ell_1=3+p.
\end{array}
$$
The first relation implies $|\ell_1|\leq 1\leq |\ell_3|$. We therefore automatically have one root which lies inside the unit circle and one outside. It remains to understand the behavior of $\ell_2(p)$. In order to see the location of $\ell_2$ in the complex plane, we can therefore take any value of $p$ such that $\Real{p}>0$. If we take $p=3$, we have 
$$
|\ell_1|\approx 0.18, \quad |\ell_2| \approx 2.34, \quad |\ell_3|\approx 2.34.
$$
Thus, $|\ell_2(p)|>1$ for any $p\in \mathbb{C}^+$ and we have shown
$$
|\ell_1(p)|<1, \quad |\ell_2(p)|>1, \quad |\ell_3(p)|>1, \quad p \in \C^+.
$$
  \section{Proof of theorem \ref{theo:discr4} \label{proof:theo:discr4}}
The proof is similar to the previous one. Let us consider a root $\ell_k$ and show that it can not be equal to $e^{i\theta}$. Let us assume that $\ell_k=e^{i\theta}$. Then
$$
e^{4i\theta}-(2-a)e^{3i\theta}+pe^{2i\theta}+(2-a)e^{i\theta}-1=0.
$$
Since $a \in \R$, this leads to
$$
  p = 2i((2-a)\sin \theta -\sin 2\theta) \in i\R.
$$
But, $\Real{p}>0$. This is a contradiction and $\ell_k\neq e^{i\theta}$.\\
We moreover have
$$
\begin{array}{l}
  \ell_1+\ell_2+\ell_3+\ell_4=2-a,\\
\ell_1\ell_2 + \ell_1\ell_3+\ell_1\ell_4+\ell_2\ell_3+\ell_2\ell_4+\ell_3\ell_4=p,\\
\ell_1\ell_2\ell_3+\ell_1\ell_2\ell_4+\ell_2\ell_3\ell_4+\ell_1\ell_3\ell_4=a-2,\\
\ell_1\ell_2\ell_3\ell_4=-1.
\end{array}
$$
If we sort the roots as $|\ell_1|\leq |\ell_2| \leq |\ell_3| \leq |\ell_4|$, the last equation leads to
$$
|\ell_1(p)|<1 \quad \text{and} \quad |\ell_4(p)|>1.
$$
Computing $\ell_2$ and $\ell_3$ for $p=1$ for example prove that $|\ell_2(p)|<1$ and $|\ell_3(p)|>1$. We therefore have the discrete separation property of the roots.

\end{appendices}



\end{document}